\documentclass[11pt,reqno]{amsart}

\usepackage{xcolor}
\usepackage[T1]{fontenc}
\usepackage{amsmath}							
\usepackage{amssymb}
\let\amslrcorner\lrcorner   %save the lrcorner command from amssymb
\usepackage{amsthm}
\usepackage{amscd}
\usepackage{amsfonts}
\usepackage{mathabx}
\usepackage{mathtools}
\usepackage{stmaryrd}
\usepackage[all]{xy}
\usepackage{array}
\usepackage{multirow}       %allows to join rows in a table to form bigger cells
\usepackage{comment}
\usepackage{subfiles}       %allows to precompile parts of the document
\usepackage[english]{babel}
\usepackage[autostyle]{csquotes}

\usepackage{caption}
\usepackage{subcaption}
\usepackage{placeins}
\usepackage{euler}

\usepackage{extarrows}

\usepackage[colorlinks, linktocpage, citecolor = purple, linkcolor = blue]{hyperref}
\usepackage{color}

\usepackage{tikz}									
\usetikzlibrary{matrix}
\usetikzlibrary{patterns}
\usetikzlibrary{positioning}
\usetikzlibrary{decorations.pathmorphing}
\usetikzlibrary{decorations.markings}
\usetikzlibrary{cd}

\usepackage{ytableau}
\usepackage{fullpage}
\usepackage[shortlabels]{enumitem}

\renewcommand{\lrcorner}{\mathrel{\amslrcorner}}        %restore lrcorner in the style of amssymb with extra spacing

\def\vertexsize {1.2pt}       % default radius of a vertex
\def\edgewidth {0.5 pt}      % default width of an edge
\def\halfedgelength {0.7}    %default length of a half-edge

\newcommand{\vertex}[2][1]{\fill (#2) circle [radius = #1 * \vertexsize];}

\newcommand{\tower}[2]{%
    \begin{center}
    \begin{tikzcd}
        \begin{tikzpicture}
            #1
        \end{tikzpicture} \arrow[d] \\
        \begin{tikzpicture}
            #2
        \end{tikzpicture} \arrow[d] \\
        \begin{tikzpicture}
            \halfedge{0}
        \end{tikzpicture}
    \end{tikzcd}
    \end{center}
}

\newcommand{\tetragonal}[1]{%
    \begin{center}
    \begin{tikzcd}
        \begin{tikzpicture}
            #1
        \end{tikzpicture} \arrow[d] \\
        \begin{tikzpicture}
            \halfedge{0}
        \end{tikzpicture}
    \end{tikzcd}
    \end{center}
}

\newcommand{\halfedge}[2][1]{%
    \vertex[#1]{0,#2}
    \path[line width= #1*\edgewidth, draw] (0, #2) -- + (\halfedgelength, 0);
}

\newcommand{\twohalfedges}[2][1]{%
    \vertex[1+#1]{0,#2}     %vertex with two half-edges has at least local degree 2, even 3 if on the of the edges is dilated
    \path[line width = #1 * \edgewidth, draw] (0, #2) -- + (15: \halfedgelength);
    \path[line width = \edgewidth, draw] (0, #2) -- + (-15: \halfedgelength);
}

\newcommand{\threehalfedges}[1]{%
    \vertex[3]{0,#1}
    \path[line width = \edgewidth, draw] (0, #1) -- + (20: \halfedgelength);
    \path[line width = \edgewidth, draw] (0, #1) -- + (0: \halfedgelength);
    \path[line width = \edgewidth, draw] (0, #1) -- + (-20: \halfedgelength);
}

\newcommand{\drawCycle}[2][0.1]{
	\begin{scope}[blue,very thick,decoration={
			markings,
			mark= between positions #1 and 1+#1 step 0.2 with {\arrow{>}}}
		] 
		#2
	\end{scope}	
}

\newcommand{\drawTildeGammaBigonal}{
    \draw (1,-0.5) .. controls (-0.2,0.5) and (-0.2,-1) .. (1,0)               coordinate[midway](vertex_left);
    \draw (1,-0.5) .. controls (1,-1) and (0,-1) .. (-0.1,-0.25)
    			.. controls (0,0.5) and (1,0.5) .. (1,0);
    \path[line width = 2*\edgewidth, draw] (1,-0.5) -- + (1,0);
    \path[line width = 2*\edgewidth, draw] (1,0) -- + (1,0);
    \draw (2.5, -0.5) ellipse (0.5 and 0.2);
    \draw (2.5, 0) ellipse (0.5 and 0.2);
    \path[line width = 2*\edgewidth, draw] (3,-0.5) -- + (1,0);
    \path[line width = 2*\edgewidth, draw] (3,0) -- + (1,0);
    \draw (4,-0.5) .. controls (5.2,0.5) and (5.2,-1) .. (4,0)                 coordinate[midway](vertex_right);;
    \draw (4,-0.5) .. controls (4,-1) and (5,-1) .. (5.1,-0.25)
    			.. controls (5,0.5) and (4,0.5) .. (4,0);
    \foreach \i in {1,2,3,4} {
        \vertex[2]{\i,0}
        \vertex[2]{\i, -0.5}}
    \vertex[2]{-0.1, -0.25}
    \vertex[2]{5.1, -0.25}
    \vertex[2]{vertex_left}
    \vertex[2]{vertex_right}
}

\newcommand{\drawExampleBigonal}{
    \path[line width = 2*\edgewidth, draw] (1,0) -- ++ (-1, -0.25) -- ++ (1, -0.25);
    \draw (1.5, 0) ellipse (0.5 and 0.2);
    \draw (1.5, -0.5) ellipse (0.5 and 0.2);
    \path[line width = 2*\edgewidth, draw] (2, 0) -- ++ (1, -0.25) -- ++ (-1, -0.25);
    \vertex[4]{0, -0.25}
    \vertex[2]{1, 0}
    \vertex[2]{2, 0}
    \vertex[2]{1, -0.5}
    \vertex[2]{2, -0.5}
    \vertex[4]{3, -0.25}
}

\newcommand{\drawExampleBigonalOutput}{
    \path[draw] (0,0) -- ++ (1, 0.5) -- ++ (1,0) -- ++ (1, -0.5) -- ++ (-1, -0.5) -- ++ (-1, 0) -- ++ (-1, 0.5);
    \path[line width = 2*\edgewidth, draw] (0,0) -- ++ (1,0);
    \draw (1.5, 0) ellipse (0.5 and 0.2);
    \path[line width = 2*\edgewidth, draw] (2, 0) -- ++ (1, 0);
    \vertex[4]{0, 0}
    \vertex{1, 0.5}
    \vertex[2]{1, 0}
    \vertex{1, -0.5}
    \vertex{2, 0.5}
    \vertex[2]{2, 0}
    \vertex{2, -0.5}
    \vertex[4]{3, 0}
}

\newcommand{\drawGammaTilde}{
	\foreach \i in {0, 1.5} {
		\path[draw] (5, 0.5+\i) -- ++ (-1, 0.5) -- ++ (-3, 0) -- ++ (-1, -0.5);
		\path[line width = 2*\edgewidth, draw] (0,0.5+\i) -- ++ (1, -0.5) ++ (1,0) -- ++ (1,0) ++ (1,0) -- ++ (1, 0.5);
		\draw (1,0+\i) arc (180:360:0.5 and 0.2);
		\draw (3,0+\i) arc (180:360:0.5 and 0.2);
		\vertex[3]{0, 0.5+\i}
		\vertex[3]{5, 0.5+\i}
		\foreach \x in {1,2,3,4} {
		    \vertex[2]{\x,\i} 
		    \vertex{\x,1+\i}
		}
	}
	\path[draw] (1,0) -- ++ (1,1.5);
	\path[draw] (1,1.5) -- ++ (1,-1.5);
	\path[draw] (3,0) -- ++ (1,1.5);
	\path[draw] (3,1.5) -- ++ (1,-1.5);
}

\newcommand{\drawPi}{
	\path[draw] (0,1.5) -- ++ (1,0) -- ++ (1, -0.5);
	\path[line width = 2*\edgewidth, draw] (2, 1) -- ++ (1,0);
	\path[draw] (3, 1) -- ++ (1, 0.5) -- ++ (1,0);
	\path[draw] (2,1) -- ++ (-2, -1) -- ++ (2, -1) -- ++ (1,0) -- ++ (2,1) -- ++ (-2,1);
	\path[line width = 2*\edgewidth, draw] (0,0) -- ++ (1, -0.5);
	\path[line width = 2*\edgewidth, draw] (4,-0.5) -- ++ (1, 0.5);
	\path[draw] (1,-0.5) -- ++ (1,0.5) -- ++ (1,0) -- ++ (1, -0.5);
	\vertex{0, 1.5}
	\vertex{1, 1.5}
	\vertex[2]{2, 1}
	\vertex[2]{3, 1}
	\vertex{4, 1.5}
	\vertex{5, 1.5}
	
	\vertex{1, 0.5}
	\vertex[3]{0, 0}
	\vertex[2]{1, -0.5}
	\vertex{2, -1}
	\vertex{3, -1}
	\vertex[2]{4, -0.5}
	\vertex[3]{5, 0}
	\vertex{4, 0.5}
	
	\vertex{2, 0}
	\vertex{3, 0}
}

\newtheorem{theorem}{Theorem}[section]
\newtheorem{lemma}[theorem]{Lemma}
\newtheorem{proposition}[theorem]{Proposition}
 
\newtheorem{conjecture}[theorem]{Conjecture} 

\theoremstyle{definition}
								
\newtheorem{definition}[theorem]{Definition}
\newtheorem{example}[theorem]{Example}

\newtheorem{remark}[theorem]{Remark}
\newtheorem*{remark*}{Remark}

\newcommand{\hooklongrightarrow}{\lhook\joinrel\longrightarrow}

\newcommand{\bigmid}{\mathrel{\big\vert}}

\newcommand{\ZZ}{\mathbb{Z}}

\newcommand{\PP}{\mathbb{P}}

\newcommand{\QQ}{\mathbb{Q}}
\newcommand{\RR}{\mathbb{R}}

\renewcommand{\div}{\operatorname{div}}

\newcommand{\oh}{\overline{h}}

\newcommand{\te}{\widetilde{e}}
\newcommand{\tf}{\widetilde{f}}
\newcommand{\tg}{\widetilde{g}}
\renewcommand{\th}{\widetilde{h}}
\newcommand{\thh}{\widetilde{h}}

\newcommand{\tp}{\widetilde{p}}
\newcommand{\tu}{\widetilde{u}}
\newcommand{\tv}{\widetilde{v}}
\newcommand{\tx}{\widetilde{x}}
\newcommand{\ty}{\widetilde{y}}
\newcommand{\tC}{\widetilde{C}}
\newcommand{\tD}{\widetilde{D}}

\newcommand{\tG}{\widetilde{G}}
\newcommand{\tK}{\widetilde{K}}

\newcommand{\tP}{\widetilde{P}}

\newcommand{\tT}{\widetilde{T}}

\newcommand{\tGa}{\widetilde{\Gamma}}
\newcommand{\tPi}{\widetilde{\Pi}}

\newcommand{\al}{\alpha}

\newcommand{\ga}{\gamma}

\newcommand{\ep}{\varepsilon}

\newcommand{\la}{\lambda}

\newcommand{\Ga}{\Gamma}
\newcommand{\La}{\Lambda}
\newcommand{\De}{\Delta}
\newcommand{\Si}{\Sigma}
\newcommand{\Om}{\Omega}

\newcommand{\tal}{\widetilde{\alpha}}
\newcommand{\tbe}{\widetilde{\beta}}
\newcommand{\tga}{\widetilde{\gamma}}

\newcommand{\tep}{\widetilde{\varepsilon}}

\newcommand{\teta}{\widetilde{\eta}}

\newcommand{\IntTor}{\mathrm{IntTor}}

\DeclareMathOperator{\Ker}{Ker}
\DeclareMathOperator{\Coker}{Coker}
\DeclareMathOperator{\Hom}{Hom}
\DeclareMathOperator{\Sym}{Sym}

\let\Im\relax
\DeclareMathOperator{\Im}{Im}
\DeclareMathOperator{\cyc}{cyc}

\newcommand{\calL}{\mathcal{L}}
\newcommand{\calM}{\mathcal{M}}

\DeclareMathOperator{\Pic}{Pic}

\DeclareMathOperator{\Prin}{Prin}

\DeclareMathOperator{\val}{val}

\DeclareMathOperator{\rk}{rk}
\DeclareMathOperator{\Id}{Id}
\DeclareMathOperator{\Div}{Div}
\DeclareMathOperator{\calDiv}{\mathcal{D}iv}

\DeclareMathOperator{\Nm}{Nm}           %Norm homomorphism Jac --> Jac

\DeclareMathOperator{\Prym}{Prym}       %Prym variety (algebraic and tropical)

\DeclareMathOperator{\Aff}{Aff}         %sheaf of affine linear functions
\DeclareMathOperator{\slope}{slope}     %slope of PL function on an edge. Should probably be changed later on.

\newcommand{\reg}{\mathrm{reg}}         %regular point in rational polyhedral spaces. In some sources people write max instead. Maybe change?

\newcommand{\aug}{\mathrm{aug}}
\newcommand{\dil}{\mathrm{dil}}         %dilation subgraph
\newcommand{\torf}{\mathrm{tf}}
\newcommand{\sat}{\mathrm{sat}}

\DeclareMathOperator{\Jac}{Jac}

\let\tilde\relax
\newcommand{\tilde}[1]{\widetilde{#1}}

\newcommand{\Dmitry}[1]{{\color{blue}{\texttt Dmitry: #1}}}

\newcommand{\Felix}[1]{{\color{magenta}{\texttt Felix: #1}}}

\title[The tropical $n$-gonal construction]{The tropical $n$-gonal construction}

\author{Felix R\"ohrle}
\address{Universit\"at T\"ubingen, Fachbereich Mathematik,  72076 T\"ubingen, Germany}
\email{\href{mailto:roehrle@math.uni-tuebingen.de}{roehrle@math.uni-tuebingen.de}}

\author{Dmitry Zakharov}
\address{Department of Mathematics, Central Michigan University, Mount Pleasant, MI 48859, USA}
\email{\href{mailto:dvzakharov@gmail.com}{dvzakharov@gmail.com}}


\begin{document}

\begin{abstract}

    We give a purely tropical analogue of Donagi's $n$-gonal construction and investigate its combinatorial properties. The input of the construction is a harmonic double cover of an $n$-gonal tropical curve. For $n = 2$ and a dilated double cover, the output is a tower of the same type, and we show that the Prym varieties of the two double covers are dual tropical abelian varieties. For $n=3$ and a free double cover, the output is a tetragonal tropical curve with dilation profile nowhere $(2,2)$ or $(4)$, and we show that the construction can be reversed. Furthermore, the Prym variety of the double cover and the Jacobian of the tetragonal curve are isomorphic as principally polarized tropical abelian varieties. Our main tool is tropical homology theory, and our proofs closely follow the algebraic versions.
    
\end{abstract}

\maketitle
\setcounter{tocdepth}{1}
\tableofcontents

%%%%%%%%%%%%%%%%%%%%%%%%%%%%%%%%%%%%%%%%%%%%%%%%%%%
\section{Introduction}
%%%%%%%%%%%%%%%%%%%%%%%%%%%%%%%%%%%%%%%%%%%%%%%%%%%

Tropical geometry aims to find polyhedral, piecewise-linear analogues of the objects studied in algebraic geometry. There are two kinds of algebraic objects for which this correspondence is particularly well-developed. The tropical analogues of algebraic curves are metric graphs, which are the subject of an extensive theory, starting with the seminal paper~\cite{MikhalkinZharkov}. The tropical analogue of an abelian variety is a real torus with additional integral structure, and tropical abelian varieties have perhaps received less attention.

There are two standard ways to associate principally polarized abelian varieties (ppavs) to algebraic curves, and both of these constructions carry over to the tropical setting. The Jacobian variety $\Jac(C)$ of a smooth algebraic curve $C$ of genus $g$ is a ppav of dimension $g$. The corresponding Torelli map $\mathcal{M}_g\to \mathcal{A}_g$ on the moduli spaces is injective, so a smooth algebraic curve can be recovered from its Jacobian. The tropical Jacobian of a metric graph was already introduced in~\cite{MikhalkinZharkov}. The tropical Torelli map $\mathcal{M}_g^{\mathrm{trop}}\to \mathcal{A}_g^{\mathrm{trop}}$ is no longer injective, and its non-injectivity locus was completely described in~\cite{CaporasoViviani}.

The Prym variety $\Prym(\tC/C)$ is a ppav of dimension $g$ associated to an \'etale double cover $\tC\to C$ of algebraic curves of genera $2g+1$ and $g+1$, respectively. It is defined as the connected component of the identity of the kernel of the norm map $\Nm:\Jac(\tC)\to \Jac(C)$, and carries a principal polarization that is half of the polarization induced from $\Jac(\tC)$. The Prym--Torelli map $\mathcal{R}_{g+1}\to \mathcal{A}_g$ on the corresponding moduli spaces is no longer injective (for example, it has positive-dimensional fibers for $g\leq 4$), and its fibers have been extensively studied~\cite{Donagi_tetragonal}. The tropical Prym variety $\Prym(\tGa/\Ga)$ associated to a harmonic double cover of metric graphs $\tGa\to\Ga$ is defined in a completely analogous manner (see~\cite{JensenLen}, ~\cite{LenUlirsch}, and~\cite{LenZakharov}). 
As for the question of describing the fibers of the tropical Prym-Torelli map $\mathcal{R}^{\mathrm{trop}}_{g+1}\to \mathcal{A}^{\mathrm{trop}}_g$ we develop some ideas in \cite{MatroidalPerspective}, however this is still far from a complete description.

There are several remarkable constructions concerning Prym varieties of double covers of curves of small gonality. Let $p:\tC\to C$ be a (possibly ramified) double cover of a smooth curve $C$ admitting a degree $n$ map $f:C\to \PP^1$. As $x$ varies over $\PP^1$, the sections of the fiber maps $(f\circ p)^{-1}(x)\to f^{-1}(x)$ glue together into a $2^n$-gonal curve $\widetilde{D}\to \PP^1$, which carries two additional structures. First, exchanging the sheets of $p$ defines an involution on $\widetilde{D}$, hence a double cover $\widetilde{D}\to D$ of a $2^{n-1}$-gonal curve $D\to \PP^1$. Second, if $p$ is \'etale, then sections have a well-defined parity and $\widetilde{D}$ decomposes as a disjoint union of two $2^{n-1}$-gonal curves (if $n$ is odd then they are exchanged by the involution, while if $n$ is even then each is itself a double cover of a $2^{n-2}$-gonal curve). For $n\leq 4$, and under certain restrictions on the ramification, the ppavs associated to these double covers satisfy a number of natural isomorphisms (due ultimately to the symmetries of the associated Weyl groups):

\begin{enumerate}
    \item[$n = 2$:] Given a double cover $p:\tC\to C$ of a hyperelliptic curve $C$, we obtain another such double cover $\tD\to D$, and applying the construction to $\tD\to D$ recovers the original double cover. The double cover $p$ is \'etale if and only if $\tD\to D$ is split. If $p$ is ramified, then the Prym varieties $\Prym(\tC/C)$ and $\Prym(\tD/D)$ (which are not in general principally polarized) are dual to one another~\cite{Pantazis}.
    
    \item[$n = 3$:] This case was the first to be described~\cite{Recillas}. Given an \'etale double cover $p:\tC\to C$ of a trigonal curve, the trigonal construction produces a tetragonal curve $D$ (the double cover $\tD\to D$ being split), which is \emph{generic} in the sense that the tetragonal map has no fibers with ramification profile $(4)$ or $(2,2)$. This construction can be reversed, and the Prym variety $\Prym(\tC/C)$ is isomorphic to the Jacobian $\Jac(D)$. The Prym variety $\Prym(\tC/C)$ is also principally polarized when $p$ is ramified at two points, and the ramified trigonal construction was described by Dalalyan in~\cite{dalalyan79} and~\cite{dalalyan84}, and was recently rediscovered in~\cite{lange2019trigonal}.
    
    \item[$n = 4$:] This is the general construction, from which the other two may be derived (see~\cite{Donagi_anouncement} and~\cite{Donagi_tetragonal}). Given an \'etale double cover $p:\tC\to C$ of a generic tetragonal curve, we obtain two more such double covers $\tC'\to C'$ and $\tC''\to C''$, and the Prym varieties $\Prym(\tC/C)$, $\Prym(\tC'/C')$, and $\Prym(\tC''/C'')$ are isomorphic. 

\end{enumerate} \medskip

The purpose of our paper is to give tropical versions of the bigonal and trigonal constructions (we outline the tropical tetragonal construction as well, but leave the details and proofs to a future paper). To describe our results, we first discuss the tropical notion of gonality. Baker and Norine~\cite{baker2007riemann} define the \emph{combinatorial rank} $r(D)$ of a divisor $D$ on a finite graph, and show that it satisfies a Riemann--Roch theorem. This was extended to metric graphs by~\cite{GathmannKerber} and~\cite{MikhalkinZharkov}. One may therefore say that a metric graph $\Ga$ is $n$-gonal if it carries a divisor $D$ of degree $n$ and rank $r(D)\geq 1$. This definition is not appropriate in our setting. Instead, following the papers~\cite{Caporaso_gonality} amd~\cite{cools2018metric}, we say that a tropical curve $\Ga$ is \emph{$n$-gonal} if it admits a finite harmonic morphism $\Ga\to K$ of degree $n$, where $K$ is a metric tree. Any fiber of such a map is a divisor of rank at least one (see~\cite{ABBRII}, Proposition 4.2), so this definition of gonality is more restrictive. We note that~\cite{cools2018metric} further require the $n$-gonal map to be  \emph{effective}, which is a numerical condition imposed on the vertices with local degree $d_f(v)\geq 2$. This condition does not play a role in the $n$-gonal construction, and we do not impose it.

In complete analogy to the algebraic case, we consider a tower
\[
\tilde \Gamma \overset{\pi}{\longrightarrow} \Gamma \overset{f}{\longrightarrow} K
\]
of harmonic morphisms of metric graphs, where $K$ is a metric tree and the degrees are $\deg \pi = 2$ and $\deg f = n$, respectively. To this tower we associate, in a purely combinatorial way, a metric graph $\tilde \Pi$ together with a harmonic map $\tilde \Pi \to K$ of degree $2^n$. This map factors as $\tilde\Pi\to \tilde K\to K$, where the \emph{orientation double cover} $\tilde K\to K$ is free (and hence split because $K$ is a tree) if and only if $\pi$ is free. In addition, there is a natural involution on $\tilde \Pi$ with quotient map $\tilde \Pi\to \Pi$. For $n=2$ we have $\Pi=\tilde K$ and hence a tower $\tilde \Pi\to \Pi\to K$ of the same kind as the original tower, which we call the \emph{tropical bigonal construction}. For $n=3$ and $\pi$ free, we instead obtain (a split double cover of) a tetragonal curve $\Pi\to K$, which we call the \emph{tropical trigonal construction}. This construction can be inverted by the \emph{tropical Recillas construction}, under certain restrictions on the fibers of the tetragonal map.

In order to state our results, we first clarify the issue of principal polarizations on tropical abelian varieties. Given a harmonic double cover of graphs $\pi:\tGa\to \Ga$, the Prym variety $\Prym(\tGa/\Ga)$ (as defined in~\cite{JensenLen}) carries a polarization induced from the principal polarization on $\Jac(\tGa)$. The induced polarization is twice a principal polarization if $\pi$ is a free double cover, or if it is dilated and the dilation subgraph is connected, but not for general dilated double covers. A related problem, explored in~\cite{2024GhoshZakharov}, is that $\Prym(\tGa/\Ga)$ behaves discontinuously under edge contractions, specifically those that create additional connected components of the dilation subgraph. This suggests that we ought to modify the definition of $\Prym(\tGa/\Ga)$, which we henceforth call the \emph{divisorial Prym variety} of the double cover $\pi:\tGa\to\Ga$ and denote $\Prym_d(\tGa/\Ga)$. To this end, we introduce an alternative object, the \emph{continuous Prym variety} $\Prym_c(\tGa/\Ga)$ of the double cover $\pi:\tGa\to \Ga$. The continuous Prym variety always carries a natural principal polarization and comes with a natural isogeny $\Prym_c(\tGa/\Ga)\to \Prym_d(\tGa/\Ga)$ of degree $2^{d-1}$, where $d$ is the number of connected components of the dilation subgraph of $\Ga$ (see Proposition~\ref{prop:Prymallpropetries}). In particular, $\Prym_c(\tGa/\Ga)=\Prym_d(\tGa/\Ga)$ if the dilation subgraph is connected or if the double cover is free. Furthermore, $\Prym_c(\tGa/\Ga)$ satisfies the expected universal property of the Prym variety (Proposition~\ref{prop:universal_property_Prymc}), while $\Prym_d(\tGa/\Ga)$ does not. Finally, we note that $\Prym_c(\tGa/\Ga)$ behaves continuously under edge contractions, though we do not use this in our paper.

%\Dmitry{need to change this} In order to state our results, we first clarify the issue of principal polarizations on tropical abelian varieties, and give a modified definition of the tropical Prym variety. Given an integral torus $\Sigma$ with a polarization, we canonically construct in Lemma~\ref{lem:induced_principal_pol} an integral torus $\Sigma^{pp}$ and a map $f:\Sigma^{pp}\to \Sigma$ such that the induced polarization on $\Sigma^{pp}$ is principal. The map $f$ is bijective on points but is dilated, hence not in general invertible as a map of integral tori.  Now let $\pi:\tGa\to \Ga$ be a harmonic double cover of tropical curves, then the kernel of the norm map $\Nm:\Jac(\tGa)\to \Jac(\Ga)$ has two connected components if $\pi$ is free and one if $\pi$ is dilated (see~\cite{JensenLen}). The even connected component $(\Ker \Nm)_0$ has a polarization induced from $\Jac(\tGa)$, which is known to not be principal in general (see~\cite{LenUlirsch}), and as an auxiliary result we compute its polarization type in Proposition~\ref{prop:pp}. We then define the Prym variety $\Prym(\tGa/\Ga)$ as the principally polarized tropical abelian variety $(\Ker \Nm)_0^{pp}$, and we work with both objects: the torus $(\Ker \Nm)_0$ (which is not principally polarized in general) and the principally polarized torus $\Prym(\tGa/\Ga)$.

We can now state the main results of our paper, which are exact analogues of the results of~\cite{Pantazis} and~\cite{Recillas}. First, we state the tropical analogue of the bigonal construction.

\begin{theorem}[Theorem~\ref{thm:bigonal_construction_restated}] \label{thm:bigonal_construction}
    Let $\tilde \Gamma \overset{\pi}{\longrightarrow} \Gamma \overset{f}{\longrightarrow} K$ be a tower of harmonic morphisms of metric graphs of degrees $\deg \pi = \deg f = 2$, where $K$ is a metric tree. Assume that there is no point $x\in K$ with the property that $|f^{-1}(x)|=2$ and $|(f\circ\pi)^{-1}(x)|=2$. Then the output $\tilde \Pi \overset{\pi'}{\longrightarrow} \Pi \overset{f'}{\longrightarrow} K$ of the bigonal construction has the same property, and applying the bigonal construction to it reproduces the original tower. If moreover $\tilde \Gamma$ and $\tilde \Pi$ are both connected, then there is an isomorphism of polarized tropical abelian varieties
    \[ \Prym_d(\tPi/\Pi)^\vee\cong \Prym_d(\tGa/\Ga),\]
    where the polarization on $\Prym_d(\tPi/\Pi)^\vee$ is the pullback of the principal polarization on $\Prym_c(\tPi/\Pi)^\vee$. 
\end{theorem}

Although the statement of the theorem involves only the divisorial Prym varieties of the two double covers, the proof requires using the continuous Prym varieties, and they are also needed to even define the polarization on the dual $\Prym_d(\tPi/\Pi)^\vee$. On the other hand, for the tropical analogue of the trigonal construction, we consider only free double covers, hence there is no need to distinguish the continuous and divisorial Pryms.

\begin{theorem}[Theorem~\ref{thm:tropical_Recillas_theorem_restated}] \label{thm:tropical_Recillas_theorem}
    Let $K$ be a metric tree. The tropical trigonal and Recillas constructions establish a one-to-one correspondence
    \begin{equation*}
        \begin{tikzcd}[column sep = 4cm, every label/.append style = {font = \normalsize}]
        \left\{ 
        \begin{minipage}{0.3\textwidth}
            \begin{center}
                Tropical curves $\Pi$ with a harmonic map of degree 4 to $K$ with dilation profiles nowhere $(4)$ or $(2,2)$.
            \end{center}
        \end{minipage}
        \right\}
        \arrow[r, "\text{Recillas construction}", bend left = 20, start anchor = {[yshift = 1ex]east}, end anchor = {[yshift = 1ex]west}]
        &
        \left\{ 
        \begin{minipage}{0.3\textwidth}
            \begin{center}
                Free double covers $\tGa \to \Gamma$ with a harmonic map of degree 3 from $\Gamma$ to $K$.
            \end{center}    
        \end{minipage}
        \right\}
        \arrow[l, "\text{trigonal construction}", bend left = 20, start anchor = {[yshift = -1ex]west}, end anchor = {[yshift = -1ex]east}]
        \end{tikzcd}
    \end{equation*}
    and under this correspondence, the Prym variety $\Prym(\tilde \Gamma / \Gamma)$ of a double cover and the Jacobian $\Jac(\Pi)$ of the associated tetragonal curve are isomorphic as principally polarized tropical abelian varieties. 
\end{theorem}

We would like to highlight the techniques that we use. Abelian varieties and maps between them are strongly constrained by intersection theory: for example, one may check that an isogeny is an isomorphism by computing its degree in homology. The tropical analogue of singular homology for rational polyhedral spaces was introduced in~\cite{IKMZ}. Tropical homology was first applied to tropical abelian varieties by Gross and Shokrieh~\cite{GrossShokrieh_Poincare}, who established a number of fundamental results about tropical abelian varieties and proved a tropical version of the Poincar\'e formula for the class of a metric graph in the homology of its Jacobian.
The techniques of tropical homology, at least as they apply to abelian varieties, turn out to be quite powerful: in the proof of our main Theorem~\ref{thm:tropical_Recillas_theorem}, we are able to translate the corresponding algebraic proof (see~\cite{Recillas} or~\cite[Theorem 12.7.2]{BirkenhakeLange}) nearly line-by-line into the tropical setting.

\subsection*{Directions for future research}
In Section~\ref{sec:tetragonal} we describe some basic properties of the tropical tetragonal construction, however we do not yet have a tropical analogue for Donagi's theorem relating the construction with Prym varieties.
Algebraically, both the Recillas and Donagi theorems were originally proved by working in cohomology (see~\cite[Theorem~12.8.2]{BirkenhakeLange} for Donagi's theorem). 
It is striking that the techniques of tropical homology allow us to translate the algebraic proof of Recillas's theorem into the tropical setting to give us the proof of Theorem~\ref{thm:tropical_Recillas_theorem}. We are confident that our techniques also work for the tetragonal construction, but a number of additional results in tropical intersection theory will first need to be established, most notably a tropical version of a formula of Macdonald for the class of a $g^1_4$ in the fourth symmetric power of a curve.

Towards a comparison of our tropical construction and the algebraic role model, we emphasize again that our work is a tropical analogue in the sense that our definition carries the underlying geometric ideas of the algebraic construction over to the tropical setting. A priori, it is not clear whether our definition can be recovered from the algebraic one via tropicalization. For example, for certain ramification profiles in the input tower, the algebraic tetragonal construction produces singular output from smooth input data (see~\cite[Local pictures~2.14]{Donagi_tetragonal}). It would be very interesting to see how this phenomenon behaves under tropicalization.

Donagi's original motivation for introducing the tetragonal construction was to study the fibers of the Prym--Torelli map $\mathcal{R}_{g+1} \to \mathcal{A}_g$. The main theorem of the tetragonal construction implies that the Prym--Torelli map is never injective, since double covers related by the tetragonal construction have isomorphic Pryms. Donagi conjectured in~\cite{Donagi_tetragonal} that the tetragonal construction fully accounts for the non-injectivity of the Prym--Torelli map in $g\geq 7$, however, Izadi and Lange~\cite{izadilange} show that the target $\mathbb{P}^1$ in the tetragonal construction can be replaced with an arbitrary curve.

A tropical analogue of Donagi's theorem would be a first step towards understanding the non-injectivity of the tropical Prym--Torelli map. We note, however, that it is elementary to construct Torelli non-injectivity loci in the tropical setting. For example, from \cite{MikhalkinZharkov,CaporasoViviani} we know that bridge edges can be contracted without changing the Jacobian of a tropical curve. It follows immediately that in a double cover $\tilde \Gamma \to \Gamma$ with $\Gamma$ having a bridge edge $e$ whose preimage edges $\tilde e^+$ and $\tilde e^-$ are bridge edges in $\tilde\Gamma$, one may contract $e$, $\tilde e^+$, and $\tilde e^-$ simultaneously without changing the Prym variety. More generally, in~\cite{MatroidalPerspective} we explored the non-injectivity of the tropical Prym--Torelli map in terms of Zaslavsky's signed graphic matroid. Giving a complete description of the fibers of the tropical Prym--Torelli map (the corresponding problem for the tropical Torelli map was solved in~\cite{CaporasoViviani}) will be the topic of future research.

\subsection*{Organization of the article}
We start by introducing the tropical $n$-gonal construction in Section~\ref{sec:graphs}. The construction is purely combinatorial and can be understood without any prior knowledge of tropical geometry or the algebraic $n$-gonal construction. To simplify the exposition, we work with graphs without edge lengths, and passing to metric graphs involves nothing more than equipping the target tree with an edge length function. We conclude by proving the first parts of Theorems~\ref{thm:bigonal_construction} (Propositions~\ref{prop:bigonal_properties} and~\ref{prop:bigonalconnected}) and~\ref{thm:tropical_Recillas_theorem} (Proposition~\ref{prop:trigonal_properties}).

To establish the isomorphisms of the tropical abelian varieties, we first introduce the necessary background on tropical curves, rational polyhedral spaces and tropical homology in Section~\ref{sec:tropical}. Section~\ref{sec:tavs} is devoted to tropical abelian varieties. We speak extensively about real tori with integral structure and develop a theory of morphisms between such objects. In particular, we refine the definition of the tropical Prym variety compared to the existing literature (see e.g.~\cite{LenUlirsch}), introduce the continuous Prym variety $\Prym_c(\tGa/\Ga)$ of a double cover $\tGa\to \Ga$, and calculate the class of $\tGa$ in the tropical homology of $\Prym_c(\tGa/\Ga)$ (Theorem~\ref{thm:tropical_Welters_criterion}, which is the Prym version of the tropical Poincar\'e formula proved in~\cite{GrossShokrieh_Poincare}). We believe this section to be of independent interest.

Finally, we prove the main parts of Theorems~\ref{thm:bigonal_construction} and~\ref{thm:tropical_Recillas_theorem} in Section~\ref{sec:proof}. The structure of the proof of Theorem~\ref{thm:tropical_Recillas_theorem} closely follows the original proof of the algebraic statement in~\cite{Recillas}, using tropical instead of singular homology. The original proof of the algebraic bigonal construction in~\cite{Pantazis} is based on the tetragonal construction, which we have not yet established, hence to prove Theorem~\ref{thm:bigonal_construction} we instead adapt the arguments that are used in Theorem~\ref{thm:tropical_Recillas_theorem}.

%. Our proof of Theorem~\ref{thm:bigonal_construction} is an adaptation thereof -- the proof for the algebraic statement \cite{Pantazis}

%cannot be ported to tropical geometry in this case because it is based on the tetragonal construction which is not yet established in the tropical language.

%We then lift the construction to actual tropical curves in Section~\ref{sec:tropical} without any effort. Furthermore, we introduce all 

%The material in this article is arranged in three parts. In the first part we introduce the tropical Recillas' construction \ref{con:inverse_trigonal} as well as the tropical $n$-gonal construction \ref{con:tropical_n_gonal}. This part only requires some familiarity with tropical curves and harmonic maps -- notions that are recalled in Section \ref{}. The first part of Theorem \ref{thm:tropical_Recillas_theorem} is proved immediately.

\subsection*{Acknowledgments} The authors would like to thank Martin Ulirsch, Paul Helminck, Dhruv Ranganathan, and Alejandro Vargas for insightful conversations. The authors would like to especially thank Madeline Brandt for extended conversations during the beginning of the project, and Andreas Gross for illuminating discussions about tropical homology. Finally, the authors would like to thank the anonymous referees for extensive comments on the original version of the paper, for suggesting a streamlined version of Section~\ref{sec:graphs}, and for pointing out a serious flaw in our original proof of Theorem~\ref{thm:bigonal_construction}.

F.~R.~has received funding from the Deutsche Forschungsgemeinschaft (DFG, German
Research Foundation) TRR 326 \enquote{Geometry and Arithmetic of Uniformized Structures}, project number
444845124, by the Deutsche Forschungsgemeinschaft (DFG, German Research Foundation) Sachbeihilfe
\enquote{From Riemann surfaces to tropical curves (and back again)}, project number 456557832, and from the
LOEWE grant \enquote{Uniformized Structures in Algebra and Geometry}. This work was partly performed during D.~Z.'s stay in Frankfurt, which was supported by the above sources. 

%\input{construction}
%\input{pptavs}
%\input{homology}

%%%%%%%%%%%%%%%%%%%%%%%%%%%%%%%%%%%%%%%%%%%%%%%%%%%
\section{Graphs and the tropical \texorpdfstring{$n$}{n}-gonal construction}
\label{sec:graphs}
%%%%%%%%%%%%%%%%%%%%%%%%%%%%%%%%%%%%%%%%%%%%%%%%%%%

In this section we recall a number of standard definitions from graph theory, and define the $n$-gonal construction for double covers of degree $n$ harmonic covers of trees. 
While we are primarily interested in the construction applied to metric graphs, the idea behind the definition as well as basic properties are best explained in the language of graphs, with an integer-valued local degree function recording the dilation factors. The description can then be lifted to metric graphs without much effort: assigning edge lengths to the base tree automatically determines edge lengths for all covers via the local degree function.

After describing the construction in general, we give more details on the \emph{bigonal}, \emph{trigonal}, and \emph{tetragonal constructions}, which are the special cases for $n = 2,3,4$, respectively.

\subsection{Graphs and harmonic morphisms}

A \emph{graph} $G$ consists of a finite set of \emph{vertices} $V(G)$, a finite set of \emph{half-edges} $H(G)$, a root map $r:H(G)\to V(G)$, and a fixed-point-free involution $H(G)\to H(G)$ denoted $h\mapsto \oh$. The set of \emph{points} of $G$ is $V(G)\cup H(G)$, which we also denote by $G$ by abuse of notation. An orbit $e=\{h,\oh\}$ of the involution is an \emph{edge} of $G$, and the set of edges is denoted $E(G)$. We allow graphs with loops and multiple edges between a pair of vertices. An \emph{orientation} on $G$ is a choice of ordering $(h,\oh)$ of the half-edges for each edge of $G$. The \emph{tangent space} $T_vG=r^{-1}(v)$ to a vertex $v\in V(G)$ is the set of half-edges rooted at $v$, and the \emph{valence} is $\val(v)=|T_vG|$ (so each loop at $v$ contributes twice). The \emph{genus} of a connected graph is defined as $g(G)=|E(G)|-|V(G)|+1$, and a \emph{tree} is a connected graph of genus zero. %\Dmitry{At the end of the day, check which of these definitions are actually necessary (valence?)}

A \emph{morphism} of graphs $f:G\to K$ is a map of the underlying sets of points which maps vertices to vertices and half-edges to half-edges and which commutes with the root and involution maps. In particular, we only consider \emph{finite} morphisms, i.e.~vertices are sent to vertices and edges to edges, and no edge is contracted. A \emph{harmonic morphism} is a pair consisting of a morphism $f:G\to K$ and a function $d_f:G\to \ZZ_{>0}$, called the \emph{local degree}, satisfying the following properties:
\begin{enumerate}
    \item The degrees on the two half-edges comprising an edge $e=\{h,\oh\}\in E(G)$ are equal, and this quantity $d_f(e)=d_f(h)=d_f(\oh)$ is the \emph{dilation factor} of the edge $e$. 
    \item For any vertex $v\in V(G)$ and any half-edge $h'\in T_{f(v)}K$ we have
\begin{equation}
d_f(v)=\sum_{h\in T_{v}G\cap f^{-1}(h')} d_f(h).
\label{eq:localdegree}
\end{equation}
\end{enumerate}
When $K$ is connected, the sum
\begin{equation}
\deg (f)=\sum_{y\in f^{-1}(x)}d_f(y)
\label{eq:globaldegree}
\end{equation}
is the same for any point $x$ (vertex or half-edge) of $K$ and is called the \emph{(global) degree} of $f$. Given harmonic morphisms $f:G\to K$ and $g:K\to L$ of graphs with degree functions $d_f$ and $d_g$, the composition $g\circ f$ is a harmonic morphism with degree function $d_{g\circ f}$ given by
\[
d_{g\circ f}(x)=d_f(x)d_g(f(x))
\]
for any point $x$ of $G$, and global degree $\deg (g\circ f)=\deg (f)\deg (g)$. We say that a harmonic morphism $f:G\to K$ is \emph{free} if it has local degree 1 everywhere; such morphisms are covering spaces in the topological sense.

A \emph{double cover} $\pi:\tG\to G$ is a harmonic morphism of degree 2. A point $x\in G$ is either \emph{free}, having two preimages $\pi^{-1}(x)=\{\tilde x^+,\tilde x^-\}$ with $d_{\pi}(\tilde x^{\pm})=1$, or \emph{dilated}, having a unique preimage that we label $\pi^{-1}(x)=\{\tilde x^{\pm}\}$ with $d_{\pi}(\tilde x^{\pm})=2$. We note that we can choose the labelings on the preimages of the free points to be \emph{vertex-trivial}, so that $r(\th^{\pm})=\tv^{\pm}$ for any $v\in V(G)$ and any $h\in T_vG$, or \emph{edge-trivial}, so that $\overline{\th^{\pm}}=\widetilde{\overline{h}}^{\pm}$ for any half-edge $h\in H(G)$ (in other words, so that the two edges of $\tG$ over any edge $\{h_1,h_2\}\in E(G)$ are $\{\th^+_1,\th^+_2\}$ and $\{\th^-_1,\th^-_2\}$). 

The double cover $\pi$ is free if it has no dilated points, otherwise we say that $\pi$ is \emph{dilated}. There is an induced graph involution $\iota:\tG\to \tG$ (not to be confused with the involution that pairs the half-edges into edges) that exchanges the preimages of the free points and fixes the preimages of the dilated points. Conversely, given a graph involution $\iota:\tG\to \tG$ that does not flip edges (meaning that $\iota(h) \neq \overline{h}$ for all $h \in H(\tG)$), the quotient map $\tG\to \tG/\iota$ naturally has the structure of a double cover, which is free if and only if the involution has no fixed points. 

%\Dmitry{convention on point labeling: $\tilde y\in \tG$, $y\in G$, $x\in K$, $D\in \tilde P$, $\tilde x\in \tilde K$}

\subsection{The tropical $n$-gonal construction.} 
\label{sec:n_gonal_harmonic}

Let $\pi:\tG\to G$ and $f:G\to K$ be harmonic morphisms of degree $2$ and $n$, respectively. We describe the tropical $n$-gonal construction, in other words we construct a graph $\tilde P$ and a harmonic morphism $\tilde p:\tilde P\to K$ of degree $2^n$, together with certain additional structures. Our exposition is indebted to and closely follows~\cite{Donagi_tetragonal}. 

Informally speaking, the points of $\tilde P$ over a fixed point $x \in K$ represent the sections $f^{-1}(x) \to (f \circ \pi)^{-1}(x)$ of $\pi:\tG\to G$ over the fiber of $x$. To make this precise and correctly assign degrees, we proceed as follows. Consider the free abelian groups $\ZZ^{\tG}$, $\ZZ^{G}$, and $\ZZ^{K}$ generated by the points (vertices and half-edges) of our graphs, and define the \emph{pushforward} and \emph{pullback} homomorphisms on the generators $\tilde y\in \tG$ and $x\in K$ as
\[
\begin{aligned}
    &\pi_*:\ZZ^{\tG}\longrightarrow \ZZ^{G},    
    &&\qquad\text{and}\qquad &
    &f^*:\ZZ^{K}\longrightarrow \ZZ^{G}, \\
    &\pi_*(\tilde y)=\pi(\tilde y),    
    &&&&f^*(x)=\sum_{y\in f^{-1}(x)}d_f(y)y.
\end{aligned}
\]
For $D\in \ZZ^{\tG}$, we write $D\geq 0$ if all coefficients are non-negative. We now set
\[
\tP=\left\{D\in \ZZ^{\tG}:D\geq 0\mbox{ and }\pi_*D=f^*(x)\mbox{ for some (necessarily unique) }x\in K\right\},
\]
and define the map $\tp:\tP\to K$ by $\tp(D)= x$. The point $x\in K$ is a vertex if and only if any $D$ lying above $x$ is a linear combination of vertices of $\tG$, and the same is true for half-edges; this defines the vertex and half-edge sets of $\tP$. Finally, the root and involution maps on $\tP$ are induced from those on $\tG$:
\begin{equation*}
    r\Big(\sum a_{\tilde h} \tilde h \Big) = \sum a_{\tilde h} r(\tilde h) 
    \qquad \text{and} \qquad
    \overline{\sum a_{\tilde h} \tilde h} = \sum a_{\tilde h} \overline{\tilde h}.
\end{equation*}
It is clear that $\tP$ is a (not necessarily connected) graph and that $\tilde p : \tilde P \to K$ is a morphism of graphs. A point $D\in \tP$ lying above a given point $x\in K$ has the form
\[
D=\sum_{y \in f^{-1}(x) \text{ free}}\left(a_{\tilde y^+}\tilde y^++a_{\tilde y^-}\tilde y^-\right)+\sum_{y \in f^{-1}(x) \text{ dilated}}d_f(y)\tilde y^{\pm}.
\]
Here $\pi^{-1}(y)=\{\tilde y^+,\tilde y^-\}$ and $\pi^{-1}(y)=\{\tilde y^{\pm}\}$ denote the preimages of a point $y\in G$ that is respectively free or dilated (with respect to the double cover $\pi$), and the $a_{\tilde y^{\pm}}$ are nonnegative coefficients satisfying $a_{\tilde y^+}+a_{\tilde y^-}=d_f(y)$.

To make $\tp$ harmonic, we define the local degree of $\tp$ at $D\in\tilde P$ to be 
\begin{equation} \label{eq:tpdegree}
    d_{\tilde p}(D) = \prod_{y \in f^{-1}(x) \text{ free}} \binom{d_f(y)}{a_{\tilde y^+}} \prod_{y \in f^{-1}(x) \text{ dilated}} 2^{d_f(y)},
\end{equation}
which does not depend on the choice of labeling of the fiber $\pi^{-1}(y)=\{\tilde y^+,\tilde y^-\}$ by the symmetry of the binomial coefficient. In Proposition~\ref{prop:harmonic} we show that the local degree function defined in Equation~\eqref{eq:tpdegree} gives $\tp$ the structure of a harmonic map, and explain the combinatorial meaning of the coefficient in~\eqref{eq:tpdegree}. For now, we note that the number of preimages of $x\in K$ is equal to
\[
    \big|\tilde p^{-1}(x) \big|=\prod_{y \in f^{-1}(x) \text{ free}} \big( d_f(y) +1 \big).
\]

%where $a_i$ is the coefficient of any $\tilde x_i$ with $\pi (\tilde x_i) = y$. This is well-defined since there are precisely two $\tilde x_i$ mapping to a free $y$ in $G$ and their coefficients in $D$ are necessarily $a_i$ and $d_f(y) - a_i$. Since $\binom{d_f(y)}{d_f(y) - a_i} = \binom{d_f(y)}{a_i}$ it does not matter which preimage $\tilde x_i$ of $y$ we use. 

We now define the \emph{orientation double cover} $k:\tK\to K$ of the tower $\tG\to G\to K$ as a quotient of $\tP$ by an equivalence relation. We say that a point $x\in K$ is \emph{dilated} if any point in $f^{-1}(x)$ is dilated (with respect to the double cover $\pi$), and \emph{free} if all points in $f^{-1}(x)$ are free. For a dilated point $x$, we set all points in $\tilde p^{-1}(x)\subset \tP$ to be equivalent. For a free point $x$ and two points
\[
D=\sum_{y \in f^{-1}(x)}\left(a_{\tilde y^+}\tilde y^++a_{\tilde y^-}\tilde y^-\right),\qquad \text{and}\qquad E=\sum_{y \in f^{-1}(x)}\left(b_{\tilde y^+}\tilde y^++b_{\tilde y^-}\tilde y^-\right)
\]
in $\tilde p^{-1}(x)\subset \tP$, we set $D\sim E$ if $\sum(a_{\tilde y^+}-b_{\tilde y^+})$ is even. This does not depend on the choice of labeling of $\pi^{-1}(y)=\{\tilde y^+,\tilde y^-\}$, because $a_{\tilde y^+}+a_{\tilde y^-}=d_f(y)=b_{\tilde y^+}+b_{\tilde y^-}$ implies that $a_{\tilde y^+}-b_{\tilde y^+}$ is even if and only if $a_{\tilde y^-}-b_{\tilde y^-}$ is even. 

We define $V(\tK)$ and $H(\tK)$ to be the quotients of respectively $V(\tP)$ and $H(\tP)$ by the equivalence relation, and induce the root and involution maps from $\tP$. Choosing either vertex-trivial or edge-trivial labelings for the preimages of the free points in $G$, we see that this is well-defined. Hence $\tK$ is a graph, and the quotient $q:\tP\to \tK$ and projection $k:\tK\to K$ maps are graph morphisms whose composition is $\tp:\tP\to K$. 

We define the local degrees of the orientation double cover $k$ as follows:
\[
d_k(\tilde x)=\begin{cases}
    2, & \text{if }k(\tx)\text{ is dilated,}\\
    1, & \text{if }k(\tx)\text{ is free.}
\end{cases}
\]
It is elementary to verify that $k$ is indeed a harmonic double cover. The orientation double cover is free if and only if the double cover $\pi:\tG\to G$ is free. If in addition $k:\tK\to K$ is a trivial double cover (which is always the case if $K$ is a tree), then we say that the tower $\tG\to G\to K$ is \emph{orientable}.

Similarly, we define the local degrees of $q$ by
\begin{equation}
    d_q(D)=
    \begin{cases}
        d_{\tp}(D)/2, & \text{if } \tp(D) \text{ is dilated},\\
        d_{\tp}(D), & \text{if } \tp(D) \text{ is free},
    \end{cases}
\end{equation}
where $d_{\tp}$ is given by~\eqref{eq:tpdegree}. We now verify that the morphisms $\tp$ and $q$ are harmonic.

\begin{proposition} \label{prop:harmonic}
    The morphisms $\tp:\tP\to K$ and $q:\tP\to \tK$ are harmonic of degrees $2^n$ and $2^{n-1}$, respectively.
\end{proposition}

\begin{proof} 
Pick a vertex $u\in V(K)$ and a half-edge $l\in T_uK$ rooted at $u$. We consider the fibers of $\pi:\tG\to G$ above the pair $(u,l)$. By~\eqref{eq:globaldegree}, the half-edge degrees $\big\{d_f(h) \bigmid h\in f^{-1}(l) \big\}$ are a partition of $n = \deg f$. Hence we can define a labeling
\[
\mu_l : \{1,\ldots,n\}\longrightarrow f^{-1}(l)
\]
of the half-edges of $G$ lying above $l$ such that for each half-edge $h\in f^{-1}(l)$ we have $|\mu^{-1}(h)| = d_f(h)$. We can similarly label the half-edges of $\tG$ lying above $l$:
\[
\tilde\mu_l : \{\pm 1,\ldots,\pm n\}\longrightarrow (f\circ \pi)^{-1}(l),
\]
so that $\tilde\mu_l^{-1}(\tilde h^+) = \big\{+i \bigmid \mu_l(i) = \pi(\tilde h^+) \big\}$ and $\tilde\mu_l^{-1}(\tilde h^-) = \big\{-i \bigmid \mu_l(i) = \pi(\tilde h^-) \big\}$ for free fibers $\{\th^+, \th^-\} = \pi^{-1}(h)$ of $\pi$ and $\tilde\mu_l^{-1}(\tilde h) = \big\{ \pm i \bigmid \mu(i) = \pi(\th^\pm) \big\}$ for dilated fibers $\{\th^\pm\} = \pi^{-1}(h)$. Composing with the root maps $(f\circ \pi)^{-1}(l)\to (f\circ \pi)^{-1}(u)$ and $f^{-1}(l)\to f^{-1}(u)$, we obtain similar labelings for the vertices above $u$:
\[
\mu_u : \{1,\ldots,n\}\longrightarrow f^{-1}(u),\qquad  \tilde\mu_u : \{\pm 1,\ldots,\pm n\}\longrightarrow (f\circ \pi)^{-1}(u).
\]
This system of labelings induces a labeling on the $n$-gonal construction $\tP \to K$ as follows:
\begin{align*}
    \tilde\nu_u : \Big\{ (\epsilon_i)_{i = 1, \ldots, n} \mathrel{\Big|} \epsilon_i \in \{+, -\} \Big\} &\longrightarrow \tp^{-1}(u) \\
    (\pm)_i &\longmapsto \sum \tilde\mu_u(\pm i),
\end{align*}
and similarly $\tilde \nu_l$ for $\tp^{-1}(l)$.
It is elementary to see that the number $\big| \tilde \nu_u^{-1}(D)\big|$ of labelings of a given vertex $D \in \tp^{-1}(u)$ is equal to the coefficient $d_{\tp}(D)$ given by Equation~\eqref{eq:tpdegree}, and similarly $d_{\tp}(E)$ is the number of labelings of a half-edge $E \in \tp^{-1}(l)$. Moreover, we note that for any vertex $D \in \tp^{-1}(u)$ and half-edge $E \in \tp^{-1}(l)$, we have that $r(E)=D$ if and only if $\tilde \nu_l^{-1}(E) \subseteq \tilde \nu_u^{-1}(D)$. Hence the harmonicity condition~\eqref{eq:localdegree} for $\tp$ at the vertex $D$ and the half-edge $l$
\[
d_{\tp}(D)=\sum_{E:r(E)=D\text{ and }\tp(E)=l}d_{\tp}(E)
\]
is now equivalent to 
\[
\tilde \nu_u^{-1}(D) = \bigsqcup_{E:r(E)=D\text{ and }\tp(E)=l} \tilde \nu_l^{-1}(E).
\]
This however is clear since the labeling functions $\tilde \nu$ are compatible with the root map of $\tP$ meaning that $r \circ \tilde\nu_l = \tilde\nu_u$, and hence $\tp$ is harmonic at $D$. Finally, the global degree of $\tp$ is the total number of sections of $\{\pm 1,\ldots,\pm n\}\to \{1,\ldots,n\}$, which is the cardinality of the domain of $\tilde\nu_u$ and equals $2^n$. The harmonicity of $q$ is proved in the same way, but taking parities into account, and we omit the details. \qedhere

\end{proof}

We introduce one additional structure on $\tP$. The involution $\iota:\tG\to \tG$ associated to the double cover $\pi:\tG\to G$ induces the pushforward involution $\iota_*:\tP\to \tP$ which acts by exchanging signs, in other words $\iota_*(\ty^{\pm})=\ty^{\mp}$.
%\Felix{In the new language this is just pushforward -- no need to spell it out.}
%\[
%\iota_p\left[
%\sum_{y \in f^{-1}(x)}\left(a_{\tilde y^+}\tilde y^++a_{\tilde y^-}\tilde y^-%\right)+\sum_{y \in f^{-1}(x)}d_f(y)\tilde y\right]=
%\sum_{y \in f^{-1}(x)}\left(a_{\tilde y^-}\tilde y^++a_{\tilde y^+}\tilde y^-%\right)+\sum_{y \in f^{-1}(x)}d_f(y)\tilde y.
%\]
%i.e. writing the orbits of $\iota$ as $\{\tilde x^+, \tilde x^- \}$ with two not necessarily different points of $\tilde G$, then we have
%\[
%\iota_P\left(\sum \big( a_{+i}\tx_{i}^+ +a_{-i}\tx_{i}^- \big) \right) = \sum \big( a_{+i}\tx_{i}^- +a_{-i}\tx_{i}^+ \big).
%\iota_P\left(\sum(b_j^+\th_j^++b_j^-\th_j^-)\right)=\sum(b_j^+\th_j^-+b_j^-\th_j^+).
%\]
Denoting the quotient by $P=\tP/\iota_*$, we obtain a double cover $\tP\to P$. The involution $\iota_*$ is fixed-point-free (and hence the double cover $\tP\to P$ is free) if and only if over every vertex $v\in V(K)$ there is a free vertex $v_i\in f^{-1}(v)$ having odd local degree $d_f(v_i)$. Hence the double cover $\tP\to P$ need not be free if $\tG\to G$ is free, and vice versa. 

Since the involution $\iota_*$ preserves the local degrees of $\tp$, there is an induced harmonic morphism $p:P\to K$ of degree $2^{n-1}$.  If $n$ is even, then the involution $\iota_*$ preserves the fibers of $\tP\to \tK$, and our morphisms factor into a tower
\[
\tP\overset{2}{\longrightarrow} P\overset{2^{n-2}}{\longrightarrow} \tK\overset{2}{\longrightarrow} K.
\]
On the other hand, if $n$ is odd, then the involution $\iota_\ast$ exchanges the equivalence classes in every fiber (if they are distinct), and we have a diagram
\begin{equation*}
    \begin{tikzcd}[sep = small]
        & \tP \arrow[dl,"2"'] \arrow[dr,"2^{n-1}"] & \\
        P \arrow[dr,"2^{n-1}"'] & & \tK \arrow[dl,"2"] \\
        & K. &
    \end{tikzcd}
\end{equation*}
In particular, if the tower $\tG\to G\to K$ is orientable, then the $n$-gonal construction $\tp:\tP\to K$ splits as two isomorphic copies of the degree $2^{n-1}$ cover $p:P\to K$ that are exchanged by $\iota_\ast$.

\medskip

We now study the $n$-gonal construction in detail for $n=2$, $3$, and $4$, and call it the \emph{bigonal}, \emph{trigonal}, and \emph{tetragonal} construction, respectively. To clarify the exposition, we first consider the case of free covers, for which we can use topological methods (but which are not directly relevant to us, since a free cover of a tree is trivial). The $n$-gonal morphism $\tp:\tP\to K$ is free if and only if both $\pi:\tG\to G$ and $f:G\to K$ are free. In this case, the entire construction may be carried out in terms of covering space theory, realizing graphs as topological spaces by gluing unit intervals in the standard way. 
%graphs as topological spaces in the obvious way \Felix{you mean by idenfifying edges with that unit interval? I don't think that \enquote{obvious} is the right word since you can also endow the set $V \cup E$ of a graph with a topology} \Dmitry{I don't know, "evident"? There is a whole chapter in Hatcher on graphs and their covers; for example, it's how you prove that a subgroup of a free group is free}. \Dmitry{There is actually a nice way to topologize a graph $G=V(G)\cup E(G)$ as a finite set: vertices are closed, and the closure of an edge $e=\{u,v\}$ is $\{e,u,v\}$. The topology is not Hausdorff, but it does satisfy the requirements for covering space theory. This came up in my work with Martin}
The tower $\tG \overset{2}{\rightarrow}G\overset{n}{\rightarrow} K$ corresponds to a monodromy representation $\pi_1(K,x_0)\to S_n^B$ (where $x_0\in K$ is a base point) in the \emph{signed permutation group} $S^B_n \subseteq S_{2n}$. Recall that $S^B_n$ consists of permutations of the $2n$-element set $\{\pm 1,\ldots,\pm n\}$ that commute with the fixed-point-free involution $i\mapsto -i$. Signed permutations act on the $2^n$ sections of the map $\{\pm 1,\ldots,\pm n\}\to \{1,\ldots,n\}$, and the corresponding monodromy representation $\pi_1(K,x_0)\to S_{2^n}$ gives the $n$-gonal cover $\tP\overset{2^n}{\rightarrow} K$. Furthermore, the tower $\tilde G \overset{2}{\rightarrow} G \overset{n}{\rightarrow} K$ is orientable if and only if the monodromy representation lies in the \emph{even signed permutation group} $S_n^D=S_n^B\cap A_{2n}$. The groups $S^B_n$ and $S^D_n$ are the Weyl groups of the root systems $B_n$ and $D_n$, respectively, and the bigonal, trigonal, and tetragonal constructions arise from the symmetries of these root systems for $n=2,3,4$.

More precisely, for $n=2$, the group $S^B_2$ is isomorphic to the dihedral group $D_4$ and carries an outer automorphism, namely conjugation by a $\pi/4$ rotation. Applying the outer automorphism to the monodromy representation $\pi_1(K,x_0)\to S^B_2$ of the tower $\tG\overset{2}{\rightarrow} G\overset{2}{\rightarrow} K$, we obtain a second tower $\tG'\overset{2}{\rightarrow} G'\overset{2}{\rightarrow} K$ of the same type. Iterating the bigonal construction reproduces the original tower, because the square of the outer automorphism is an inner automorphism.

In $n = 3$, the equality $A_3=D_3$ of root systems corresponds to an isomorphism $S^D_3\cong S_4$ of the Weyl groups, which we now describe. The map $\{\pm 1,\pm 2,\pm 3\}\to \{1,2,3\}$ has four even sections (meaning sections such that the product of the images is positive), and the group $S^D_3$ acts on this four-element set. This defines a homomorphism $S^D_3\to S_4$. Conversely, each pair of even sections defines a unique element of $\{\pm 1,\pm 2,\pm 3\}$, namely the intersection of the images, and this gives the inverse homomorphism. Given an orientable tower $\tG\overset{2}{\rightarrow} G\overset{3}{\rightarrow} K$, the $3$-gonal construction $\tP\overset{8}{\rightarrow} K$ splits as two isomorphic copies of a degree 4 cover $P\overset{4}{\rightarrow} K$, and the monodromy representation of the latter is obtained by composing the monodromy representation $\pi_1(K,x_0)\to S^D_3$ of the tower with the isomorphism $S^D_3\cong S_4$ defined above. Conversely, any degree 4 cover $P\overset{4}{\rightarrow} K$ defines an orientable tower $\tG\overset{2}{\rightarrow} G\overset{3}{\rightarrow} K$; we call this inverse correspondence the \emph{Recillas construction}.

Finally, for $n = 4$ the Dynkin diagram $D_4$ has automorphism group $S_3$, which is also the group of outer automorphisms of the Weyl group $S_4^D$ modulo inner automorphisms. Given an orientable tower $\tG\overset{2}{\rightarrow} G\overset{4}{\rightarrow} K$, the orientation double cover $\tK\overset{2}{\rightarrow}K$ splits, hence the degree 16 morphism $\tP\overset{2}{\rightarrow} P\overset{4}{\rightarrow} \tK \overset{2}{\rightarrow} K$ splits into two towers $\tG'\overset{2}{\rightarrow} G'\overset{4}{\rightarrow} K$ and $\tG''\overset{2}{\rightarrow} G''\overset{4}{\rightarrow} K$ of the same type as the original tower. The monodromy representations of the three towers are related by the outer automorphisms, and the construction is a triality, in other words each of the towers reproduces the other two.

We now move on to a description in the general, dilated case.

\subsection{The tropical bigonal construction} 
Consider a tower $\tG\to G\to K$ of harmonic double covers of graphs. We first introduce a classification system for points of $K$:

\begin{definition} 
    Let $\tG\to G\to K$ be a tower of harmonic double covers. A point $x\in K$ is called

    \begin{enumerate}
        \item \emph{Type I} if it has a unique preimage in $G$, which in turn has a unique preimage in $\tG$.
        \item \emph{Type II} if it has a unique preimage in $G$, which has two preimages in $\tG$.
        \item \emph{Type III} if it has two preimages in $G$, one having two preimages in $\tG$, the other having one. 
        \item \emph{Type IV} if it has two preimages in $G$, each having two preimages in $\tG$.
        \item \emph{Type V} if it has two preimages in $G$, each having a unique preimage in $\tG$.
    \end{enumerate}

\end{definition}
We note that for types I through IV, the type is the number of preimages in $\tG$. Let $\tP\to P\to K$ be the bigonal construction (see Subsection~\ref{sec:n_gonal_harmonic} with $n = 2$) associated to the tower $\tG\to G\to K$. A case-by-case verification shows that the types of points of $K$ with respect to the two towers change as follows (see Figure~\ref{fig:bigonaltypes}):
\[
I\to I,\qquad II\to III,\qquad III\to II,\qquad IV\to IV,\qquad V\to I.
\]

\begin{figure}[htb]
    \centering
    {\renewcommand{\arraystretch}{4}%
\setlength{\tabcolsep}{1em}%
\begin{tabular}{c|c}
    $\tilde G \longrightarrow G \longrightarrow K$ & 
    $\tilde P \longrightarrow P \longrightarrow K$ \\
    \hline
    \begin{tikzpicture}
        \draw (0, 0) node[anchor = east] {Type I};
        
        \vertex[4]{1.5, 0}
        \draw (1.4, 0) node[anchor = east] {\small$\tilde y_1^{\pm}$};
        \vertex[2]{3, 0}
        \draw (2.9, 0) node[anchor = east] {\small$y_1$};
        \vertex{4.5, 0}
        \draw (4.4, 0) node[anchor = east] {\small$x$};
    \end{tikzpicture} 
    & 
    \begin{tikzpicture}
        \draw (-0.5, 0) node[anchor = east] {Type I};
        
        \vertex[4]{1.5, 0}
        \draw (1.4, 0) node[anchor = east] {\small$2\tilde y_1^{\pm}$};
        \vertex[2]{3, 0}
        \vertex{4.5, 0}
        \draw (4.4, 0) node[anchor = east] {\small$x$};
    \end{tikzpicture} 
    \\ \hline
    \begin{tikzpicture}
        \draw (0, 0) node[anchor = east] {Type II};
        
        \vertex[2]{1.5, 0.5}
        \draw (1.4, 0.5) node[anchor = east] {\small$\tilde y_1^+$};
        \vertex[2]{1.5, -0.5}
        \draw (1.4, -0.5) node[anchor = east] {\small$\tilde y_1^-$};
        
        \vertex[2]{3, 0}
        \draw (2.9, 0) node[anchor = east] {\small$y_1$};
        
        \vertex{4.5, 0}
        \draw (4.4, 0) node[anchor = east] {\small$x$};
    \end{tikzpicture} 
    & 
    \begin{tikzpicture}
        \draw (-0.5, 0) node[anchor = east] {Type III};
        
        \vertex{1.5, 0.7}
        \draw (1.4, 0.7) node[anchor = east] {\small$2 \tilde y_1^+$};
        \vertex{1.5, 0.3}
        \draw (1.4, 0.3) node[anchor = east] {\small$2\tilde y_1^-$};
        \vertex[2]{1.5, -0.5}
        \draw (1.4, -0.5) node[anchor = east] {\small$\tilde y_1^+ + \tilde y_1^-$};
        
        \vertex{3, 0.5}
        \vertex{3, -0.5}
        
        \vertex{4.5, 0}
        \draw (4.4, 0) node[anchor = east] {\small$x$};
    \end{tikzpicture} 
    \\ \hline
    \begin{tikzpicture}
        \draw (0, 0) node[anchor = east] {Type III};
        
        \vertex{1.5, 0.7}
        \draw (1.4, 0.7) node[anchor = east] {\small$\tilde y_1^+$};
        \vertex{1.5, 0.3}
        \draw (1.4, 0.3) node[anchor = east] {\small$\tilde y_1^-$};
        \vertex[2]{1.5, -0.5}
        \draw (1.4, -0.5) node[anchor = east] {\small$\tilde y_2^{\pm}$};
        
        \vertex{3, 0.5}
        \draw (2.9, 0.5) node[anchor = east] {\small$y_1$};
        \vertex{3, -0.5}
        \draw (2.9, -0.5) node[anchor = east] {\small$y_2$};
        
        \vertex{4.5, 0}
        \draw (4.4, 0) node[anchor = east] {\small$x$};
    \end{tikzpicture} 
    & 
    \begin{tikzpicture}
        \draw (-0.5, 0) node[anchor = east] {Type II};
        
        \vertex[2]{1.5, 0.5}
        \draw (1.4, 0.5) node[anchor = east] {\small$\tilde y_1^+ + \tilde y_2^{\pm}$};
        \vertex[2]{1.5, -0.5}
        \draw (1.4, -0.5) node[anchor = east] {\small$\tilde y_1^- + \tilde y_2^{\pm}$};
        
        \vertex[2]{3, 0}
        
        \vertex{4.5, 0}
        \draw (4.4, 0) node[anchor = east] {\small$x$};
    \end{tikzpicture} \\ \hline
    \begin{tikzpicture}
        \draw (0, 0) node[anchor = east] {Type IV};
        
        \vertex{1.5, 0.7}
        \draw (1.4, 0.7) node[anchor = east] {\small$\tilde y_1^+$};
        \vertex{1.5, 0.3}
        \draw (1.4, 0.3) node[anchor = east] {\small$\tilde y_1^-$};
        \vertex{1.5, -0.3}
        \draw (1.4, -0.3) node[anchor = east] {\small$\tilde y_2^+$};
        \vertex{1.5, -0.7}
        \draw (1.4, -0.7) node[anchor = east] {\small$\tilde y_2^-$};
        
        \vertex{3, 0.5}
        \draw (2.9, 0.5) node[anchor = east] {\small$y_1$};
        \vertex{3, -0.5}
        \draw (2.9, -0.5) node[anchor = east] {\small$y_2$};
        
        \vertex{4.5, 0}
        \draw (4.4, 0) node[anchor = east] {\small$x$};
    \end{tikzpicture} 
    &
    \begin{tikzpicture}
        \draw (-0.5, 0) node[anchor = east] {Type IV};
        
        \vertex{1.5, 0.7}
        \draw (1.4, 0.7) node[anchor = east] {\small$\tilde y_1^+ + \tilde y_2^+$};
        \vertex{1.5, 0.3}
        \draw (1.4, 0.3) node[anchor = east] {\small$\tilde y_1^- + \tilde y_2^-$};
        \vertex{1.5, -0.3}
        \draw (1.4, -0.3) node[anchor = east] {\small$\tilde y_1^+ + \tilde y_2^-$};
        \vertex{1.5, -0.7}
        \draw (1.4, -0.7) node[anchor = east] {\small$\tilde y_1^- + \tilde y_2^+$};
        
        \vertex{3, 0.5}
        \vertex{3, -0.5}
        
        \vertex{4.5, 0}
        \draw (4.4, 0) node[anchor = east] {\small$x$};
    \end{tikzpicture} \\ \hline
    \begin{tikzpicture}
        \draw (0, 0) node[anchor = east] {Type V};
        
        \vertex[2]{1.5, 0.5}
        \draw (1.4, 0.5) node[anchor = east] {\small$\tilde y_1^{\pm}$};
        \vertex[2]{1.5, -0.5}
        \draw (1.4, -0.5) node[anchor = east] {\small$\tilde y_2^{\pm}$};
        
        \vertex{3, 0.5}
        \draw (2.9, 0.5) node[anchor = east] {\small$y_1$};
        \vertex{3, -0.5}
        \draw (2.9, -0.5) node[anchor = east] {\small$y_2$};
        
        \vertex{4.5, 0}
        \draw (4.4, 0) node[anchor = east] {\small$x$};
    \end{tikzpicture} 
    & 
    \begin{tikzpicture}
        \draw (-0.5, 0) node[anchor = east] {Type I};
        
        \vertex[4]{1.5, 0}
        \draw (1.4, 0) node[anchor = east] {\small$\tilde y_1^{\pm} + \tilde y_2^{\pm}$};
        \vertex[2]{3, 0}
        \vertex{4.5, 0}
        \draw (4.4, 0) node[anchor = east] {\small$x$};
    \end{tikzpicture}
\end{tabular}}
    \caption{The structure of the bigonal construction locally over a point $x \in K$. The size of points indicates the dilation factor with respect to $K$.}
    \label{fig:bigonaltypes}
\end{figure}

\noindent
We immediately observe that the tropical bigonal construction is not invertible, since type V and type I points with respect to $\tG\to G\to K$ both produce type I points with respect to $\tP\to P\to K$. This phenomenon, which we call \emph{dilation collapse}, also occurs for the trigonal construction (see Remark~\ref{rem:dilationcollapse}) and forces us to introduce restrictions on the dilation of the $n$-gonal map.

\begin{definition} 
    A tower $\tG\to G\to K$ of harmonic double covers is called \emph{generic} if $K$ has no points of type V.

\end{definition}

Restricted to generic towers, the bigonal construction is an involution, and dilation behavior in fibers is exchanged (this is the tropical analogue of Lemma 2.7 in~\cite{Donagi_tetragonal}).

%\Dmitry{Think about retaining Figure 2}
%\begin{figure}[htb]
%    \centering
%    \subfile{tikz-pictures/bigonal_construction}
%    \caption{Half-edge local pictures of the tropical bigonal construction. Thickness of edges and %vertices corresponds to dilation.}
%    \label{tab:bigonal_construction}
%\end{figure}

\begin{proposition} 
    \label{prop:bigonal_properties}
    Let $\tG\to G\to K$ be a generic tower of harmonic double covers, and let $\tP\to P\to K$ be the bigonal construction.
    
    \begin{enumerate}
        \item The tower $\tP\to P\to K$ is also generic.
        \item Points of $K$ that are dilated with respect to $G \to K$ are in 1:1-correspondence with points of $P$ that are dilated with respect to $\tilde P \to P$, and the same is true for $\tilde G \to G$ and $P \to K$.
        \item The bigonal construction applied to $\tP\to P\to K$ reproduces the original tower.

    \end{enumerate}
    
\end{proposition}

\begin{proof} 
    The first two statements follow directly from the type classification shown on Figure \ref{fig:bigonaltypes}. For the last part, let $\tilde G' \to G' \to K$ be the bigonal construction of $\tilde P \to P \to K$. We claim that there is a canonical bijection $\tG\to \tG'$ that is equivariant with respect to the double covers.
    
    Let $x\in K$ be a type IV point, with preimages $y_1,y_2\in G$ and $\tilde y_1^{\pm},\tilde y_2^{\pm}\in \tG$. The points of $\tP$ over $x$ are the linear combinations $\ty_1^++\ty_2^+$ and $\ty_1^-+\ty_2^-$ (mapping to one point of $P$ over $x$) and $\ty_1^++\ty_2^-$ and $\ty_1^-+\ty_2^+$ (mapping to the other point). The points of $\tG'$ over $x$ are certain linear combinations of these linear combinations, namely
    \[
    2\ty_1^++\ty_2^++\ty_2^-,\qquad 2\ty_1^-+\ty_2^++\ty_2^-,\qquad \ty^+_1+\ty^-_1+2\ty^+_2,\qquad \ty^-_1+\ty^-_1+2\ty^-_2,
    \]
    with the first two and the last two mapping to the same point of $G'$. Hence the map
    \[
    \ty^{\pm}_i\longmapsto 2\ty^{\pm}_i+\ty^+_{3-i}+\ty^-_{3-i}
    \]
    is a canonical equivariant bijection of the fibers of $\tG$ and $\tG'$ over $x$. To establish the bijection of the fibers over a point $x\in K$ of types I to III, we use the same formula and set $\ty_1^{\pm}=\ty_2^{\pm}$ (for type II), $\ty_2^+=\ty_2^-$ (for type III), or both (for type I). Choosing vertex-trivial labelings for the double cover $\tG\to G$, we see that the bijection $\tG\to \tG'$ commutes with the root map, while choosing edge-trivial labelings shows that it commutes with the involution. Hence $\tG\to \tG'$ is an isomorphism of graphs. \qedhere

%    Recall that the points 

%    Again it is obvious from Figure \ref{fig:bigonaltypes} that the fibers of $\tilde G \to G \to K$ and $\tilde G' \to G' \to K$ are the same over every $x \in K$. To complete the proof it remains to check that this identification commutes with the root map. This can be done case by case: all possible half-edge local pictures of the bigonal construction are shown in Figure \ref{tab:bigonal_construction}. The computation is left to the avid reader.
\end{proof}

We now restrict our attention to generic towers $\tG\to G\to K$, where the graph $K$ is a tree, so that the double cover $f:G\to K$ is a \emph{hyperelliptic graph}. 
To state the next proposition we introduce the following notation. The set of dilated edges and vertices form the~\emph{dilation subgraph} $G_{\dil}\subseteq G$. The \emph{dilation index} of the double cover $\pi:\tilde G \to G$ is
\[
d = d(\tilde G / G)=\begin{cases} \mbox{number of connected components of } G_{\dil},& \mbox{if $\pi$ is dilated,}\\
1, & \mbox{if $\pi$ is free.}
\end{cases}
\]

\begin{proposition} 
    \label{prop:bigonalconnected} 
    Let $\tG\to G\to K$ be a generic tower of harmonic double covers, where $\tG$ is connected and $K$ is a tree, and let $\tP\to P\to K$ be the bigonal construction. Then $\tP$ is connected if and only if the double cover $\pi:\tG\to G$ is not free. Furthermore, in this case, the dilation indices $d$ and $d'$ of the double covers $\pi:\tG\to G$ and $\pi':\tP\to P$ and the genera of the four graphs are related as follows:
    \begin{equation}
        d+d'-2=g(\tG)-g(G)=g(\tP)-g(P).
        \label{eq:bigonalgenera}
    \end{equation}
\end{proposition}

\begin{proof} 
    If the double cover $\pi$ is free, then so is the orientation double cover $\tK\to K$, which is then trivial since $K$ is a tree. Therefore $P=\tK$ is disconnected, and hence so is $\tP$.
    
    Conversely, suppose that $\pi$ is not free. If $K$ contains a point $x$ of type I with respect to the tower $\tG\to G\to K$, then $x$ has type I with respect to the tower $\tP\to P\to K$ as well, and therefore $\tP$ is connected. If $K$ contains no points of type I, then it must contain a point of type II (otherwise $G\to K$ is free and hence disconnected) and a point of type III (otherwise $\pi$ is free). The bigonal construction exchanges types II and III, hence the tower $\tP\to P\to K$ also contains both type II and III points. Therefore both $\tP\to P$ and $P\to K$ are dilated, so $\tP$ is connected.

    Let $e_i$ and $v_i$ denote respectively the number of edges and vertices of $K$ having type $i=\mathrm{I},\ldots,\mathrm{IV}$ with respect to the tower $\tG\to G \to K$. We claim that the dilation subgraphs $G_{\dil}\subset G$ and $P_{\dil}\subset P$ are unions of trees. Indeed, if $\gamma\subset G_{\dil}$ is a nontrivial simple cycle, then the restriction of the map $G\to K$ folds $\gamma$ in two because $K$ is simply connected, and hence the image of $\gamma$ in $K$ contains points of type V. Hence each connected component of $G_{\dil}$, and similarly $P_{\dil}$, is a tree. It follows that
    \begin{align*}
        d&=|V(G_{\dil})|-|E(G_{\dil})|=v_{\mathrm{I}}+v_{\mathrm{III}}-e_{\mathrm{I}}-e_{\mathrm{III}},\\
        d'&=|V(P_{\dil})|-|E(P_{\dil})|=v_{\mathrm{I}}+v_{\mathrm{II}}-e_{\mathrm{I}}-e_{\mathrm{II}}.
    \end{align*}
    On the other hand, looking at Figure~\ref{fig:bigonaltypes}, we see that
    \[
    |V(\tG)|-|V(G)|=|V(\tP)|-|V(P)|=v_{\mathrm{II}}+v_{\mathrm{III}}+2v_{\mathrm{IV}},
    \]
    and similarly
    \[
    |E(\tG)|-|E(G)|=|E(\tP)|-|E(P)|=e_{\mathrm{II}}+e_{\mathrm{III}}+2e_{\mathrm{IV}}.
    \]
    Using $\sum e_i-\sum v_i+1=g(K)=0$, we obtain Equation~\eqref{eq:bigonalgenera}. \qedhere
\end{proof}

\begin{figure}[htb]
        \centering
        \begin{subfigure}[b]{0.4\textwidth}
            \centering
            \begin{tikzcd}[column sep = tiny]
    \tilde G \arrow[d] & 
    \begin{tikzpicture}[baseline = (O.base)]
        \node (O) at (0, -0.25) {};
        \drawTildeGammaBigonal
    \end{tikzpicture} 
    \arrow[d] \\
    G \arrow[d] &
    \begin{tikzpicture}
    	\draw (0.5, 0) ellipse (0.5 and 0.2);
    	\path[line width = 2*\edgewidth, draw] (1,0) -- + (1,0);
    	\draw (2.5, 0) ellipse (0.5 and 0.2);
    	\path[line width = 2*\edgewidth, draw] (3,0) -- + (1,0);
    	\draw (4.5, 0) ellipse (0.5 and 0.2);
    	\foreach \i in {0,1,2,3,4,5} {\vertex[2]{\i,0}}
    \end{tikzpicture}
    \arrow[d] \\
    K & 
    \begin{tikzpicture}
        \path[draw] (0,0) -- + (5,0);
    	\foreach \i in {0,1,2,3,4,5} {\vertex{\i,0}}
    \end{tikzpicture}
\end{tikzcd}
        \end{subfigure}
        \qquad \qquad
        \begin{subfigure}[b]{0.4\textwidth}
            \centering
            \begin{tikzcd}[column sep = tiny]
    \tilde P \arrow[d] &
    \begin{tikzpicture}[baseline = (O.base)]
        \node (O) at (0, 0.5) {};
        \path[draw] (1,0) -- + (-1, 0.2);
    	\path[draw] (1,0) -- + (-1, -0.2);
    	\path[line width = 2*\edgewidth, draw] (1,0) -- + (1,0);
    	\draw (2.5, 0) ellipse (0.5 and 0.2);
    	\path[line width = 2*\edgewidth, draw] (3,0) -- + (1,0);
    	\path[draw] (4,0) -- + (1, 0.2);
    	\path[draw] (4,0) -- + (1, -0.2);
    	
    	\path[draw] (0,1) -- ++ (1,0.2) -- ++ (3,0) -- ++ (1,-0.2) -- ++ (-1, -0.2) -- ++ (-3,0) -- ++ (-1, 0.2);
    	
    	\foreach \i in {1,2,3,4} {
    	    \vertex[2]{\i, 0}
    	    \vertex{\i, 0.8}
    	    \vertex{\i, 1.2}
    	    }
    	\vertex{0, -0.2}
    	\vertex{0, 0.2}
    	\vertex[2]{0, 1}
    	\vertex{5, -0.2}
    	\vertex{5, 0.2}
    	\vertex[2]{5, 1}
    \end{tikzpicture}
    \arrow[d] \\
    P \arrow[d] &
    \begin{tikzpicture}[baseline = (O.base)]
        \node (O) at (0, 0.15) {};
        \path[draw] (0,0) -- +(5, 0);
    	\path[draw] (0,0.3) -- +(5, 0);
    	\foreach \i in {0,1,2,3,4,5} {
    	    \vertex{\i, 0}
    	    \vertex{\i, 0.3}
    	    }
    \end{tikzpicture}
    \arrow[d] \\
    K &
    \begin{tikzpicture}
        \path[draw] (0,0) -- + (5,0);
        \foreach \i in {0,1,2,3,4,5} {\vertex{\i,0}}
    \end{tikzpicture}
    \\
\end{tikzcd}
        \end{subfigure}
        \caption{Example of the bigonal construction with thickness indicating dilation with respect to the base tree $K$. The involution on $\tilde G$ is reflection along the horizontal axis and similarly for each of the two components of $\tilde P$.}
        \label{fig:example_bigonal_construction}
    \end{figure}
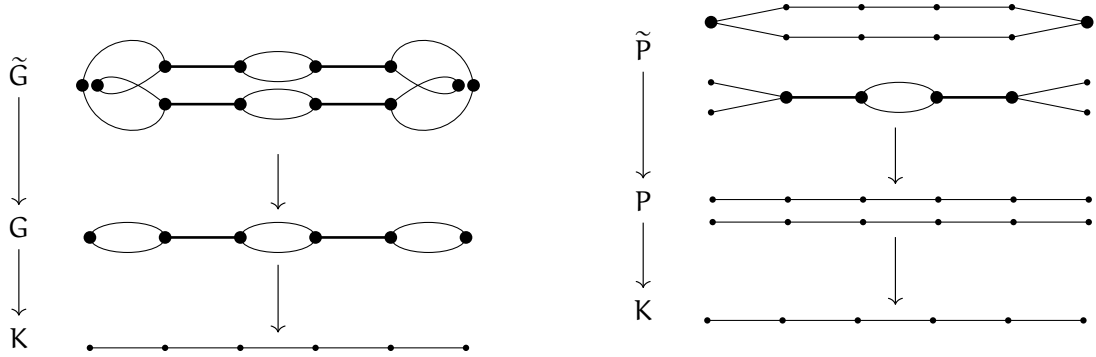
    
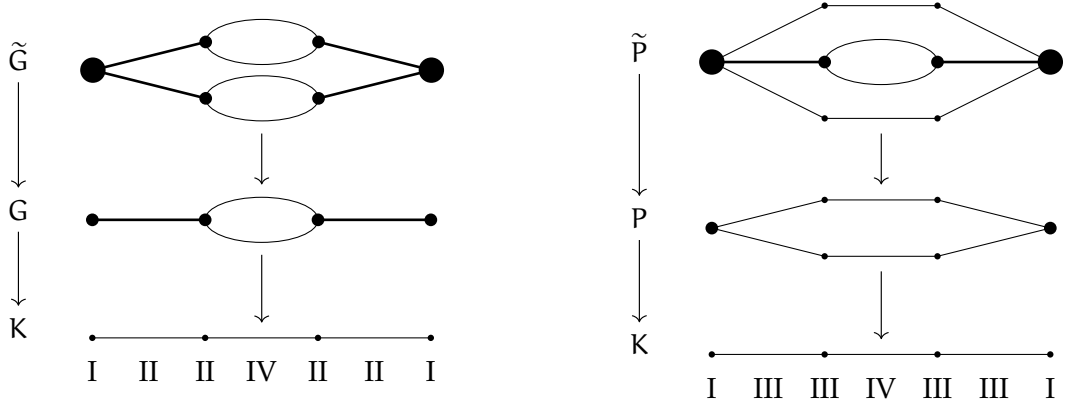
\begin{figure}[htb]
    \centering
    \begin{subfigure}[b]{0.4\textwidth}
            \centering
            $\begin{tikzcd}[column sep = tiny]
\tilde G
\arrow[d]
& 
\begin{tikzpicture}[x=15mm, y=15mm, baseline = (O.base)]
    \drawExampleBigonal
    \node (O) at (0, -0.25) {};
\end{tikzpicture} 
\arrow[d]
\\
G
\arrow[d]
&
\begin{tikzpicture}[x=15mm, y=15mm]
    \vertex[2]{0,0}
    \path[line width = 2*\edgewidth, draw] (0,0) -- ++ (1, 0);
    \vertex[2]{1,0}
    \draw (1.5, 0) ellipse (0.5 and 0.2);
    \vertex[2]{2,0}
    \path[line width = 2*\edgewidth, draw] (2, 0) -- ++ (1, 0);
    \vertex[2]{3, 0}
\end{tikzpicture} 
\arrow[d]
\\
K
&
\begin{tikzpicture}[x=15mm, y=15mm]
    \path[draw] (0,0) -- ++ (3, 0);
    \vertex{0,0}
    \draw (0, -0.1) node[anchor = north] {I};
    \draw (0.5, -0.1) node[anchor = north] {II};
    \vertex{1,0}
    \draw (1, -0.1) node[anchor = north] {II};
    \draw (1.5, -0.1) node[anchor = north] {IV};
    \vertex{2,0}
    \draw (2, -0.1) node[anchor = north] {II};
    \draw (2.5, -0.1) node[anchor = north] {II};
    \vertex{3,0}
    \draw (3, -0.1) node[anchor = north] {I};
\end{tikzpicture}
\end{tikzcd}$
        \end{subfigure}
        \qquad \qquad
        \begin{subfigure}[b]{0.4\textwidth}
            \centering
            $\begin{tikzcd}[column sep = tiny]
\tilde P
\arrow[d]
& 
\begin{tikzpicture}[x=15mm, y=15mm]
    \drawExampleBigonalOutput
\end{tikzpicture} 
\arrow[d]
\\
P
\arrow[d]
&
\begin{tikzpicture}[x=15mm, y=15mm]
    \path[draw] (0,0) -- ++ (1, 0.25) -- ++ (1, 0) -- ++ (1, -0.25) -- ++ (-1, -0.25) -- ++ (-1, 0) -- ++ (-1, 0.25);
    \vertex[2]{0,0}
    \vertex{1,0.25}
    \vertex{2,0.25}
    \vertex{1,-0.25}
    \vertex{2,-0.25}
    \vertex[2]{3, 0}
\end{tikzpicture} 
\arrow[d]
\\
K
&
\begin{tikzpicture}[x=15mm, y=15mm]
    \path[draw] (0,0) -- ++ (3, 0);
    \vertex{0,0}
    \draw (0, -0.1) node[anchor = north] {I};
    \draw (0.5, -0.1) node[anchor = north] {III};
    \vertex{1,0}
    \draw (1, -0.1) node[anchor = north] {III};
    \draw (1.5, -0.1) node[anchor = north] {IV};
    \vertex{2,0}
    \draw (2, -0.1) node[anchor = north] {III};
    \draw (2.5, -0.1) node[anchor = north] {III};
    \vertex{3,0}
    \draw (3, -0.1) node[anchor = north] {I};
\end{tikzpicture}
\end{tikzcd}$
        \end{subfigure}
    \caption{Example of the bigonal construction with connected input and output. Thickness corresponds to dilation and the involutions are given by reflection along the horizontal axis. For each point of $K$ we indicate its type.}
    \label{fig:bigonal_example_2}
\end{figure}

\begin{example} 
    Figure~\ref{fig:example_bigonal_construction} shows an example of the bigonal construction. Note that $\tilde G \to G$ is a free double cover and therefore $\tilde P$ is necessarily disconnected. Nevertheless, one can check that applying the bigonal construction to the tower $\tilde P \to P \to K$ reproduces the input tower. Let us modify this example by contracting the extremal edges of $K$ and everything lying above them. The result is shown in Figure~\ref{fig:bigonal_example_2}. This time the input (and hence the output as well) contains a vertex of type I. In particular, input and output are connected. Again one can check that the bigonal construction applied to $\tilde P \to P \to K$ reproduces the original tower. 
\end{example}

\subsection{The tropical trigonal and Recillas construction.} \label{sec:Recillas}

Let $\tG\to G\to K$ be an orientable tower of covers of a graph $K$ of degrees $2$ and $3$, respectively. In particular, the double cover $\pi : \tG \to G$ is free. The trigonal construction associates to $\tG\to G\to K$ a degree 8 harmonic morphism $\tP\to K$ that splits as a tower $\tP\to P\to K$, where $p:P\to K$ has degree 4 and $\tP\to P$ is the trivial double cover. The \emph{Recillas construction} inverts this correspondence, by associating to a degree $4$ cover $P\to K$ (satisfying certain restrictions on local degrees) an orientable tower $\tG\to G\to K$ of covers of degrees $2$ and $3$.

Similarly to the bigonal case, we classify the points of $K$ according to the structure of the fibers. We recall that the points of $\tP$ over a given point $x\in K$ have the form $D=\sum (a^+_i\ty^+_i+a^-_i\ty^-_i)$, where the $y_i$ are the preimages of $x$ in $G$, and the $a^{\pm}_i$ satisfy $a_i^{\pm}\geq 0$ and $a^+_i+a^-_i=d_f(y_i)$ for all $i$. Since the tower $\tG\to G\to K$ is orientable, we can represent points of $P$ as those linear combinations $D=\sum (a^+_i\ty^+_i+a^-_i\ty^-_i)$ for which the quantity $\sum a^+_i$ has a fixed parity, and we choose $\sum a^+_i\equiv 1$ mod $2$.

\begin{definition} \label{def:trigonal_types} 
    Let $\tG\to G\to K$ be an orientable tower consisting of a free double cover $\pi:\tG\to G$ and a degree 3 harmonic morphism $f:G\to K$, and let $p:P\to K$ be the associated harmonic morphism of degree $4$. A point $x\in K$ is said to be

    \begin{enumerate}

        \item \emph{Type I} if it has a unique preimage $y_1$ in $G$ with $d_f(y_1)=3$. The corresponding points of $P$ are $3\ty_1^+$ and $\ty_
        1^++2\ty_1^-$, at which $p$ has degrees $1$ and $3$, respectively.
        
        \item \emph{Type II} if it has two preimages in $G$, namely $y_1$ with $d_f(y_1)=1$ and $y_2$ with $d_f(y_2)=2$. The corresponding points of $P$ are $\ty_1^++2\ty_2^+$, $\ty_1^++2\ty_2^-$, and $\ty_1^-+\ty_2^++\ty_2^-$, at which $p$ has degrees $1$, $1$, and $2$, respectively.        
        
        \item \emph{Type III} if it has three preimages $y_1$, $y_2$, and $y_3$. The corresponding points of $P$ over $x$ are $\ty_1^++\ty_2^++\ty_3^+$, $\ty_1^++\ty_2^-+\ty_3^-$, $\ty_1^-+\ty_2^++\ty_3^-$, and $\ty_1^-+\ty_2^-+\ty_3^+$. The local degrees of $f$ and $p$ are all equal to one.

    \end{enumerate}
\end{definition}

We note that the type of a point $x\in K$ is the number of its preimages in $G$, and that a half-edge may be rooted at a vertex of equal or lower type. There are thus six possible pairings of a half-edge and a vertex, and the local structure of the tower and the degree four map are shown on Figure~\ref{tab:trigonal_construction}. We also observe that the degree $4$ harmonic morphism $p:P\to K$ has the property that the fibers of $p$ cannot have degree profiles $(4)$ or $(2,2)$. Note that, in the algebraic setting, an identical restriction is imposed on the \emph{ramification} profile of the degree four map. We give a corresponding definition:

\begin{figure}[htb]
    \centering
    
\begin{tabular}{b{1.5cm}| *{6}{b{1.5cm}}}
    \begin{center}
        \begin{tikzcd}
            \tilde G \arrow[d, "\pi"] \\
            G \arrow[d, "f"] \\
            K
        \end{tikzcd}
    \end{center}
     & 
     \tower{
        \halfedge{0}
        \halfedge{0.2}
        \halfedge{0.7}
        \halfedge{0.9}
        \halfedge{1.4}
        \halfedge{1.6}
     }{
        \halfedge{0}
        \halfedge{0.3}
        \halfedge{0.6}
     }
     & 
     \tower{
        \halfedge{0}
        \halfedge{0.2}
        \twohalfedges{0.7}
        \twohalfedges{1.2}
     }{
        \halfedge{0}
        \twohalfedges{0.4}
     }
     & 
     \tower{
        \halfedge{0}
        \halfedge{0.2}
        \halfedge[2]{0.7}
        \halfedge[2]{0.9}
     }{
        \halfedge{0}
        \halfedge[2]{0.3}
     }
     & 
     \tower{
        \twohalfedges[2]{0}
        \twohalfedges[2]{0.5}
     }{
        \twohalfedges[2]{0}
     }
     &  
     \tower{
        \threehalfedges{0}
        \threehalfedges{0.6}
     }{
        \threehalfedges{0}
     }
     & 
     \tower{
        \halfedge[3]{0}
        \halfedge[3]{0.2}
     }{
        \halfedge[3]{0}
     }
     \\
     \hline
     \begin{center}
        \begin{tikzcd}
            P \arrow[d, "p"] \\
            K
        \end{tikzcd}
     \end{center}
     & 
     \tetragonal{\halfedge{0} \halfedge{0.3} \halfedge{0.6} \halfedge{0.9}}
     & 
     \tetragonal{\halfedge{0} \twohalfedges{0.5} \halfedge{1}}
     & 
     \tetragonal{\halfedge{0} \halfedge[2]{0.3} \halfedge{0.6}}
     & 
     \tetragonal{\halfedge{0} \twohalfedges[2]{0.5}}
     & 
     \tetragonal{\halfedge{0} \threehalfedges{0.6}}
     & 
     \tetragonal{\halfedge{0} \halfedge[3]{0.3}}
\end{tabular}

    \caption{Overview of the trigonal construction and its inverse, locally over a half-edge of $K$. Thickness of edges and vertices corresponds to dilation.}
    \label{tab:trigonal_construction}
\end{figure}
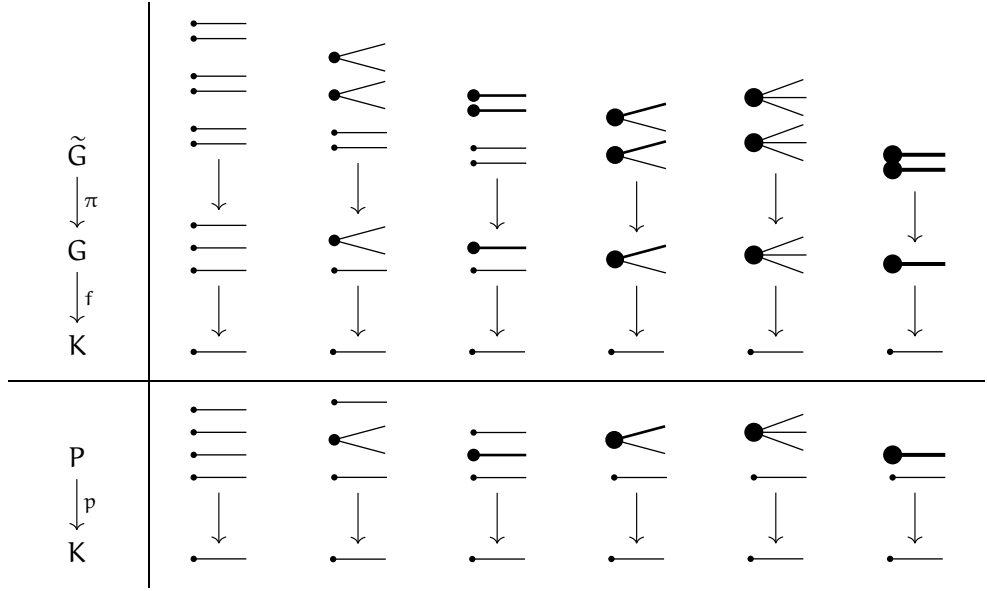

\begin{definition} 
    \label{def:generic_tetragonal}
    A degree $4$ harmonic morphism $p:P\to K$ is called \emph{generic} if every point $x\in K$ has a preimage in $P$ at which the local degree of $p$ is equal to one. Given a generic $p:P\to K$, a point $x$ of $K$ is said to be of \emph{type I, II, or III} if the degree profile of the fiber is $(3,1)$, $(2,1,1)$, or $(1,1,1,1)$, respectively.

\end{definition}

We now invert the tropical trigonal construction for a generic degree $4$ harmonic morphism.

\begin{definition}
\label{def:Recillas_construction}

Let $p:P\to K$ be a generic degree $4$ harmonic morphism. The \emph{Recillas construction} associated to $p:P\to K$ is a tower consisting of a free double cover $\pi:\tG\to G$ and a harmonic morphism $f:G\to K$ of degree $3$, defined as follows. The points of $\tG$ over a given point $x \in K$ are unordered pairs of points of $P$ lying in the same fiber of $p$:
\[
\tG=\left\{D_1+D_2\in \ZZ^P : D_1, D_2 \in P \text{ and } D_1+D_2\leq p^*(x)\mbox{ for some (necessarily unique) }x\in K\right\}.
\]
The vertices and half-edges of $\tG$ lie over the vertices and half-edges of $K$, respectively, and the root and involution maps are induced from $P$. It is clear that $\tG$ is a graph and $\tf:\tG\to K$, $D_1+D_2\mapsto x$ is a morphism of graphs. 

There is a natural involution on $\tG$ given by 
\[ D_1+D_2 \longmapsto p^*(x) - D_1-D_2, \]
and the assumption that $p$ is generic implies (indeed, is equivalent to assuming) that this involution is fixed-point-free. Hence $\tf:\tG\to K$ factors as a free quotient double cover $\pi:\tG\to G$ and an induced map $f:G\to K$.

Finally, we define the local degrees of $f$ and $\tf$. Let $x\in K$. If $x$ is a type III point, then $p^*(x)=D_1+D_2+D_3+D_4$ for distinct $D_i$ and we set $d_{\tf}(D_i+D_j)=1$ for all pairs. If $x$ is a type II point, then $p^*(x)=D_1+D_2+2D_3$ and we set
\[
d_{\tf}(D_1+D_2)=d_{\tf}(2D_3)=1,\qquad d_{\tf}(D_1+D_3)=d_{\tf}(D_2+D_3)=2.
\]
If $x$ is a type I point, then $p^*(x)=D_1+3D_2$ and we set
\[
d_{\tf}(D_1+D_2)=d_{\tf}(2D_2)=3.
\]
For any $\ty\in \tG$ we set $d_f(\pi(\ty))=d_{\tf}(\ty)$, and it is elementary to verify that $f$ is a harmonic morphism of degree 3.

\end{definition}

\begin{remark} 
    Generalizing Definition~\ref{def:Recillas_construction}, we can naturally associate to a harmonic morphism $P\to K$ of degree $n$ a degree $\binom{n}{k}$ harmonic morphism $\tG\to K$, for any $k\leq n$.
\end{remark}

We now assume that $K$ is a tree, and refer to harmonic morphisms $G\to K$ and $P\to K$ of degrees 3 and 4 as respectively \emph{trigonal} and \emph{tetragonal} graphs (we note that this condition is stronger than carrying a $g^1_3$ or a $g^1_4$ in the sense of Baker and Norine~\cite{baker2007riemann}, see~\cite{ABBRII}). 
The following establishes the first part of Theorem \ref{thm:tropical_Recillas_theorem} in the case of graphs.

\begin{proposition} 
    \label{prop:trigonal_properties}
    Let $K$ be a tree. The trigonal construction and the Recillas construction establish a bijection between free double covers of trigonal graphs $\tG\to G\to K$ and generic tetragonal graphs $P\to K$. The graph $\tG$ is connected if and only if $P$ is connected, in which case
    \begin{equation}
        \label{eq:trigonalgenera}
    g(P)=g(G)-1.
    \end{equation}
\end{proposition}

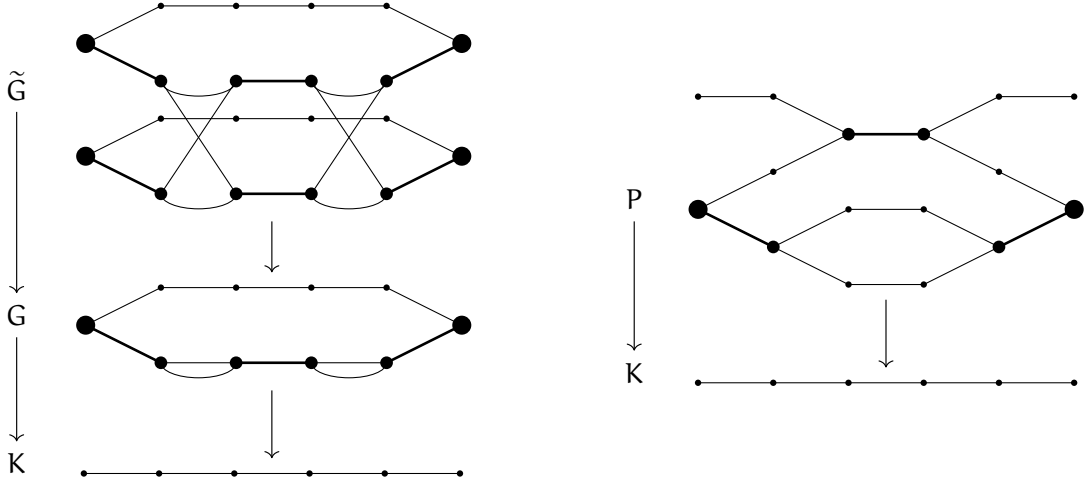
\begin{figure}[htb]
    \centering
    \begin{subfigure}{0.4\textwidth}
        \centering
        \begin{tikzcd}[column sep = tiny]
    		\tilde{G} \arrow[d] & 
    		\begin{tikzpicture}[baseline = (O.base)]
                \node (O) at (0, 1.25) {};
    			\drawGammaTilde
    		\end{tikzpicture} \arrow[d] \\
    		G \arrow[d] &
    		\begin{tikzpicture}[baseline = (O.base)]
                \node (O) at (0, 0.5) {};
    			\path[draw] (0,0.5) -- ++ (1, -0.5) -- ++ (3,0) -- ++ (1, 0.5) -- ++ (-1, 0.5) -- ++ (-3, 0) -- ++ (-1, -0.5);
    			\path[line width = 2*\edgewidth, draw] (0, 0.5) -- ++ (1, -0.5) ++ (1, 0) -- ++ (1, 0) ++ (1, 0) -- ++ (1, 0.5);
    			\draw (1,0) arc (180:360:0.5 and 0.2);
    			\draw (3,0) arc (180:360:0.5 and 0.2);
    			\vertex[3]{0, 0.5}
    			\vertex[3]{5, 0.5}
    			\foreach \x in {1,2,3,4} {\vertex[2]{\x,0} \vertex{\x,1}}
    		\end{tikzpicture} \arrow[d] \\
    		K & 
    		\begin{tikzpicture}
    			\path[draw] (0,0) -- + (5,0);
    			\foreach \x in {0,1,2,3,4,5} {\vertex{\x,0}}
    		\end{tikzpicture}
    	\end{tikzcd}
    \end{subfigure}
    \qquad \qquad
    \begin{subfigure}{0.4\textwidth}
    	\begin{tikzcd}[column sep = tiny]
    		P \arrow[d] & 
    		\begin{tikzpicture}
    			\drawPi			
    		\end{tikzpicture} \arrow[d] \\
    		K &
    		\begin{tikzpicture}
    			\path[draw] (0,0) -- + (5,0);
    			\foreach \x in {0,1,2,3,4,5} {\vertex{\x,0}}
    		\end{tikzpicture}
    	\end{tikzcd}
    \end{subfigure}
    \caption{Example of a tower and a tetragonal graph corresponding to each other under the tropical trigonal construction. Thickness indicates dilation factors.}
    \label{fig:example_trigonal}
\end{figure}

\begin{example} 

    Consider the tower $\tilde G \to G \to K$ on the left of Figure~\ref{fig:example_trigonal}. Applying the tropical trigonal construction to it produces the generic tetragonal graph on the right. Conversely, applying the tropical Recillas construction to $P \to K$ recovers the original tower. We will verify by hand in Example~\ref{ex:trigonal_construction_computation} that the Prym variety $\Prym(\tG/G)$ and the Jacobian variety $\Jac(P)$ (which are defined after $K$ has been equipped with arbitrary edge lengths) are isomorphic.
\end{example}

\begin{proof}[Proof of Prop.~\ref{prop:trigonal_properties}] 

Let $\pi:\tG\to G$ be a free double cover of a trigonal graph $f:G\to K$. The orientation double cover $\tK\to K$ is free and hence trivial, hence the tower $\tG\to G\to K$ is orientable. We have already seen that the associated tetragonal graph $p:P\to K$ is generic, and let $\tG'\to G'\to K$ denote its Recillas construction.

As in the bigonal case, we construct a canonical equivariant bijection $\tG\to \tG'$. Let $x\in K$ be a type III point, with preimages $y_1,y_2,y_3$ in $G$ and $\ty_1^{\pm},\ty_2^{\pm},\ty_3^{\pm}$ in $\tG$. The points of $P$ over $x$ are the four linear combinations
\[
\ty_1^++\ty_2^++\ty_3^+,\qquad \ty_1^++\ty_2^-+\ty_3^-,\qquad\ty_1^-+\ty_2^++\ty_3^-,\qquad \ty_1^-+\ty_2^-+\ty_3^+,
\]
and the points of $\tG'$ over $x$ are the six pairwise sums of these combinations. Hence we see that
\[
\ty^{\pm}_i\longmapsto 2\ty^{\pm}_i+\sum_{j\neq i}(\ty^+_j+\ty^-_j)
\]
is the required bijection, which is equivariant since it preserves signs. For points of types II and I, we use the same formula specialized to $\ty^{\pm}_2=\ty^{\pm}_3$ and $\ty^{\pm}_1=\ty^{\pm}_2=\ty^{\pm}_3$, respectively. As in the proof of Proposition~\ref{prop:bigonal_properties}, choosing either vertex-trivial or edge-trivial labelings for $\tG\to G$ shows that the bijection $\tG\to \tG'$ is an isomorphism of graphs.

Conversely, let $P'\to K$ be the trigonal construction of the Recillas construction $\tG\to G\to K$ of a generic tetragonal graph $P\to K$. Then the map
\[
D_i\longmapsto 3D_i+\sum_{j\neq i}D_j
\]
establishes a bijection between the fibers of $P$ and $P'$ over $x$, preserves root and involution maps, and is hence a graph isomorphism $P\to P'$. Hence the trigonal and Recillas constructions are inverses.

    To show the second part of the claim, suppose that the tower $\tG \to G \to K$ and the tetragonal graph $p:P\to K$ correspond to one another under the trigonal and Recillas construction. Assume that $P=P_1\sqcup P_2$ is disconnected, and further assume without loss of generality that $\deg (p|_{P_1})\geq 2$. Viewing points of $\tG$ as 
    linear combinations $D_1+D_2\in \ZZ^P$ (with possibly $D_1=D_2$), we define the function
    \[ \deg_1 : \tG \longrightarrow \ZZ, \qquad  D_1 + D_2 \longmapsto \big| \{ i : D_i \in P_1 \} \big|. \]
    Each fiber of $p$ contains at least two points of $P_1$ (counted with multiplicity) and at least one point of $P_2$, hence $\deg_1$ takes values $1$ and $2$. On the other hand, $\deg_1$ is preserved by the root and involution maps, and is therefore a locally constant function on $\tG$. Hence $\tG$ is disconnected. 
    %\Felix{Note that this can be viewed as a continuity argument for a suitable topology...}

    Conversely, assume that $\tG$ is disconnected. We distinguish two cases. If $G$ is connected, then the double cover $\tG\to G$ is trivial, i.e. $\tG = \tG^+ \sqcup \tG^-$ and $\tG^+ \cong \tG^- \cong G$, and we can globally label the preimages of $y\in G$ in $\tG$ as $\ty^+\in \tG^+$ and $\ty^-\in\tG^-$. Viewing points of $\tP$ as sections of the double cover $\tG\to G$, the map
    \[ \deg_2 : \tP \longrightarrow \ZZ, \qquad \sum_{i=1}^3 \left(a_i^+ \ty_i^++a_i^-\ty_i^-\right) \longmapsto \sum_{i=1}^3 a_i^+ \]
    defines a locally constant function on $\tP$ with image $\{0, 1, 2, 3\}$. By definition, two points $D,D'\in\tP$ over the same point of $K$ map to the same point in $P$ if and only if $\deg_2(D) \equiv \deg_2(D')$ mod $2$. Hence $\deg_2$ induces a locally constant, surjective $\ZZ/2\ZZ$-valued function on the quotient $P$ of $\tP$, and therefore $P$ is disconnected.

    For the second case we assume that $G$ is disconnected as well. Since $f : G \to K$ is harmonic of degree 3, at least one of the connected components of $G$ maps with degree 1 to $K$ and is therefore isomorphic to $K$. We write $G = G_1 \sqcup K$ and correspondingly $\tG = \tG_1 \sqcup K \sqcup K$ for the double cover. Any point $D\in \tP $
    has the form $D = D_1 + y$, where $D_1$ is a point of the bigonal construction of $\tG_1 \to G_1 \to K$ and $y$ is a point in one of the two copies of $K \subseteq \tG$, with $D_1$ and $y$ lying over the same point of the target tree. In other words, $\tP$ consists of two copies of the bigonal construction of $\tG_1 \to G_1 \to K$ and these are exchanged by the involution on $\tP$. Hence $P$ itself is precisely the bigonal construction of $\tG_1 \to G_1 \to K$. By Proposition~\ref{prop:bigonal_properties}, we see that it is disconnected, because $\tG_1\to G_1$ is free.
    
    To complete the proof, we determine the relationship between the genera. Looking at the local structure in Figure \ref{tab:trigonal_construction}, we see that $p^{-1}(x)$ has one more element than $f^{-1}(x)$ for each point $x\in K$. It follows that 
    \[
        g(P)=|E(P)|-|V(P)|+1=|E(G)|+|E(K)|-|V(G)|-|V(K)|+1=g(G)-1,
    \]
    because $K$ is a tree and hence $|E(K)|-|V(K)|=-1$. This completes the proof. \qedhere
\end{proof}

\begin{remark} \label{rem:dilationcollapse}
    The trigonal construction can be applied to towers $\tG\to G\to K$ where $\tG\to G$ is not free. The result is a tower $\tP\to P\to K$ where, depending on the degree profiles of the first tower, $\tP\to P$ may be dilated, $P\to K$ non-generic, or both. Conversely, the Recillas construction can be extended to non-generic tetragonal graphs $P\to K$ to produce towers $\tG\to G\to K$ where $\tG\to G$ is not free (for example, a point $x$ of $K$ with degree profile $(2,2)$ with respect to the tetragonal map has two preimages in $G$, one having degree $2$ and a single preimage in $\tG$, the other having degree $1$ and two preimages). However, applying the trigonal construction to the resulting tower $\tG\to G\to K$ does not reproduce the original tetragonal graph $P\to K$. Hence the bijective correspondence fails for these generalizations.
    
\end{remark}

\subsection{The tetragonal construction} \label{sec:tetragonal}
In this section, we briefly summarize the harmonic tetragonal construction, which we plan to study in detail in a future paper.
Let $\tG\to G\to K$ be a tower of harmonic morphisms of degrees $2$ and $4$, and let $\tP\to K$ be the outcome of the tetragonal construction. In general, the graph $\tP$ does not split (if $\tG\to G\to K$ is not orientable), and the involution $\iota:\tP\to \tP$ may have fixed points. However, imposing some natural restrictions on the dilation produces an outcome that exactly corresponds to the free case, and mirrors the algebraic construction.

\begin{proposition} 
    \label{prop:tetragonal_properties}
    Let $K$ be a tree, and let $\tG\to G\to K$ be a free double cover of a generic (in the sense of Definition \ref{def:generic_tetragonal}) tetragonal graph $G\to K$. The tetragonal construction applied to $\tG\to G\to K$ splits as a disjoint union of $\tG_i\to G_i\to K$ for $i=1,2$, where each tower is a free double cover of a generic tetragonal graph.
    
\end{proposition}

\begin{proof} 
    If $\tG\to G$ is free then so is the orientation double cover $\tK\to K$, which is then trivial since $K$ is a tree. Hence $\tP\to K$ splits as a disjoint union of morphisms $\tP_1\to K$ and $\tP_2\to K$ of degree eight. Since $G\to K$ is generic, each point $x\in K$ has a preimage in $G$ at which the tetragonal map has odd degree (specifically, equal to one) and which has two preimages in $\tG$ (since $\tG\to G$ is free). Hence the sign involution $\iota:\tP\to \tP$ acts without fixed points. Since the involution restricts to each connected component, we take quotients and obtain two towers $\tP_i\to P_i\to K$ for $i=1,2$. It is then a direct verification to show that if a point $x \in K$ has type $I$, $II$, or $III$ in the sense of Definition~\ref{def:generic_tetragonal} with respect to the original tower, then $x$ has the same type with respect to each of the two new towers. In particular, the $P_i\to K$ are generic tetragonal graphs. \qedhere
\end{proof}

%%%%%%%%%%%%%%%%%%%%%%%%%%%%%%%%%%%%%%%%%%%%%%%%%%%
\section{Rational polyhedral spaces and tropical homology}
\label{sec:tropical}
%%%%%%%%%%%%%%%%%%%%%%%%%%%%%%%%%%%%%%%%%%%%%%%%%%%

In this section, we review a number of standard notions of tropical geometry. The ambient category containing all tropical objects that we consider in this article is the category of \emph{rational polyhedral spaces}. In particular, tropical curves are purely 1-dimensional rational polyhedral spaces satisfying a smoothness condition. We rephrase the tropical $n$-gonal construction and the tropical Recillas construction for tropical curves (thus justifying the name). Finally, we summarize some basic properties of tropical cycles and tropical homology that will serve as essential tools to prove Theorems~\ref{thm:bigonal_construction} and~\ref{thm:tropical_Recillas_theorem}. Tropical homology was introduced in~\cite{IKMZ}, but our exposition closely follows the sheaf-theoretic approach~\cite{GrossShokrieh_homology}.

\subsection{Rational polyhedral spaces}

A \emph{(rational) polyhedron} in $\RR^n$ is a finite intersection of half-spaces of the form
\[ \big\{x \in \RR^n \bigmid \langle m, x \rangle \leq a \big\} \]
for some $m \in (\ZZ^n)^\ast$ and $a \in \RR$. Consider the partial compactification $\overline{\RR}^n$ of $\RR^n$, where $\overline{\RR} = \RR \cup \{\infty\}$, endowed with the order topology. It is stratified by sets of the form
\[ \overline{\RR}_I^n = \big\{ (x_i)_{i = 1, \ldots, n} \bigmid x_i = \infty \text{ if and only if } i \in I \big\} \]
for $I\subset\{1,\ldots,n\}$. A \emph{polyhedron} in $\overline{\RR}^n$ is the closure of a polyhedron in one of the $\overline{\RR}_I^n$. A \emph{polyhedral subset} $X \subseteq \overline{\RR}^n$ is a finite union of polyhedra.

Let $X \subseteq \overline{\RR}^n$ be a polyhedral subset. An \emph{integral affine linear function} on $X$, or \emph{affine function} for short, is a function $f : X \to \RR$ such that locally at every point of $X$ it is of the form $x \mapsto \langle m, x \rangle + a$ for some $m \in (\ZZ^n)^\ast$ and $a \in \RR$. An affine function is not allowed to take an infinite value, so the local expression $\langle m, x\rangle + a$ is required to satisfy $m_i=0$ near any point $x$ of $X$ with $x_i=\infty$. In other words, at infinite points, affine functions are locally constant in those directions in which infinite coordinates become finite. Affine functions form a sheaf on $X$, denoted $\Aff_X$.

A \emph{rational polyhedral space} $X$ is a second countable Hausdorff topological space together with a sheaf of continuous functions $\Aff_X$ such that for every point $x \in X$ there is an open neighborhood $x \in U \subseteq X$, an open subset $V \subseteq Y$ of a polyhedral set $Y \subseteq \overline{\RR}^n$, and a homeomorphism $\phi : U \to V$ such that pullback of affine functions along $\phi$ is an isomorphism $\phi^{-1} \Aff_V \to \Aff_U$.
If all these polyhedral sets $Y$ can be taken to live in $\RR^n$, then $X$ is called \emph{boundaryless}.
A point $x \in X$ is called \emph{regular} if it has an open neighborhood isomorphic to an open subset of $\RR^n$, where $n$ is the \emph{local dimension} at $x$. The subset of regular points in $X$ is denoted $X^\reg$. A rational polyhedral space $X$ is said to be \emph{purely $n$-dimensional} if each point of $X^\reg$ has local dimension $n$. 

A \emph{morphism of rational polyhedral spaces} is a continuous map $f : X \to Y$ that induces a morphism $f^{-1}\Aff_Y \to \Aff_X$. It is \emph{proper} if the preimage of every compact subset of $Y$ is compact. In particular, if $X$ is compact, then $f$ is proper.

The \emph{cotangent sheaf} $\Omega_X^1$ of a rational polyhedral space $X$ is the quotient of the sheaf $\Aff_X$ by the subsheaf of locally constant functions. 
A morphism $f : X \to Y$ of rational polyhedral spaces induces a morphism of cotangent sheaves
\[ f^{-1}\Omega_Y^1 \longrightarrow \Omega_X^1. \]
For $x \in X$, the dual $T_x^\ZZ X = \Hom_\ZZ(\Omega_{X, x}^1, \ZZ)$ of the stalk of the cotangent sheaf is the \emph{integral tangent space} at $x$. 
Dualizing the morphism of cotangent sheaves gives the differential $d_xf : T_x^\ZZ X \to T_{f(x)}^\ZZ Y$.

\subsection{Tropical curves and harmonic morphisms} Let $\Ga$ be a connected and compact purely 1-dimensional rational polyhedral space. The underlying topological space of $\Ga$ has the combinatorial structure of a connected and finite graph $G$. Furthermore, we can define an edge length function $\ell:E(G)\to (0,\infty]$ by setting $\ell(e)$ to be the smallest positive increment of an affine function along $e$. The pair $(G,\ell)$, which is a \emph{metric graph}, is called a \emph{model} for $\Ga$. We say that $\Ga$ is a \emph{tropical curve} if all edges of infinite length are extremal (in other words, there are finitely many vertices at infinity, and they are univalent). We always implicitly choose a model when talking about a tropical curve, and we note that the model recovers the curve if the latter is smooth, which we always assume (see below). The genus of $\Ga$ is the genus of any graph model. By a point $x$ of $\Ga$ we mean a point in the metric space, which may correspond to either a vertex or an interior point of an edge with respect to a chosen model (an edge $e\in E(G)$ of the model may be viewed as a generic point of the corresponding edge in $\Ga$).

An affine function $f$ on a tropical curve $\Ga$ has a well-defined \emph{slope} along each oriented edge $e$ of $\Ga$. A tropical curve $\Gamma$ is called \emph{smooth} if locally around every finite vertex it is isomorphic to
\begin{equation} \bigcup_{i = 0, \ldots, n} e_i\RR_{\geq 0} \subseteq \RR^{n+1} / (1, \ldots, 1)\RR \label{eq:smoothlocalmodel}\end{equation}
for some $n \geq 1$.
This condition ensures that $\Gamma$ has sufficiently many affine functions, in the following sense: given an affine function defined near a point $x\in\Ga$, the only condition on the slopes along the outgoing edges is that they sum to zero. Furthermore, the univalent vertices of a smooth tropical curve $\Gamma$ are at infinity. We note that an arbitrary metric graph can be augmented to produce a smooth tropical curve in a canonical way, by attaching a compact infinite ray to each finite univalent vertex.

Conversely, let $G$ be a finite graph and let $\ell:E(G)\to (0,\infty]$ be an edge length assignment such that an edge is infinite if and only if it is extremal. We construct a smooth tropical curve $\Ga$ with model $(G,\ell)$ by gluing real intervals $[0, \ell(e)]$ for every edge $e \in E(G)$ according to the adjacency determined by $G$, and choosing the smooth local model~\eqref{eq:smoothlocalmodel} at each vertex of valence $3$ or higher (cf.~Proposition~3.6 in~\cite{MikhalkinZharkov}). We henceforth assume that all tropical curves are smooth, in other words we consider metric graphs $(G,\ell)$ having finitely many univalent points at infinity. In particular, trees are assumed to have all of their leaf vertices at infinity.  

Let $\Gamma$ be a (smooth) tropical curve with a chosen orientation. A \emph{harmonic 1-form} is a global section of $\Omega^1_\Gamma$. More explicitly, a harmonic 1-form $\omega = \sum_{e\in E(\Gamma)} a_e d e$ is given by the choice of a coefficient $a_e \in \ZZ$ subject to the condition 
\[ \sum_{e \text{ entering } v} a_e - \sum_{e \text{ leaving } v} a_e = 0 \]
at every vertex $v$ of $\Ga$. 

A \emph{harmonic function} on $\Gamma$ is a section of $\Aff_\Gamma$, i.e. a continuous function $f : U \to \RR$ on an open subset $U\subset \Ga$ that is linear with integer slope on every edge, and such that the sum of outgoing slopes of $f$ at every vertex is zero. Recording the slopes of $f$, we obtain a harmonic $1$-form on $U$. We note that harmonic functions are constant near infinite extremal vertices.

Let $\Ga_1=(G_1,\ell_1)$ and $\Ga_2=(G_2,\ell_2)$ be tropical curves, and let $f:G_1\to G_2$ be a harmonic morphism of graphs with degree function $d_f$, where we recall that we do not allow graph morphisms to contract edges. If $d_f(e)=\frac{\ell_2(f(e))}{\ell_1(e)}$ for each edge $e\in E(G_1)$, then we define the associated \emph{harmonic morphism of tropical curves} $f:\Ga_1\to \Ga_2$, which is an affine linear map on each edge $e$ of $\Ga_1$ with slope or \emph{dilation factor} equal to $d_f(e)$. It is elementary to verify that $f$ is a morphism of rational polyhedral spaces. Conversely, a surjective morphism of rational polyhedral spaces $f:\Ga_1\to \Ga_2$ induces a harmonic morphism of graphs (with respect to appropriately chosen models) if it does not contract any edges, and such a morphism has a well-defined global degree $\deg f$.
Given a harmonic morphism $f:\Ga_1\to \Ga_2$ and a point $x\in \Ga_1$, we denote $d_f(x)=d_f(e)$ if $x$ lies in the interior of the edge $e$, and $d_f(x)=d_f(v)$ if $x$ corresponds to a vertex $v$ (with respect to an appropriate model). We observe that, given a harmonic morphism of graphs $f:G_1\to G_2$, a choice of edge lengths on $G_2$ uniquely determines edge lengths on $G_1$ in such a way that $f$ induces a harmonic morphism of tropical curves.

In parallel to the graph case, we say that a harmonic morphism $f:\Ga_1\to \Ga_2$ of tropical curves is \emph{free} if $d_f(x)=1$ for all $x\in \Ga_1$ (equivalently, if it is a covering isometry), \emph{dilated} if it is not free, and a \emph{double cover} if it has global degree 2.

We note that an arbitrary harmonic morphism $f:(G_1,\ell_1)\to (G_2,\ell_2)$ of metric graphs can be augmented to a harmonic morphism of smooth tropical curves in the following way. For each univalent vertex $v_2\in V(G_2)$, attach a compact infinite ray $l$ to $v_2$, then for each $v_1\in f^{-1}(v_2)$ attach $d_f(v_1)$ compact infinite rays to $v_1$ and map them with degree $1$ to $l$. From now on we will always assume that harmonic morphisms are in fact harmonic morphisms between smooth tropical curves.

\begin{remark}
    \label{rem:noweights}
    A tropical curve $\Ga$ arising as the tropicalization of an algebraic curve naturally comes with a \emph{vertex weight function}, recording the genera of the irreducible components of the special fiber. These vertex weights appear to play no role in the tropical $n$-gonal construction, and we do not consider them.
\end{remark}

\subsection{Divisors on tropical curves}
 
Let $\Gamma$ be a tropical curve. A \emph{divisor} $D$ on $\Gamma$ is a finite formal $\ZZ$-linear combination of points on $\Gamma$, i.e.~$D = \sum a_x \cdot x$ with $a_x = 0$ for almost all $x \in \Gamma$. Denote the group of divisors on $\Gamma$ by $\Div(\Gamma)$. The \emph{degree} of a divisor $D$ is $\deg D = \sum a_x$. A divisor is called \emph{effective} if $a_x \geq 0$ for all $x \in \Gamma$, in which case we write $D \geq 0$. Denote the set of effective divisors of degree $n$ by $\Div^+_n(\Gamma)$. 

A \emph{rational function} $f : \Gamma \to \RR$ is a piecewise linear function with integer slopes. A rational function induces a divisor as follows
\[ \div f = \sum_{x \in \Gamma} \big(\text{sum of outgoing slopes at $x$} \big) \cdot x. \] 
Any divisor of the form $\div f$ is called a \emph{principal divisor}, and the subgroup of principal divisors is denoted by $\Prin(\Gamma) \subseteq \Div(\Gamma)$. Two divisors $D_1$ and $D_2$ are \emph{linearly equivalent} if $D_1 - D_2$ is a principal divisor, in which case we write $D_1 \sim D_2$. The \emph{Picard group} of $\Gamma$ is defined as 
\[ \Pic(\Gamma) = \Div(\Gamma) / \Prin(\Gamma) \qquad \text{and} \qquad \Pic_k(\Gamma) = \big\{ [D] \in \Pic(\Gamma) \mid \deg D = k \big\} . \]
Every $\Pic_k(\Gamma)$ is a torsor over the \emph{Picard variety} $\Pic_0(\Gamma)$.

\subsection{The $n$-gonal and Recillas construction for tropical curves} In this section, we extend the $n$-gonal and Recillas construction to tropical curves. We recall that our definition of gonality involves maps to trees instead of Baker--Norine rank.

\begin{definition}
    An \emph{$n$-gonal} tropical curve is a tropical curve $\Gamma$ together with a harmonic map $f : \Gamma \to K$ of degree $n$ to a metric tree $K$. For $n = 2, 3, 4$ we will also use the terms \emph{hyperelliptic}, \emph{trigonal}, and \emph{tetragonal}, respectively. A tetragonal curve is called \emph{generic} if for all $x\in K$ the fiber $f^{-1}(x)$ has dilation profile $(3,1)$, $(2,1,1)$, or $(1,1,1,1)$. A double cover $\tGa\to \Ga\to K$ of a hyperellipic curve is called \emph{generic} if $K$ has no point which has two preimages in $\Ga$, each of which has a unique preimage in $\tGa$. 

\end{definition}

A harmonic morphism of tropical curves is uniquely specified by giving a harmonic morphism of graphs together with an edge length function on the target. This observation allows us to directly lift the tropical $n$-gonal construction and the tropical Recillas construction from graphs to tropical curves:

\begin{definition}
    Let $\tilde \Gamma \to \Gamma \to K$ be a double cover of an $n$-gonal tropical curve. The \emph{$n$-gonal construction} is the degree $2^n$ harmonic morphism of tropical curves $\tilde \Pi \to K$ that arises by running the construction of Section~\ref{sec:n_gonal_harmonic} on the underlying tower of graphs, and then endowing $\tPi$ with the edge-length function induced by $K$. Similarly, the tropical Recillas construction associates to a generic tetragonal curve $\Pi\to K$ a tower $\tGa\to\Ga\to K$ consisting of a free double cover of a trigonal curve. 

\end{definition}

All results from Section \ref{sec:graphs} carry over to the setting of tropical curves. In particular, the fiberwise description over a point $x \in K$ still holds, where $x$ may now correspond to a vertex or an interior edge point. We may view an edge $e$ of a graph model of $K$ as a generic point for the points of that edge in $K$. In this sense our definition of the construction for generic fibers ensures that the fibers of the $n$-gonal construction depend continuously on the metric realization. With this point of view, the root map at the level of graphs can now be understood as the continuous limit for $x \in K$ approaching a vertex.

\begin{remark} We can restate the $n$-gonal and Recillas constructions directly for tropical curves in the language of divisors, in a way that is consistent with the graph-theoretic definitions given in Sections~\ref{sec:n_gonal_harmonic} and~\ref{sec:Recillas}. Let $\pi:\tilde \Gamma \to \Gamma$ be a free double cover of an $n$-gonal curve $f:\Gamma\to K$, and let $\tilde \Pi \to K$ be the output of the $n$-gonal construction. As a set $\tPi$ is given by the construction of Section~\ref{sec:n_gonal_harmonic} applied pointwise in the metric graphs: 
\begin{equation}
    \label{eq:tildePi_in_Symn}
    \tilde \Pi = \Big\{ x_1 + \cdots + x_n \in \Div_n^+(\tilde \Gamma) \mathrel{\Big|} \exists x \in K : \pi(x_1) + \cdots + \pi(x_n) = \sum_{y \in f^{-1}(x)} d_f(y)\cdot y \Big\}.
\end{equation}
Even more is true: since $K$ is a tree, the natural graph model for $\tilde \Gamma$ is loop-free and hence \cite{BrandtUlirsch} gives a polyhedral structure on $\Div_n^+(\tilde\Gamma)$. The graph structure of $\tilde\Pi$ is precisely the restriction of this polyhedral structure. 
A similar description can be given for the tropical Recillas construction. Let $k : \Pi \to K$ be a generic tetragonal tropical curve, and let $\tilde \Gamma \to \Gamma \to K$ be the output of the tropical Recillas construction. Then as a set
\begin{equation} \label{eq:tildeGamma_in_Sym2}
    \tilde \Gamma = \Big\{ x_1 + x_2 \in \Div^+_2(\Pi) \mathrel{\Big|} \exists x_3 + x_4 \in \Div^+_2(\Pi) \text{ and } x \in K \text{ such that } x_1 + x_2 + x_3 + x_4 = \sum_{y \in k^{-1}(x)} d_k(y) \cdot y \Big\} ,
\end{equation}
and the involution whose quotient is the double cover $\tGa\to \Ga$ is the one sending $x_1+x_2$ to $x_3+x_4$. However, it turns out that defining $\tPi$ and $\tGa$ in this way does not naturally induce the correct edge lengths. 
For example, in the context of the trigonal construction, consider a $(2,3)$-tower $\tilde \Gamma \overset{2}{\rightarrow} \Gamma \overset{3}{\rightarrow} K$ with an edge $e\in E(K)$ of Type I (in the language of Definition~\ref{def:trigonal_types}). Then both edges $\tilde e^+$ and $\tilde e^-$ in $\tilde\Gamma$ above $e$ are of the same length, say $a$ (so that the length of $e$ is $3a$). For the construction of $\Pi$ there are two relevant cells in $\Div_3^+(\tilde\Gamma)$, namely 
\[ (\tilde e^+ \times \tilde e^+)/S_2 \times \tilde e^- \qquad \text{and} \qquad (\tilde e^- \times \tilde e^- \times \tilde e^-)/ S_3. \]
The two edges of $\Pi$ arising from $\tilde e^+$ and $\tilde e^-$ are given as the diagonals of these cells, both of which have lattice length $a$, in contrast to the correct edge lengths $a$ and $3a$. This is one of the main reasons that we define the $n$-gonal and Recillas constructions for combinatorial graphs first, and then import the construction into the setting of tropical curves.

\end{remark}

\subsection{Tropical cycles}

We will now recall the definition of the tropical cycle class groups $Z_k(X)$ associated to a rational polyhedral space $X$. We first recall the definition in affine space, following~\cite{GrossShokrieh_homology}.

\begin{definition}
    Let $\RR^n$ be endowed with the integral structure $\ZZ^n \subseteq \RR^n$. A \emph{tropical fan $k$-cycle} on $\RR^n$ is a function $A : \RR^n \to \ZZ$ satisfying the following properties:
    \begin{enumerate}
        \item For all $x \in \RR^n$ and $\lambda \in \RR_{>0}$ we have $A(\lambda x ) = A(x)$.
        \item The \emph{support} $|A| = \overline{\{ x \in \RR^n \mid A(x) \neq 0 \}}$ is a purely $k$-dimensional rational polyhedral set.
        \item $A$ is locally constant on $|A|^\reg$ and is equal to 0 on $|A| - |A|^\reg$.
        \item $A$ satisfies the so-called \emph{balancing condition}. This is a condition at every codimension 1 cell $\tau$ of the rational polyhedral structure of $|A|$ and requires the sum of the outwards facing lattice normal vectors of incident maximal cells, weighted by the values of $A$, to be contained in the tangent space of $\tau$. The condition does not depend on the chosen fan structure of $|A|$ and we will only need to check balancing when $k = 1$, in which case the only codimension one cell is the origin $\tau=0$. Let $e$ range over the 1-dimensional cones in $|A|$ and for each $e$ let $\eta_e \in e \cap \ZZ^n$ be an outwards facing primitive tangent vector of $e$. The balancing condition at 0 is 
        \begin{equation}
            \label{eq:balancing}
            \sum_e A(\eta_e) \eta_e = 0 .
        \end{equation}
        We note that for a $1$-cycle $A$ (whose support $|A|$ is a graph), verifying balancing does not involve the values at the vertices (points of $|A|-|A|^\reg$), and we will often ignore condition (3) and allow $A$ to have arbitrary vertex values.
    \end{enumerate}
\end{definition}

The idea is now to define a tropical $k$-cycle on an arbitrary rational polyhedral space by requiring it to look locally like a tropical fan $k$-cycle. 

\begin{definition}
    Let $X$ be a rational polyhedral space. A \emph{local face structure} at a point $x \in X$ is a finite polyhedral complex $\Sigma$ such that 
    \begin{enumerate}
        \item $x$ is contained in the topological interior of $|\Sigma|$,
        \item there exists a chart $U \to V \subseteq \overline{\RR}^n$ of $X$ such that $|\Sigma| \subseteq U$,
        \item $x$ is contained in every inclusion maximal cell of $\Sigma$.
    \end{enumerate}
\end{definition}

Face structures are higher-dimensional analogues of graph models for tropical curves.

\begin{definition}
    Let $X$ be a rational polyhedral space. A \emph{tropical $k$-cycle} is a function $A : X \to \ZZ$ such that the following properties hold.
    \begin{enumerate}
        \item $A$ is \emph{locally constructible}, i.e. for every $x \in X$ there is a local face structure $\Sigma$ at $x$ such that the restriction to the relative interior $A|_{\operatorname{relint}(\sigma)}$ is locally constant for every $\sigma \in \Sigma$.
        \item For every $x \in X$ the germ of $A$ at $x$ (extended to be constant along all lines through the origin) defines a tropical fan $k$-cycle on the real tangent space $T^\ZZ_x (X) \otimes_\ZZ \RR$ at $x$.
    \end{enumerate}
    The set of tropical $k$-cycles on $X$ is denoted $Z_k(X)$.

\end{definition}

The following two examples show that smooth tropical curves behave well from the viewpoint of intersection theory. Ultimately this is the reason for our standing smoothness assumption.

\begin{example}
    \label{ex:fundamental_cycle}
    Let $\Gamma$ be a tropical curve (which we assume be to smooth as always). A function $A : \Gamma \to \ZZ$ with value 0 on every vertex of valence $>2$ is balanced if and only if its value on every edge is the same. 
    If this value is 1 everywhere, then this defines the \emph{fundamental cycle} $[\Gamma]\in Z_1(\Ga)$ of $\Gamma$. 
    
    We emphasize that if $\Gamma$ is not smooth then it does not admit a fundamental cycle. For example, the balancing condition \eqref{eq:balancing} cannot be satisfied at a finite 1-valent vertex. At a 1-valent vertex at infinity, however, the balancing condition is trivially satisfied, since all affine functions are locally constant and hence the tangent space is 0. In particular, the primitive tangent vector of an infinite edge at the vertex at infinity is 0 and Equation \eqref{eq:balancing} is satisfied. 
\end{example}

\begin{example} \label{ex:diagonal_cycle}
    Let $\Gamma$ be a (smooth) tropical curve and define $[\Delta_\Gamma] : \Gamma^2 \to \ZZ$ by $[\Delta_\Gamma](x,y) = 1$ if $x = y$ and $x$ is not a vertex of valence $> 2$ and $[\Delta_\Gamma](x,y) = 0$ otherwise. We claim that this is a tropical 1-cycle on $\Gamma^2$, the \emph{diagonal cycle}.
    
    A choice of graph model on $\Gamma$ provides a polyhedral complex structure on $\Gamma^2$. Subdividing cells of the form $e \times e$ for any edge $e$ of $\Gamma$ provides a face structure that shows that $[\Delta_\Gamma]$ is (locally) constructible. We check balancing.
    Locally at a point $(x, x)$, where $x$ is not a vertex of $\Gamma$, the support of $[\Delta_\Gamma]$ looks like the diagonal in $\RR^2$. This is clearly balanced because the sum of outwards facing primitive tangent vectors at the origin is 0.
    Now assume $x$ is a vertex of $\Gamma$. Denote the primitive tangent vectors of the edges of $\Gamma$ incident to $x$ by $\eta_1, \ldots, \eta_n$. Then the primitive tangent vectors of $[\Delta_\Gamma]$ at $(x,x)$ are 
    \[\begin{pmatrix}
        \eta_i \\ \eta_i
    \end{pmatrix} \in T^\ZZ_{(x,x)}(\Gamma^2) \cong T^\ZZ_x(\Gamma) \times T^\ZZ_x(\Gamma) \qquad \qquad \text{for } i = 1, \ldots, n,\]
    and again we see that the sum is 0 because $\Gamma$ was assumed smooth, i.e. $\sum_i \eta_i = 0$.
\end{example}

Let $A, B : X \to \ZZ$ be tropical $k$-cycles on $X$. The sum function $A+B:X\to \ZZ$ is not, in general, a tropical $k$-cycle. However, there exists a tropical $k$-cycle agreeing with the algebraic sum $A+B$ away from the non-regular locus $|A|\setminus |A|^\reg \cup |B|\setminus |B|^\reg$. Denoting this cycle $A+B$ by abuse of notation, we obtain a group structure on $Z_k(X)$.

Now let $f : X \to Y$ be a proper and surjective morphism of $k$-dimensional rational polyhedral spaces. There is a \emph{pushforward} $f_\ast : Z_k(X) \to Z_k(Y)$ defined as follows (see Definition~3.6 of \cite{GrossShokrieh_homology}). Let $A \in Z_k(X)$. At $y \in Y^\reg \setminus f(X \setminus X^\reg)$ define
\begin{equation}
    \label{eq:definition_push_forward}
    f_\ast A(y) = \sum_{x \in f^{-1}(y)} \big[ T_y^\ZZ Y : d_x f(T_x^\ZZ X) \big] A(x),
\end{equation}
where the lattice index is taken to be 0 if it is not finite and we set $f_\ast A(y) = 0$ for all other points of $Y$. Similar to the case of addition of tropical cycles, this $f_\ast A$ is in general not a tropical cycle, but there is a tropical pushforward cycle (also denoted $f_*A$) that coincides with $f_*A$ away from a locus of dimension at most $k-1$.

\begin{example}
    Let $\Gamma$ and $\Pi$ be smooth tropical curves and let $\pi : \Gamma \to \Pi$ be a harmonic map of degree $d$. Then $\pi_*[\Gamma] = d[\Pi]$. To see this, it suffices to check that $\pi_*[\Gamma](y) = d$ for every $y$ in the interior of an edge of $\Pi$. By definition of pushforward and Equation~\eqref{eq:globaldegree} we have
    \[ \pi_*[\Gamma](y) = \sum_{x \in \pi^{-1}(y)}  \big[ T_y^\ZZ \Pi : d_x\pi (T_x^\ZZ \Gamma) \big] = \sum_{x \in \pi^{-1}(y)} d_\pi(x) = n, \]
    where $d_\pi(x)$ is the dilation factor on the edge to which $x$ belongs and $d_x\pi$ is the differential of $\pi$ at $x$.
\end{example}

\subsection{Tropical homology} Finally, we give a brief introduction to tropical homology and cohomology groups (which were introduced in~\cite{IKMZ}), following the paper~\cite{GrossShokrieh_homology}.
Let $X$ be a rational polyhedral space, possibly with boundary.

\begin{definition}
	Let $p \geq 1$. Define the sheaf $\Omega^p_X$ of \emph{tropical $p$-forms} to be the image of 
	\[ \bigwedge^p \Omega^1_X \longrightarrow i_\ast \big( \bigwedge^p \Omega_X^1|_{X^\mathrm{reg}} \big), \]
	where $i : X^\mathrm{reg} \to X$ is the inclusion. 
\end{definition}

There is a maximal stratification of $X$ such that $\Omega^1_X$ is locally constant on every stratum. A singular $q$-simplex, i.e. a continuous map $\sigma : \Delta^q \to X$, is \emph{allowable} if every open face of $\Delta^q$ maps into a single stratum of $X$. 
Denote by $\ZZ_{\sigma(\Delta^q)}$ the constant sheaf on $\sigma(\Delta^q)$ with values in $\ZZ$. Then the $(p,q)$-th chain group is defined as
\[ C_{p,q}(X) = \bigoplus_{\sigma : \Delta^q \to X \text{ allowable}} \Hom\big(\Omega^p_X, \ZZ_{\sigma(\Delta^q)} \big). \]
Elements in $C_{p,q}(X)$ are denoted as $\sum_\sigma \sigma \otimes \eta_\sigma$, where $\eta_\sigma \in \Hom\big(\Omega^p_X, \ZZ_{\sigma(\Delta^q)} \big)$. The boundary map of $C_{p, \bullet}(X)$ is given as
\begin{equation*}
    C_{p, q+1}(X) \longrightarrow C_{p,q}(X), \qquad
    \sigma \otimes \eta \longmapsto \sum_{\tau \in \partial \sigma} \tau \otimes (r_\tau \circ \eta),
\end{equation*}
where $\partial$ is the boundary map from singular homology and $r_\tau : \ZZ_{\sigma(\Delta^{q+1})} \to \ZZ_{\tau(\Delta^q)}$ is the restriction map. By abuse of notation we denote the boundary maps of $C_{p, \bullet}(X)$ again by $\partial$ and we define the $(p,q)$-th \emph{tropical homology group} $H_{p,q}(X) = H_q(C_{p, \bullet}(X))$. The cochain complexes are $C^{p, \bullet} = \Hom(C_{p, \bullet}, \ZZ)$ and the \emph{tropical cohomology groups} are $H^{p,q}(X) = H^q(C^{p, \bullet}(X))$.

Recall that there is a \emph{tropical cycle class map} 
$$ \cyc : Z_k(X) \longrightarrow H_{k,k}(X) $$ 
which assigns to any tropical $k$-cycle a class in tropical homology. We often suppress $\cyc$ from the notation and identify a tropical cycle in $Z_k(X)$ with its image in $H_{k,k}(X)$. This map is defined for rational polyhedral spaces with boundary in~\cite[Section 5]{GrossShokrieh_homology}, and is given a convenient description in the $k=1$ case for boundaryless polyhedral spaces in~\cite{GrossShokrieh_Poincare}. We recall the latter formula for $A\in Z_1(X)$, generalized to the case when the support of $|A|$, which is a graph, is allowed to have boundary vertices, but no boundary edges. 

For every edge $e$ in $|A|$ choose a generator for $T^{\ZZ}_x|A|$ for some $x \in e$. By parallel transport this gives rise to a generator for any $T^{\ZZ}_y|A|$ with $y \in e$ and hence a morphism $\Omega^1_{|A|} \to \ZZ_e$. Precomposing with $\Omega^1_X \to \Omega^1_{|A|}$ induced by the inclusion $|A| \hookrightarrow X$ we obtain $\eta_e \in \Hom(\Omega^1_X, \ZZ_e)$. Let $\gamma_e : \Delta^1 \to X$ be a parametrization of $e$ in the direction given by $\eta_e$. Then 
\[ \cyc(A) = \sum_e A(e) \gamma_e \otimes \eta_e \in C_{1,1}(X). \]
Let us check that $\partial \cyc (A) = 0$. To do so, let $v$ be a finite vertex of $|A|$ and assume that all edges $e$ incident to $v$ have been oriented away from $v$. Then each $e$ contributes $-A(e)\eta_e|_{v}$ to $\partial \cyc(A)$. But now the sum over these contributions is 0 because $A$ was assumed to be a tropical cycle and in particular balanced. On the other hand, if $v$ is an infinite vertex, then the stalk of $\Omega^1_{|A|}$ at $v$ is trivial, hence the element $\eta_e$ (corresponding to any incident edge $e$) is equal to 0.

There are natural pushforward maps of tropical homology classes and pullback maps of tropical cohomology classes along morphisms of rational polyhedral spaces. These maps are compatible with the cycle class map in the sense that proper pushforward of tropical cycles and pushforward of tropical homology classes commute with the cycle class map \cite[Proposition~5.6]{GrossShokrieh_homology}. Finally, there are cup and cap products \cite[Section~4.6]{GrossShokrieh_homology}, 
\[
\smile:H^{p,q}(X)\times H^{p',q'}(X)\longrightarrow H^{p+p',q+q'}(X),\qquad 
\frown:H_{p,q}(X)\times H^{i,j}(X)\longrightarrow H_{p-i,q-j}(X),
\]
and the latter gives rise, if $X$ is \emph{smooth}, to Poincar\'e duality \cite[Corollary~6.8]{GrossShokrieh_homology}
\[
H^{p,q}(X)\cong H_{n-p,n-q}(X),\qquad c \longmapsto \cyc [X]\frown c.
\]
Smoothness is defined in Definition 6.1 in \cite{GrossShokrieh_homology} and we do note repeat the general definition, but simply note that we only use this isomorphism when $X$ is a power of a smooth tropical curve, which is smooth or when $X$ is a tropical abelian variety, which is also smooth. Moreover, note that if $X$ is smooth, it has a fundamental class.

\subsection{Tropical Cartier divisors} \label{sec:Cartier}

Let $X$ be a rational polyhedral space. A \emph{rational function} on $X$ is a continuous function $f : X \to \RR$ that is piecewise affine on every chart. Denote the sheaf of rational functions by $\calM_X$. Clearly, every affine function is rational, so there is an inclusion $\Aff_X \to \calM_X$. Denote the quotient $\calM_X / \Aff_X$ by $\calDiv(X)$, so that there is a short exact sequence of sheaves
\begin{equation}
    \label{eq:ses_Cartier_divisors}
    0 \longrightarrow \Aff_X \longrightarrow \calM_X \longrightarrow \calDiv(X) \longrightarrow 0.
\end{equation}
The group of global sections $\Div(X)= \Gamma \big(X, \calDiv(X) \big)$ is the group of \emph{Cartier divisors} on $X$. If $f : X \to Y$ is a morphism of rational polyhedral spaces, then there is an induced pullback map on Cartier divisors $f^* : \Div(Y) \to \Div(X)$.

The group $H^1(X, \Aff_X)$ classifies \emph{tropical line bundles} on $X$, and the short exact sequence \eqref{eq:ses_Cartier_divisors} gives rise to a boundary homomorphism $\Div(X) = H^0\big(X, \calDiv(X)\big) \to H^1(X, \Aff_X)$ that associates to a Cartier divisor $D$ a tropical line bundle $\calL(D)$. 
Pullback of Cartier divisors commutes with this association \cite[Proposition 3.15]{GrossShokrieh_homology}. Furthermore, the short exact sequence defining $\Omega^1_X$
\[ 0 \longrightarrow \RR_X \longrightarrow \Aff_X \longrightarrow \Omega^1_X \longrightarrow 0 \]
gives rise to the \emph{first Chern class map}
$c_1 : H^1(X, \Aff_X) \to H^1(X, \Omega^1_X) = H^{1,1}(X)$.

Finally, there is a natural intersection pairing
\[ \Div(X) \times Z_k(X) \longmapsto Z_{k-1}(X), \qquad (D,A)\longmapsto D\cdot A. \]
If $X$ is smooth, then in particular $X$ admits a fundamental cycle and for $k = \dim X$ this gives an isomorphism $\Div(X) \cong Z_{\dim X -1}(X)$. Again, we note that we only use this isomorphism when $X$ is a power of a tropical curve, which is smooth. We recall from \cite[Proposition 5.12]{GrossShokrieh_homology} that $\cyc(D \cdot A) = \cyc(A) \frown c_1(\calL(D))$ for every Cartier divisor $D$ and tropical cycle $A$.

%%%%%%%%%%%%%%%%%%%%%%%%%%%%%%%%%%%%%%%%%%%%%%%%%%%%%%
\section{Tropical abelian varieties}
%and Welters' criterion
\label{sec:tavs}
%%%%%%%%%%%%%%%%%%%%%%%%%%%%%%%%%%%%%%%%%%%%%%%%%%%%%%

In this section, we recall the theory of tropical abelian varieties.  The definitions that we use were introduced in~\cite{LenZakharov} and differ slighly from the standard definitions (see e.g.~\cite{foster2018non}), but are equivalent to them. We prove a number of elementary results concerning morphisms of tropical abelian varieties.
We then recall the Jacobian variety of a tropical curve (already introduced in~\cite{MikhalkinZharkov}) and introduce the continuous Prym variety of a double cover (modifying the original construction from~\cite{JensenLen}), and show that they satisfy natural universal properties. Finally, we prove Theorem~\ref{thm:tropical_Welters_criterion}, which is a tropical version of the homological formula for the pushforward of the fundamental class under the Abel--Prym map (which is classically a part of Welters's criterion characterizing Prym varieties).

\subsection{Integral tori} Let $\Lambda$ and $\Lambda'$ be finitely generated free abelian groups of the same rank and let $[\cdot, \cdot] : \Lambda \times \Lambda' \to \RR$ be a non-degenerate pairing. The triple $(\Lambda, \Lambda', |\cdot, \cdot])$ defines a \emph{real torus with integral structure}  $\Sigma=\Hom(\Lambda, \RR) / \Lambda'$, or simply an \emph{integral torus}, where the inclusion $\Lambda' \subseteq \Hom(\Lambda, \RR)$ is given by $\lambda' \mapsto [\cdot,\lambda']$.
The \emph{dual torus} $\Sigma^\vee = \Hom(\Lambda', \RR) / \Lambda$ is defined by the transposed triple $(\Lambda', \Lambda, [\cdot, \cdot]^t)$. 
The \emph{dimension} of an integral torus is $\dim_{\RR}\Si=\rk \La=\rk \La'$. We note that integral tori admit a group structure, which is the descent of the group structure of the universal cover $\Hom (\Lambda, \RR)$. 
From now on we abuse notation and refer to the triples as integral tori as well. 

\begin{remark} \label{rem:affine_fcn_on_tori}
    Integral tori are rational polyhedral spaces as follows. Identifying $\Hom(\Lambda, \ZZ)$ with $\ZZ^g$ endows the universal cover $\Hom(\Lambda, \RR)$ with the structure of a rational polyhedral space. An integral affine linear function on $\Hom(\Lambda, \RR)$ is the sum of a linear function taking integer values on the lattice $\Hom(\Lambda, \ZZ)$ and a constant real shift. In other words, affine functions are precisely elements of $\Hom \big(\Hom(\Lambda, \ZZ), \ZZ \big) \oplus \RR \cong \Lambda \oplus \RR$.  
    The torus inherits the rational polyhedral structure from the universal covering via the quotient map.     
    Note that $\Omega^1_\Sigma (\Sigma) = \Lambda$ and $H_1(\Sigma, \ZZ) = \Lambda'$ for any integral torus $\Sigma = (\Lambda, \Lambda', [\cdot,\cdot])$. 
    Moreover, it is easy to see e.g. via \cite[Lemma~6.3]{GrossShokrieh_homology} that integral tori are smooth rational polyhedral spaces.
\end{remark}

We first classify morphisms of integral tori (as rational polyhedral spaces). 
Recall that a holomorphic map of complex tori factors as a group homomorphism followed by a translation. An identical classification holds for integral tori. We first define the two types of maps.

\begin{definition} Let $\Si$ be an integral torus. For every $y \in \Sigma$ the \emph{translation} $t_y : \Sigma \to \Sigma$ is given by $t_y(x)=x + y$.
\end{definition}

It is clear that translations are morphisms of rational polyhedral spaces, inducing identity maps on $\Omega^1_{\Sigma}(\Sigma)$ and $H^{p,q}(\Sigma)$. 

\begin{definition}
    A \emph{homomorphism of integral tori} $f = (f^\#, f_\#) : (\Lambda_1, \Lambda_1', [\cdot, \cdot]_1) \to (\Lambda_2, \Lambda_2', [\cdot, \cdot]_2)$ consists of a pair of homomorphisms $f^\# : \Lambda_2 \to \Lambda_1$ and $f_\# : \Lambda_1' \to \Lambda_2'$ satisfying the relation 
    \begin{equation}
        \label{eq:homomorphism_condition}
        \big[f^\#(\lambda_2),\lambda_1'\big]_1 = \big[ \lambda_2,f_\#(\lambda_1') \big]_2
    \end{equation}
    for all $\la'_1\in \La'_1$ and $\la_2\in \La_2$.    
    The maps $f^\#$ and $f_\#$ necessarily have the same rank, denoted $\rk f$.
    The \emph{dual homomorphism} of $f$ is given by the transposed pair $f^\vee = (f_\#, f^\#) : \Sigma_2^\vee \to \Sigma_1^\vee$. Given another homomorphism of integral tori $g = (g^\#, g_\#) : (\Lambda_2, \Lambda_2', [\cdot, \cdot]_2) \to (\Lambda_3, \Lambda_3', [\cdot, \cdot]_3)$, the \emph{composition} with $f$ is given by $g \circ f = (f^\# \circ g^\#, g_\# \circ f_\#)$.
\end{definition}

Note that for a homomorphism $f = (f^\#, f_\#)$, the $\Hom$-dual $\Hom(\La_1,\RR)\to \Hom(\La_2,\RR)$ of $f^\#$ restricts to $f_\#$ on $\La'_1$ and hence descends to a group homomorphism on the underlying tori
\begin{equation} \label{eq:morphism_on_tori}
    \Si_1=\Hom(\Lambda_1, \RR) / \Lambda_1' \longrightarrow \Si_2=\Hom(\Lambda_2, \RR) / \Lambda_2',
\end{equation}
which is a map of rational polyhedral spaces.
By abuse of notation we denote the map in Equation~\eqref{eq:morphism_on_tori} again by $f$ and from now on we will conflate the representation of a homomorphism as a pair $(f^\#, f_\#)$ and its underlying description as a (point-wise) map of rational polyhedral spaces. This is compatible with composition of homomorphisms.

\begin{lemma} \label{lem:homomorphism_translation}
    Let $\Sigma_i = (\Lambda_i, \Lambda_i', [\cdot, \cdot]_i)$ for $i = 1,2$ be integral tori and let $f : \Sigma_1 \to \Sigma_2$ be a map of rational polyhedral spaces. Then $f$ factors uniquely as a homomorphism $g = (g^\#, g_\#) : \Sigma_1 \to \Sigma_2$ followed by a translation $t : \Sigma_2 \to \Sigma_2$.
\end{lemma}

\begin{proof}
    Define $g = t_{-f(0)} \circ f$, then clearly $g$ is a map of rational polyhedral spaces with $g(0) = 0$.
    In particular, $g$ pulls back affine linear functions defined in a neighborhood of $0 \in \Sigma_2$ to affine linear functions on a neighborhood of $0 \in \Sigma_1$. Since $g(0) = 0$, the pullback of a linear function is in fact linear. 
    But for any integral torus $\Sigma = (\Lambda, \Lambda', [\cdot, \cdot])$, linear functions in a neighborhood of $0$ are simply given by elements of $\Hom \big(\Hom(\Lambda, \ZZ), \ZZ \big) \cong \Lambda$, see Remark~\ref{rem:affine_fcn_on_tori}. Hence pullback defines a group homomorphism $g^\# : \Lambda_2 \to \Lambda_1$ whose $\Hom$-dual induces the map $g$ on the tori, which is therefore a homomorphism of integral tori. \qedhere
\end{proof}

\begin{lemma} \label{lem:swap_homomorphism_translation}
    Let $\Sigma_i = (\Lambda_i, \Lambda_i', [\cdot, \cdot]_i)$ for $i = 1,2$ be integral tori and let $f = (f^\#, f_\#)$ be a homomorphism $\Sigma_1 \to \Sigma_2$. Then for any $y \in \Sigma_1$ the diagram
    \begin{equation*}
        \begin{tikzcd}
            \Sigma_1 \arrow[d, "t_y"] \arrow[r, "f"] & \Sigma_2 \arrow[d, "t_{f(y)}"] \\
            \Sigma_1 \arrow[r, "f"] & \Sigma_2
        \end{tikzcd}
    \end{equation*}
    commutes.
\end{lemma}

\begin{proof}
    This is clear because $f$ is a homomorphism with respect to the group structures of the integral tori $\Sigma_1$ and $\Sigma_2$, in other words $f(x + y) = f(x) + f(y)$. \qedhere
\end{proof}

We now define the kernel, cokernel, and image of a homomorphism of integral tori. For an abelian group $A$, we denote by $A^{\torf}=A/A^{\mathrm{tor}}$ the quotient by its torsion subgroup. For a pair of lattices $\Lambda\subset \Lambda'$, we define the \emph{saturation} of $\Lambda$ in $\Lambda'$ as 
\[
\Lambda^{\sat}=\La'\cap (\La\otimes \QQ)=\big\{\lambda'\in \Lambda':n\lambda'\in \Lambda\mbox{ for some }n\in \ZZ \big\},
\]
and note that $(\La'/\La)^{\torf}\cong \La'/\La^{\sat}$.
Let $\Sigma_i = (\Lambda_i, \Lambda_i', [\cdot, \cdot]_i)$ for $i = 1,2$ be integral tori and let $f = (f^\#, f_\#)$ be a homomorphism $\Sigma_1 \to \Sigma_2$. We consider the following integral tori:
\begin{equation}
(\Ker f)_0  = \Big( (\Coker f^\#)^\mathrm{tf}, \quad \Ker f_\#, \quad [\cdot, \cdot]_K \Big),\qquad     \Coker f = \Big( \Ker f^\#, \quad (\Coker f_\#)^\torf, \quad [ \cdot, \cdot]_C \Big),
\label{eq:kercoker}
\end{equation}
\begin{equation}
\Im f=\Big( \La_2/\Ker f^\#, \quad (\Im f_\#)^{\mathrm{sat}}, \quad [\cdot, \cdot]_I \Big),
\label{eq:im}
\end{equation}
where the pairings $[\cdot,\cdot]_K$, $[\cdot, \cdot]_C$, and $[\cdot, \cdot]_I$ are induced by the pairings $[\cdot, \cdot]_1$, $[\cdot, \cdot]_2$, and $[\cdot, \cdot]_2$, respectively. The natural maps on the lattices induce the following sequence of homomorphisms of integral tori:
\begin{equation}
(\Ker f)_0\overset{i}{\longrightarrow}\Sigma_1\overset{p}{\longrightarrow}\Im f\overset{j}{\longrightarrow}\Sigma_2\overset{q}{\longrightarrow}\Coker f.
\label{eq:kerimcoker}
\end{equation}

We now compare these definitions to their group-theoretic counterparts. We start by showing that $\Im f$ is in fact the group-theoretic image of $f$.

\begin{lemma} \label{lem:image_of_tori}
    Let $f:\Si_1\to \Si_2$ be a homomorphism of integral tori. The homomorphism $j:\Im f\to \Sigma_2$ identifies $\Im f$ with the group-theoretic image $\big\{f(x)\in \Si_2:x\in \Si_1\big\}$.
\label{lem:im}
\end{lemma}
\begin{proof} Let $f^*:\Hom(\La_1,\RR)\to \Hom(\La_2,\RR)$ denote the $\Hom$-dual of $f^\#$.
    The group-theoretic image of $f$ is the image of $f^*$ descended to the torus $\Sigma_2$. This image vector space is easily seen to be 
    \[
        \Im f^* = \big\{ \phi \in \Hom(\Lambda_2, \RR) : \phi(\Ker f^\#) = 0 \big\} 
        \cong{} \Hom \big(\Lambda_2 / \Ker f^\#, \RR \big).
    \]
    Descending this to $\Sigma_2$ amounts to taking the quotient by the lattice 
    \begin{align*}
        L = [-, \Lambda_2']_2 \cap \Im f^* = \big\{ \lambda_2' \in \Lambda_2' \bigmid [-, \lambda_2'] \text{ vanishes on } \Ker f^\# \big\}.
    \end{align*}
    For any $\lambda_1' \in \Lambda_1'$ the map $[-, f_\#(\lambda_1')]_2 = [f^\#(-), \lambda_1']_1$ clearly vanishes on $\Ker f^\#$ and hence we see that $[-, \Im f_\#]_2 \subseteq L$. Note that this is in fact a full rank sublattice of $L$:
    \[ \rk L = \dim f^* = \rk (\Lambda_2 / \Ker f^\#) = \rk f^\# = \rk f = \rk f_\# = \rk [-, \Im f_\#]_2  \]
    where the final equality holds since $[-,-]_2$ is a non-degenerate pairing. But this means that $\lambda_2' \in L$ if and only if $\lambda_2' \in \Lambda_2' \cap (\Im f_\#) \otimes \QQ = (\Im f_\#)^\sat$. In summary we see that the group-theoretic image of $f$ is the integral torus $\Im f$. \qedhere
\end{proof}

We carry on to $(\Ker f)_0$ and we claim that it is the connected component of the identity of the group-theoretic kernel, and that it is the kernel object in the category of integral tori.

\begin{proposition}
    Let $f:\Si_1\to \Si_2$ be a homomorphism of integral tori. 
    
    \begin{enumerate} \item The map $i:(\Ker f)_0\to \Sigma_1$ identifies the integral torus $(\Ker f)_0$ with the connected component of the identity of the group-theoretic kernel $\Ker f=\big\{x\in \Si_1:f(x)=0\big\}$. The group of connected components of $\Ker f$ is naturally identified with the quotient $f_{\#}(\La'_1)^{\sat}/f_{\#}(\La'_1)$.
    
    \item Given an integral torus $\Pi = (\Delta, \Delta', [\cdot, \cdot]_\Delta)$ and a homomorphism $g : \Pi \to \Sigma_1$ such that $f \circ g = 0$, there exists a unique homomorphism $u : \Pi \to (\Ker f)_0$ such that 
    \begin{equation*}
        \begin{tikzcd}
            (\Ker f)_0 \arrow[r,"i"] & \Sigma_1 \\
            \Pi \arrow[u, "u"] \arrow[ur, "g"'] & 
        \end{tikzcd}
    \end{equation*}
    commutes. 
    \end{enumerate}
    \label{prop:ker}
\end{proposition}

\begin{proof} 
    \begin{enumerate}
        \item We compute the group theoretic-kernel explicitly. As in the proof of Lemma~\ref{lem:image_of_tori}, let $f^*$ denote the $\Hom$-dual of $f^\#$. The kernel $\Ker f$ is the quotient of $(f^*)^{-1} \big([-, \Lambda_2']_2\big)$ by the lattice $[-, \Lambda_1']_1$. By basic linear algebra we know that if there exists a $\phi \in \Hom(\Lambda_1, \RR)$ such that $\phi \circ f^\# = [-, \lambda_2']_2$, then the fiber $(f^\ast)^{-1}\big([-, \lambda_2']_2 \big)$ is given by the coset $\phi + \Ker f^*$. 
        The set of all $[-, \lambda_2']_2$ which lie in the image of $f^*$ was determined in the proof of Lemma~\ref{lem:image_of_tori}: it is precisely $[-, (\Im f_\#)^\sat]_2$, which as a set is in bijection with $(\Im f_\#)^\sat$. 
        Using 
        \[ V = \Ker f^* \cong \Hom(\Coker f^\#, \RR) \cong \Hom\big((\Coker f^\#)^\torf, \RR \big), \]
        we see that the group-theoretic kernel of $f$ is isomorphic to $V \times (\Im f_\#)^\sat$ modulo $[-, \Lambda_1']_1$. 
        Under this quotient a connected component $V \times \{\lambda'\}$ gets identified with the image of $V \times \{0\}$ if and only if $[-, \lambda']_2$ is the $f^\ast$-image of a lattice point in $[-, \Lambda_1']_1$. This in turn is the case if and only if the defining point $\lambda'$ lies in $\Im f_\#$. 
        Finally note that $V \cap [-, \Lambda_1']_1 = [-, \Ker f_\#]_1$. Therefore we see that 
        \[ \Ker f = \underbrace{\frac{\Hom\big((\Coker f^\#)^\torf, \RR \big)}{[-, \Ker f_\#]_1}}_{= (\Ker f)_0} \times \frac{(\Im f_\#)^\sat}{\Im f_\#}. \]
        
        \item The assumption $f \circ g = 0$ means that $g^\# \circ f^\# = 0$ and $f_\# \circ g_\# = 0$. By the universal properties of kernels and cokernels of abelian groups, and using the fact that $\Delta$ is torsion free, we obtain unique morphisms $u^\#$ and $u_\#$ such that the diagrams
        \begin{equation*}
            \begin{tikzcd}
                (\Coker f^\#)^\mathrm{tf} \arrow[d, "u^\#"'] & \Lambda_1 \arrow[l,"i^\#"'] \arrow[dl, "g^\#"] \\
                \Delta & 
            \end{tikzcd}
            \qquad
            \text{and}
            \qquad
            \begin{tikzcd}
                \Ker f_\# \arrow[r,"i_\#"] & \Lambda_1' \\
                \Delta' \arrow[u, "u_\#"] \arrow[ur, "g_\#"'] & 
            \end{tikzcd}
        \end{equation*}
        commute. 
        We need to verify that $u = (u^\#, u_\#)$ is a homomorphism of integral tori. Let $\lambda \in \Lambda_1$ and $\delta' \in \Delta'$. Denote the class of $\lambda$ in $(\Coker f^\#)^\mathrm{tf}$ by $\overline{\lambda}$. Then
        \begin{equation*}
            \big[\overline{\lambda}, u_\#(\delta') \big]_K 
            = \big[\lambda, \underbrace{ { u_\#(\delta') } }_{\mathrlap{={} g_\#(\delta')}} \big]_1
            = \big[ \underbrace{g^\#(\lambda)}_{\mathrlap{= u^\#(\overline{\lambda})}}, \delta' \big]_\Delta ,
        \end{equation*}
        as required. \qedhere
    \end{enumerate}    
\end{proof}

The universal property for $\Coker f$ is stated and proved in complete analogy, and it is elementary to verify that
\[
(\Ker f)_0^\vee \cong \Coker (f^\vee).
\]
We may classify homomorphisms of integral tori according to two properties: the structure of the induced map on the underlying groups, and the dilation properties of the map of rational polyhedral spaces.

\begin{definition} \label{def:types_of_homomorphisms}
    Let $\Sigma_i = (\Lambda_i, \Lambda_i', [\cdot, \cdot]_i)$ for $i=1,2$ be integral tori of dimensions $g_i$. A homomorphism $f=(f^\#,f_\#): \Sigma_1 \to \Sigma_2$ of integral tori is said to be 
    \begin{enumerate}
        \item \emph{surjective} if $\rk f=g_2$ (equivalently, if $f^\#$ is injective), 
        \item \emph{finite} if $\rk f=g_1$ (equivalently, if $f_\#$ is injective),
        \item \emph{injective} if it is finite and $f_{\#}(\Lambda'_1)$ is saturated in $\Lambda'_2$, 
        \item an \emph{isogeny} if it is surjective and finite (equivalently, if $\rk f=g_1=g_2$),
        \item a \emph{free isogeny} if it is an isogeny and $f^\#(\La_2)=\La_1$ (equivalently, if $f^\#$ is an isomorphism),
        \item a \emph{dilation} if is an isogeny and injective (equivalently, if $f_\#$ is an isomorphism), and
        \item an \emph{isomorphism} if $f_\#$ and $f^\#$ are isomorphisms. 
    \end{enumerate}    
\end{definition}

In the sequence~\eqref{eq:kerimcoker}, the maps $i$ and $j$ are injective, while $p$ and $q$ are surjective. By Proposition~\ref{prop:ker} and Lemma~\ref{lem:im}, a homomorphism of integral tori $f$ is surjective, finite, injective, or a dilation if and only if the corresponding group homomorphism $f:\Si_1\to \Si_2$ is respectively surjective, has finite kernel (identified with the quotient $f_\#(\La'_1)^{\mathrm{sat}}/f_\#(\La'_1)$), injective, or is an isomorphism. %Also, $f$ is an isogeny if and only if $\Si_1$ and $\Si_2$ have the same dimension and the maps $f^\#$ and $f_\#$ have maximal rank.

For an isogeny $f:\Sigma_1\to \Sigma_2$, we can define several notions of degree. The image lattices $f_{\#}(\Lambda'_1)$ and $f^{\#}(\Lambda_2)$ have finite index in $\Lambda'_2$ and $\Lambda_1$, respectively, and we define the \emph{dilation degree} $\deg_d f$, the \emph{geometric degree} $\deg_g f$, and the \emph{total degree} $\deg f$ as
\[
\deg_d f=\big[\Lambda_1:f^\#(\Lambda_2)\big],\qquad \deg_g f =\big[\Lambda_2':f_\#(\Lambda_1')\big],\qquad \deg f=\deg_d f\cdot\deg_g f. 
\]
By Proposition~\ref{prop:ker}, the geometric degree is the order of the group-theoretic kernel of $f$, while the dilation degree is the factor by which $f$ stretches volume. All three degrees are multiplicative in compositions, and an isogeny $f$ is a free isogeny, a dilation, and an isomorphism if and only if respectively $\deg_d f=1$, $\deg_g f=1$, and $\deg f=1$. We note that the dual of a surjective homomorphism is finite (and vice versa), the dual of an isogeny is an isogeny (the geometric and dilation degrees are exchanged), and the dual of a free isogeny is a dilation (and vice versa). Finally, we may canonically factor any isogeny as follows:

\begin{lemma} 
    \label{lem:factoring_isogenies}
    Let $f=(f^\#,f_\#):(\Lambda_1, \Lambda_1', [\cdot, \cdot]_1)\to (\Lambda_2, \Lambda_2', [\cdot, \cdot]_2)$ be an isogeny of integral tori. Then there exists an integral torus $\Sigma_3$ together with a free isogeny $h : \Sigma_1 \to \Sigma_3$ and a dilation $g : \Sigma_3 \to \Sigma_2$ such that $f$ factors as $f = g \circ h$. Moreover, this factorization is unique in the sense that for any other factorization $\Sigma_1 \to \Pi \to \Sigma_2$ into a free isogeny followed by a dilation, there is a unique isomorphism $\phi$ of integral tori such that 
    \begin{equation*}
        \begin{tikzcd} [row sep = tiny]
            & \Sigma_3 \arrow[dr, "g"] & \\
            \Sigma_1 \arrow[ur, "h"] \arrow[dr] & & \Sigma_2 \\
            & \Pi \arrow[ur] \arrow[uu, "\phi"] & 
        \end{tikzcd}
    \end{equation*}
    commutes.
\end{lemma}

We note that we can also uniquely factor an isogeny as a dilation followed by a free isogeny. 

\begin{proof} 
    In order to factor $f$, we define a third integral torus $\Sigma_3 = (\Lambda_1, \Lambda_2', [\cdot, \cdot]_3)$. To define the pairing $[\cdot, \cdot]_3$ note the following. By the existence of a Smith normal form for $f_\#$, we may choose $\ZZ$-bases $e_1, \ldots, e_g$ of $\Lambda_1'$ and $f_1, \ldots, f_g$ of $\Lambda_2'$ such that $f_\#(e_i) = a_if_i$ for nonzero integers $a_i \in \ZZ$. We now define 
    \[ [\lambda, f_i]_3 = \frac{1}{a_i}[\lambda, e_i]_1 \]
    for any $\lambda \in \Lambda_1$ and factor $f$ as
    \begin{equation*}
        \begin{tikzcd} [row sep = tiny, column sep = huge]
            \Sigma_1 \arrow[r, "h"] & \Sigma_3 \arrow[r, "g"] & \Sigma_2 \\
            \Lambda_1 & \Lambda_1 \arrow[l, "\Id"'] & \Lambda_2 \arrow[l, "f^\#"'] \\
            \Lambda_1' \arrow[r, "f_\#"] & \Lambda_2' \arrow[r, "\Id"] & \Lambda_2'.
        \end{tikzcd}
    \end{equation*}
    It is easy to check Condition \eqref{eq:homomorphism_condition} for $h$ and $g$, hence they are homomorphisms. Moreover, $h$ is a free isogeny because $h^\# = \Id$ is an isomorphism, and $g$ is a dilation because $g_\# = \Id$ is surjective. 
    
    In order to show uniqueness, let $\Pi = (\Delta, \Delta', [\cdot, \cdot]_\Pi)$ be another integral torus and let $\Sigma_1 \overset{\tilde h}{\longrightarrow} \Pi \overset{\tilde g}{\longrightarrow} \Sigma_2$ be another factorization of $f$ as a free isogeny $\tilde h$ followed by a dilation $\tilde g$.
    In particular, $\tilde h^\#$ and $\tilde g_\#$ are isomorphisms.
    In order to define an isomorphism $\phi$ as in the claim, we need to ensure that $h = \phi \circ \tilde h$ and $\tilde g = g \circ \phi$, which means that the only candidate is $\phi = ((\tilde h^\#)^{-1}, \tilde g_\#)$. We verify that this is indeed a homomorphism, i.e.~we need to check Condition \eqref{eq:homomorphism_condition}. 
    Let $\delta_1, \ldots, \delta_g$ be the basis of $\Delta'$ such that $\tilde g_\#(\delta_i) = f_i$.
    Then $\tilde h_\#(e_i) = a_i\delta_i$ because $\Id \circ h_\# = \tilde g_\# \circ \tilde h_\#$ and it becomes clear that
    \[ \big[ (\tilde h^\#)^{-1}(\lambda), \delta_i \big]_\Pi 
        = \frac{1}{a_i} [\lambda, e_i]_1 
        = [\lambda, f_i]_3
        = \big[ \lambda, \tilde g_\#(\delta_i) \big]_3
    \]
    for all $\lambda \in \Lambda_1$ and $i = 1, \ldots, g$. \qedhere
\end{proof}

In some sense the key observation to prove Lemma \ref{lem:factoring_isogenies} was simply that an isogeny $f$ is free if and only if $f^\#$ is an isomorphism and it is a dilation if and only if $f_\#$ is an isomorphism.

\subsection{Polarizations and tropical homology} Let $\Sigma = (\Lambda, \Lambda', [\cdot, \cdot])$ be an integral torus. A \emph{polarization} on $\Sigma$ is a group homomorphism $\zeta : \Lambda' \to \Lambda$ such that $(\cdot,\cdot) = [\zeta(\cdot), \cdot] : \Lambda'_\RR \times \Lambda'_\RR \to \RR$ is a symmetric and positive definite bilinear form. A polarization is necessarily injective, and is called \emph{principal} if it is bijective. The invariant factors $(a_1,\ldots,a_g)$ of the Smith normal form (where $a_i\geq 1$ and $a_i|a_{i+1}$ for $i=1,\ldots,g-1$) define the \emph{type} of a polarization $\zeta$, and a polarization is principal if and only if all $a_i=1$. A polarization defines an isogeny $(\zeta, \zeta) : \Sigma \to \Sigma^\vee$ to the dual, which is an isomorphism if and only if the polarization is principal. 

Let $f : \Sigma_1 \to \Sigma_2$ be a finite homomorphism of integral tori. Given a polarization $\zeta_2$ on $\Si_2$, we define a polarization $f^*(\zeta_2)=f^\# \circ \zeta_2 \circ f_\#$ on $\Sigma_1$, called the \emph{induced polarization}. We note that a polarization induced from a principal polarization need not itself be principal.

\begin{definition}
    An integral torus together with a (principal) polarization is called a \emph{(principally) polarized tropical abelian variety} or \emph{(p)ptav} for short. 
    An \emph{isomorphism of ptavs} is an isomorphism $f : \Sigma_1 \to \Sigma_2$ of integral tori such that the polarization on $\Sigma_1$ is the polarization induced from $\Sigma_2$.    
\end{definition}

Finally, we recall the tropical homology and cohomology groups of an integral torus $\Sigma = (\Lambda, \Lambda', [\cdot, \cdot])$, and the relationship with polarizations. The groups can be computed explicitly (see \cite[Section 6]{GrossShokrieh_Poincare}):
\begin{align}
	H_{p, p}(\Sigma) &\cong \bigwedge^p \Lambda^*\otimes_\ZZ \bigwedge^p \Lambda',  \label{eq:torushomology}\\
	H^{p, p}(\Sigma) &\cong \bigwedge^p \Lambda \otimes_\ZZ \bigwedge^p(\Lambda')^* \label{eq:toruscohomology}.
\end{align}
Here $(\cdot)^*=\Hom(\cdot,\ZZ)$ denotes the dual lattice. In particular, $H^{1,1}(\Sigma) = \La\otimes(\La')^*=\Hom(\La',\La)$. Via this identification, we may view a polarization $\zeta:\La'\to \La$ as an element $\zeta\in H^{1,1}(\Sigma)$. Alternatively, we may define $\zeta=\cyc[\Theta]$, where $\Theta$ is the \emph{theta divisor} (see Section~\ref{sec:Cartier}, and see~\cite{MikhalkinZharkov} for the definition of the theta divisor in terms of tropical theta functions). By Poincar\'e duality, we may view the cup product $\zeta^p$ as an element of either $H^{p,p}(\Sigma)$ or $H_{g-p,g-p}(\Sigma)$, where $g$ is the dimension of $\Sigma$.

%\Dmitry{I would avoid calling this element $[\Theta]$, it should be either $\cyc [\Theta]$ or $c_1(\mathcal{L}(\Theta))$. I think that we need to carefully distinguish cycles and homology classes.}

We now state a homological criterion that allows us to check whether two pptavs are isomorphic, which is our principal reason for introducing tropical homology.

\begin{proposition} \label{prop:degree_pptavs_via_homology}
    Let $f : \Sigma_1 \to \Sigma_2$ be an isogeny of pptavs of dimension $g$ with principal polarizations $\zeta_i\in H^{1,1}(\Sigma_i)$. If
    \[ f_\ast (\zeta_1^{g-1}) = N \zeta_2^{g-1} \in H_{1,1}(\Sigma_2) \]
    for some integer $N\geq 1$, then the total degree of $f$ is equal to $\deg f = N^g$, and furthermore
    \[
    f^*(\zeta_2)=N\zeta_1\in H^{1,1}(\Sigma_1).
    \]
    In particular, if this condition holds for $N = 1$, then $f$ is an isomorphism of pptavs.
\label{prop:isomorphic_pptavs_via_homology}
\end{proposition}

\begin{proof}
    Via the Smith normal form construction we may pick bases $\Lambda_1' = \langle \lambda_1', \ldots, \lambda_g'\rangle$ and $\Lambda_2' = \langle \mu_1', \ldots, \mu_g' \rangle$ such that $f_\#(\lambda_i') = a_i \mu_i'$ for some integers $a_i$. Setting $\lambda_i = \zeta_1(\lambda_i')$ and $\mu_i = \zeta_2(\mu_i')$, we obtain corresponding bases of $\Lambda_1$ and $\Lambda_2$, respectively. With this notational setup we have
    \[
    \zeta_1=\sum_{i=1}^g\la_i\otimes (\la'_i)^*\in H^{1,1}(\Sigma_1).
    \]
    In terms of the explicit descriptions~\eqref{eq:torushomology}-\eqref{eq:toruscohomology}, the cup and cap products are given by the formulas 
    \[
    (\alpha\otimes \omega^*)\smile (\beta\otimes \xi^*) =(\alpha\wedge \beta)\otimes (\omega^*\wedge \xi^*),\qquad
    (\alpha^*\otimes \omega)\frown (\beta\otimes \xi^*)=(\al^* \lrcorner\beta)\otimes (\omega\lrcorner\xi^*),
    \]
    where $\lrcorner$ denotes the interior product on the exterior algebra \cite[Section~6]{GrossShokrieh_Poincare}. Therefore
    \[
    \zeta_1^{g-1}=(g-1)!\sum_{i=1}^g\left(\la_1\wedge \cdots\wedge \widehat{{\la_i}}\wedge \cdots\wedge\la_g\right)\otimes \left((\la_1')^*\wedge \cdots\wedge \widehat{(\la_i')^*}\wedge \cdots\wedge(\la_g')^*\right)\in H^{g-1,g-1}(\Sigma_1).
    \]
    The fundamental cycle of a pptav was computed in~\cite{GrossShokrieh_Poincare} (see Lemma 9.6, where the calculation was done for the Jacobian of a tropical curve, but the same argument works for any pptav): 
    \[
    \cyc[\Sigma_1]=\left(\la^*_1\wedge\cdots\wedge\la^*_g\right)\otimes \left(\la'_1\wedge\cdots\wedge\la'_g\right)\in H_{g,g}(\Sigma_1).
    \]
    Hence the Poincar\'e dual of $\zeta_1^{g-1}$ is
    \[
    \zeta_1^{g-1}=(g-1)!\sum_{i=1}^g\la^*_i\otimes \la'_i\in H_{1,1}(\Sigma_1).
    \]
    In terms of~\eqref{eq:torushomology}, the pushforward map $f_*:H_{1,1}(\Sigma_1)\to H_{1,1}(\Sigma_2)$ is $f_*=(f^\#)^*\otimes f_\#$, hence
    \[
    f_*(\zeta_1^{g-1})=(g-1)!\sum_{i=1}^g (f^\#)^*(\la^*_i)\otimes f_\#(\la'_i) =(g-1)!\sum_{i=1}^g (f^\#)^*(\la^*_i)\otimes a_i \mu'_i \in H_{1,1}(\Sigma_2).
    \]
    On the other hand, computing the Poincar\'e dual of $\zeta_2^{g-1}$ in the same way, we see that the condition  $f_*(\zeta_1^{g-1})=N \zeta_2^{g-1}$ implies that
    \[
    \sum_{i=1}^g (f^\#)^*(\la^*_i)\otimes a_i \mu'_i = N \sum_{i=1}^g \mu_i^*\otimes \mu'_i\in H_{1,1}(\Sigma_1).
    \]
    If follows that each $a_i$ divides $N$, that $f^\#$ is diagonalized by our choice of bases for $\Lambda_1$ and $\Lambda_2$, and that $f^\#(\mu_i) = \frac{N}{a_i} \la_i$. Hence the dilation, geometric, and total degrees of $f$ (see Definition~\ref{def:types_of_homomorphisms}) are equal to
    \[ \deg_d f = \prod_{i = 1}^g \frac{N}{a_i}, \qquad \deg_g f = \prod_{i = 1}^g a_i,\qquad \deg f=\deg_d f\cdot \deg_g f=N^g. \]
    Passing back to the dual lattices, the pullback map $f^*:H^{1,1}(\Sigma_2)\to H^{1,1}(\Sigma_1)$ is $f^*=f^\#\otimes (f_\#)^*$, hence $f^*(\zeta_2)=N\zeta_1$. \qedhere
\end{proof}

%\Dmitry{I've commented out the statement that uses pullback instead of pushforward}

%The following proposition then allows us to check homologically whether two principally polarized tropical abelian varieties are isomorphic.

%\begin{proposition}
%    \label{prop:isomorphic_pptavs_via_homology}
%    Let $\Sigma_i = (\Lambda_i, \Lambda_i', [\cdot, \cdot]_i)$ for $i = 1,2$ be tropical abelian varieties with principal polarizations $\xi_i$ repectively. Let $f : \Sigma_1 \to \Sigma_2$ be an isomorphism of integral tori. Then the following are equivalent.
%    \begin{enumerate}
%        \item $f$ is an isomorphism of pptavs.
%        \item The following diagram commutes
%        \begin{equation*}
%            \begin{tikzcd}
%                    \Lambda_1 & \Lambda_2 \ar[l, "f^\#"'] \\
%                    \Lambda_1' \ar[r, "f_\#"] \ar[u, "\xi_1"] & \Lambda_2' \ar[u, "\xi_2"]
%            \end{tikzcd}
%        \end{equation*}
%        \item The bilinear forms $[\xi_i(\cdot), \cdot]_i$ commute with the isomorphism $f_\# : \Lambda'_1 \to \Lambda'_2$. 
%        \item The induced isomorphism $f^*:H^{1,1}(\Sigma_2) \to H^{1,1}(\Sigma_1)$ identifies $\xi_2 \in H^{1,1}(\Sigma_2)$ with $\xi_1 \in H^{1,1}(\Sigma_1)$.
%    \end{enumerate}
%\end{proposition}

\subsection{Tropical Jacobians} Let $\Gamma$ be a tropical curve (which is assumed smooth as always). Given an oriented model $G$ of $\Ga$, the simplicial chain group $C_1(G,\ZZ)$ is the free abelian group on the edges of $G$, containing the simplicial homology group $H_1(G,\ZZ)$. These groups fit into a directed system with respect to refinements of models (see~\cite{baker2011metric} for details), and we denote the direct limit by $C_1(\Ga,\ZZ)$. The images of the $H_1(G,\ZZ)$ are all equal and are denoted $H_1(\Ga,\ZZ)$.

There is a natural isomorphism $d:H_1(\Ga,\ZZ)\to \Omega^1_{\Ga}(\Ga)$ sending a cycle $\sum a_e\, e$ to the $1$-form $\sum a_e \,de$. In addition, the \emph{integration pairing}
\begin{align*}
    [\cdot, \cdot] : \Omega^1_\Ga(\Ga)\times C_1(\Gamma, \ZZ)   &\longrightarrow \RR \\
    (\omega, \gamma) = \left(\sum a_e de,\sum b_e e\right) &\longmapsto \int_\gamma \omega=\sum a_eb_e\ell(e)
\end{align*}
restricts to a non-degenerate pairing $\Omega^1_\Ga(\Ga)\times H_1(\Gamma, \ZZ) \to \RR$.  Hence we have a pptav 
\[ \Jac(\Gamma) = \Big(\Omega^1_\Ga(\Ga), \quad H_1(\Gamma, \ZZ), \quad [\cdot, \cdot]\Big) =\Hom\big(\Om^1_{\Ga}(\Ga),\RR\big)/H_1(\Ga,\ZZ) \]
of dimension equal to the genus $g(\Ga)$, the \emph{tropical Jacobian variety} of $\Gamma$. 

\begin{definition}
    Fix a base point $q \in \Gamma$. The \emph{Abel--Jacobi map} relative to $q$ is given by
    \begin{align*}
        \phi_q : \Gamma &\longrightarrow \Jac(\Gamma) \\
        p &\longmapsto \Big( \omega \mapsto \int_{\gamma_p} \omega \Big) ,
    \end{align*}
    where $\gamma_p$ denotes any path from $q$ to $p$. 
\end{definition}

The Abel--Jacobi map naturally extends to symmetric powers of $\Ga$ and hence to divisors. This map respects linear equivalence, and the tropical Abel--Jacobi theorem (see Theorem 6.3 in~\cite{MikhalkinZharkov}) states that the induced map $\Pic_0(\Ga)\to \Jac(\Ga)$ (which does not depend on the choice of a base point) is an isomorphism. Under this identification, the Abel--Jacobi map can also be described as $p\mapsto p-q$.

We now prove that the Abel--Jacobi map enjoys a universal property among morphisms of rational polyhedral spaces to integral tori, which is an exact analogue of the algebraic property (see Proposition~11.4.1 in~\cite{BirkenhakeLange}, and see Section 1.4 in~\cite{baker2007riemann} for the corresponding property of the Jacobian of a finite graph). We first make the following elementary observation, which does not appear to have a proof in the literature.

\begin{proposition} \label{prop:abel_jacobi_rpc}
    The Abel--Jacobi map $\phi_q : \Gamma \to \Jac(\Gamma)$ is a map of rational polyhedral spaces. Moreover, the identification $\Omega_{\Jac (\Gamma)}^1(\Jac (\Ga)) \cong \Omega_\Gamma^1(\Ga)$ (see Remark~\ref{rem:affine_fcn_on_tori}) is induced by pullback along $\phi_q$.
\end{proposition}

\begin{proof}
    Let $\eta \in \Omega^1_{\Jac(\Gamma)}(\Jac(\Ga))$ be a 1-form.  We need to show that its pullback along $\phi_q$ is a 1-form on $\Gamma$. Recall that $\Omega^1_{\Jac (\Gamma)}(\Jac (\Ga)) \cong \Omega^1_{\Gamma}(\Ga)$ and denote the 1-form on $\Gamma$ that corresponds to $\eta$ under this identification by $\omega = \sum_{e\in E(\Gamma)} a_e d e$ with $a_e \in \ZZ$. We show that the pullback of $\eta$ along $\phi_q$ is $\omega$. 
    Thinking of $\eta$ as a linear function on $\Hom(\Omega^1_\Gamma(\Ga), \RR)$, we easily see that
    \[ \eta \big(\phi_q(p) \big) = \eta \Big(\int_{\gamma_p} - \Big) = \int_{\gamma_p} \omega \mod H_1(\Ga,\ZZ)\]
    for any $p \in \Gamma$ and $\gamma_p$ a path from $q$ to $p$. The second equality is simply the identification of $\eta \in \Hom\big(\Hom(\Omega^1_\Gamma(\Ga), \RR), \RR\big)$ with $\omega \in \Omega^1_\Gamma(\Ga)$. But now we are already done because the coefficients $a_e$ of $\omega$ are precisely $\frac{1}{\ell(e)}\int_e \omega$, with $e$ parametrized with the orientation indicated by $de$. \qedhere
\end{proof}

\begin{proposition} \label{prop:universal_property_AbelJacobi_map}
    Let $\chi:\Gamma\to X$ be a morphism of rational polyhedral spaces from a tropical curve $\Gamma$ to an integral torus $X=(\La,\La',[\cdot,\cdot])$. Then there exists a unique homomorphism $\mu:\Jac(\Gamma)\to X$ of integral tori such that the diagram
    \begin{equation} \label{eq:univ_property_Jac}
        \begin{tikzcd}
            \Gamma \arrow[r,"\chi"] \arrow[d,"\phi_q"] & X \arrow[d,"t_{-\chi(q)}"]\\
            \Jac(\Gamma) \arrow[r, "\mu"] & X 
        \end{tikzcd}
    \end{equation}
    commutes for all $q\in \Gamma$. 
\end{proposition}

\begin{proof} Using the Abel--Jacobi isomorphism $\Pic_0(\Ga)=\Jac(\Ga)$, we could simply define $\mu$ by the rule
\[
\mu\left(\sum a_i p_i\right)= \sum a_i \chi(p_i).
\]
However, we want to carefully define $\mu$ as a homomorphism of integral tori.

    Fix $q \in \Gamma$. Clearly, the composition $\chi_0 = t_{-\chi(q)} \circ \chi$ maps $q$ to 0. We show that $\chi_0$ factors via the Abel--Jacobi map. To do so, we need to describe maps $\mu^\# : \Lambda \to \Omega_\Ga^1(\Ga)$ and $\mu_\# : H_1(\Gamma, \ZZ) \to \Lambda'$ which are compatible with the pairings on the integral tori such that $\mu = (\mu^\#,\mu_\#)$ makes Diagram~\eqref{eq:univ_property_Jac} commute.
    
    The morphism $\chi_0$ of rational polyhedral spaces induces a morphism of cotangent sheaves
    \[ \chi_0^\ast : (\chi_0)^{-1}\Omega_X^1 \longrightarrow \Omega_\Gamma^1 . \]
    We recall from Remark~\ref{rem:affine_fcn_on_tori} that $\Omega_X^1(X) \cong \Lambda$, and recall from Proposition~\ref{prop:abel_jacobi_rpc} that $\phi_q^\ast$ is the isomorphism $\Omega_{\Jac (\Gamma)}^1(\Jac(\Ga)) \cong \Omega_\Gamma^1(\Ga)$. Hence, the only possible choice for $\mu^\#$ which ensures commutativity of Diagram~\eqref{eq:univ_property_Jac} is $\mu^\# =   (\phi_q^\ast)^{-1}\circ \chi_0^\ast$. Passing to singular homology, the continuous map $\chi_0$ induces a pushforward map $H_1(\Gamma, \ZZ) \to H_1(X, \ZZ)$. Identifying $H_1(X, \ZZ) = \Lambda'$, we let $\mu_\#$ be the pushforward map. To show that the pair $\mu = (\mu^\#, \mu_\#)$ defines a homomorphism of integral tori we need to verify compatibility with pairings, in other words we need to show that 
    \begin{equation} \label{eq:compatibility_univ_property_Jac}
        \big[\lambda, \mu_\#(\gamma)\big] = \int_\gamma \mu^\#\lambda
    \end{equation}
    for all $\lambda \in \Lambda$ and $\gamma \in H_1(\Gamma, \ZZ)$.
    
    We now show that the pairing $[ \cdot, \cdot]$ on $X$ can be interpreted as integration as well. Indeed, let $\lambda \in \Lambda = \Omega^1_X(X)$, which we view as the differential of a locally well-defined affine linear function $f$ on $X$. Choose the integration constant so that $f(0)=0$, then $f$ is linear and can be extended to a globally well-defined linear function on the universal cover $\Hom(\Lambda, \RR)$, namely
    \begin{align*}
        F :\Hom(\Lambda, \RR) &\longrightarrow \RR \\
        u &\longmapsto u(\lambda) .
    \end{align*}
    Let $\gamma \in H_1(X, \ZZ) = \Lambda'$. Choose a piecewise smooth representative and lift $\ga$ to a path $\gamma' : [0, 1] \to \Hom(\Lambda, \RR)$ on the universal cover of $X$ going from 0 to some point $\lambda' \in \Lambda' \subseteq \Hom(\Lambda, \RR)$.
    Then $\int_\gamma \lambda = \int_{\gamma'} \lambda = F(\gamma'(1)) - F(\gamma'(0)) = F(\lambda') = \lambda'(\lambda)$. 
    But now recall that $\Lambda'$ is embedded in $\Hom(\Lambda, \RR)$ as $[\cdot, \lambda']$. This shows that $[\cdot, \cdot]$ is just integration of 1-forms along closed paths, and Equation~\eqref{eq:compatibility_univ_property_Jac} is simply the change-of-variables formula for line integrals. Hence $\mu$ is a homomorphism of integral tori.

    Finally, we show that $\mu$ makes the diagram~\eqref{eq:univ_property_Jac} commute for any $q' \in \Gamma$ (and not just the $q$ that we fixed at the beginning of the proof). Indeed, it is clear that $\phi_{q'} = t_{-\phi_{q}(q')} \circ \phi_{q}$. 
    By Lemma~\ref{lem:swap_homomorphism_translation} we obtain 
    \begin{equation} \label{eq:diagram_q_to_q'}
        \begin{tikzcd}
            \Jac(\Gamma) \arrow[r, "\mu"] \arrow[d, "t_{-\phi_{q}(q')}"'] & X \arrow[d, "t_{\mu(-\phi_{q}(q'))}"] \\
            \Jac(\Gamma) \arrow[r, "\mu"] & X
        \end{tikzcd}
    \end{equation}
    But now $\mu\big(-\phi_{q}(q')\big) = -\mu\big(\phi_q(q')\big) = -t_{-\chi(q)}\big(\chi(q')\big) = \chi(q) - \chi(q')$, where the second equality uses the already established commutativity of Diagram~\eqref{eq:univ_property_Jac} for $q$. Combining Diagram~\eqref{eq:univ_property_Jac} for $q$ with Diagram~\eqref{eq:diagram_q_to_q'}, we obtain the claim for $q'$ and hence we are done.

    The uniqueness of $\mu$ follows from the fact that the image of $\Ga$ under the Abel--Jacobi map $\phi_q$ spans $\Jac(\Ga)$ as a torus, and hence any two homomorphisms on $\Jac(\Ga)$ that agree on the image of $\Ga$ are equal.  \qedhere
\end{proof}

\begin{remark}
    For the proof of Theorem \ref{thm:tropical_Recillas_theorem} it is essential that all tropical curves be smooth and hence carry a fundamental cycle. Given a metric graph with finite univalent vertices, this is achieved by adding compact infinite rays to such vertices (see Example~\ref{ex:fundamental_cycle}). We emphasize that the universal property in Proposition~\ref{prop:universal_property_AbelJacobi_map} is still valid in this context: any morphism of polyhedral spaces $\Gamma \to X$ maps each infinite ray to a single point in $X$, because any affine linear function on the ray is eventually constant.
\end{remark}

\subsection{Tropical Prym variety} In this section, we recall the definition of the tropical Prym variety $\Prym(\tGa/\Ga)$ of a harmonic double cover $\pi : \tilde \Gamma \to \Gamma$ of tropical curves. The tropical Prym variety was defined in~\cite{JensenLen} and further studied in~\cite{LenUlirsch} and~\cite{LenZakharov}.  We also define an alternative object, the \emph{continuous Prym variety} $\Prym_c(\tGa/\Ga)$, which is more naturally suited to our purposes, and investigate its relationship with $\Prym(\tGa/\Ga)$, which we henceforth call the \emph{divisorial Prym variety} and denote $\Prym_d(\tGa/\Ga)$.

%which we adopt slightly for our purposes.  Recall that an edge or vertex of $\Ga$ is called~\emph{free} if it has two preimages in $\tGa$ each of which has dilation factor equal to $1$, and~\emph{dilated} if it has a unique preimage with dilation factor equal to $2$. The set of dilated edges and vertices form the~\emph{dilation subgraph} of $\Gamma$. We say that $\pi$ is~\emph{free} if the dilation subgraph is empty and~\emph{dilated} otherwise.

We assume a choice of graph model for $\pi$. Recall from Section~\ref{sec:graphs} that an edge or vertex of $\Ga$ is called~\emph{free} if it has two preimages in $\tGa$ each of which has dilation factor equal to $1$, and~\emph{dilated} if it has a unique preimage with dilation factor equal to $2$. The set of dilated edges and vertices form the~\emph{dilation subgraph} of $\Gamma_{\dil}\subseteq\Gamma$. We say that $\pi$ is~\emph{free} if the dilation subgraph is empty and~\emph{dilated} otherwise. The \emph{dilation index} of the double cover $\pi:\tGa\to \Ga$ is
\[
d(\tGa/\Ga)=\begin{cases} \mbox{number of connected components of }\Ga_{\dil},& \mbox{if $\pi$ is dilated,}\\
1, & \mbox{if $\pi$ is free,}
\end{cases}
\]
where we note that a free double cover has dilation index $1$, not the expected $0$. 

\begin{remark} The tropicalization of an algebraic \'etale double cover has the additional property of being~\emph{unramified}. This condition involves vertex weights, which we do not use, and also imposes a restriction on the dilation subgraph: each vertex must have even valence (see \cite[Corollary~5.5]{JensenLen}). This restriction does not naturally arise in the tropical setting, and we do not impose it. 
\end{remark}

A free edge $e$ of $\Ga$ has two distinct preimages that we arbitrarily label $\te^+$ and $\te^-$, while a dilated edge $e$ has a unique preimage that we denote $\te^+=\te^-$ by abuse of notation. The double cover $\pi:\tGa\to\Ga$ has an associated involution $\iota:\tGa\to\tGa$ defined by $\iota(\te^{\pm})=\te^{\mp}$, in other words $\iota$ exchanges the preimages of a free edge and fixes the preimage of a dilated edge. The double cover $\pi$ and the involution $\iota$ induce maps 
\begin{equation}
    \begin{aligned}
        &\pi^\ast :& \Omega^1_\Gamma(\Ga) &\longrightarrow \Omega^1_{\tilde \Gamma}(\tGa) \\
        &&\sum a_e de &\longmapsto \sum a_e(d \tilde e^+ + d \tilde e^-)\\[1em]
        &\pi_\ast :& H_1(\tilde \Gamma, \ZZ) &\longrightarrow H_1(\Gamma, \ZZ) \\
        &&\sum a_{\tilde e^{\pm}}\tilde e^{\pm} &\longmapsto\sum (a_{\tilde e^+} + a_{\tilde e^-}) e \\
    \end{aligned}
    \qquad \text{and} \qquad
    \begin{aligned}
        &\iota^*:&\Omega^1_{\tilde \Gamma}(\tGa) & \longrightarrow  \Omega^1_{\tilde \Gamma}(\tGa) \\
        &&\sum a_{\tilde e^{\pm}}d\tilde e^{\pm} &\longmapsto\sum a_{\tilde e^{\pm}}d\tilde e^{\mp}\\[1em]
        &\iota_*:&H_1(\tilde \Gamma, \ZZ) & \longrightarrow H_1(\tilde \Gamma, \ZZ)\\
        &&\sum a_{\tilde e^{\pm}}\tilde e^{\pm} &\longmapsto\sum a_{\tilde e^{\pm}}\tilde e^{\mp}\\
    \end{aligned}
    \label{eq:pushpull}
\end{equation}
We have two associated maps on the Jacobians. The pair of maps $(\pi^*,\pi_*)$ defines the \emph{tropical norm homomorphism}  $\pi: \Jac(\tilde \Gamma) \to \Jac(\Gamma)$. Under the Abel--Jacobi identifications $\Pic_0(\tGa)=\Jac(\tGa)$ and $\Pic_0(\Ga)=\Jac(\Ga)$, the norm homomorphism is the pushforward map on divisors. The involution $\iota:\tGa\to \tGa$ induces an involution $\iota=(\iota^*,\iota_*):\Jac(\tGa)\to \Jac(\tGa)$, which again is simply the involution acting on divisors, and the composition $\pi\circ (\Id-\iota):\Jac(\tGa)\to \Jac(\Ga)$ is the zero map. 
%\Dmitry{I decided to denote the norm homomorphism by $\pi$} \Felix{I can see why, but I don't think the gains justify changing it everywhere.}

We now define the two Prym varieties associated to the double cover $\pi:\tGa\to \Ga$.

\begin{definition} The \emph{divisorial Prym variety} of the double cover $\pi:\tGa\to\Ga$ is the connected component of the identity of the kernel of $\pi$ (see~\eqref{eq:kercoker} and Proposition~\ref{prop:ker}):
\begin{equation}
\Prym_d(\tilde \Gamma / \Gamma) = (\Ker \pi)_0= \Big( (\Coker \pi^*)^\mathrm{tf}, \quad \Ker \pi_*, \quad [\cdot, \cdot]_P \Big).
\label{eq:Prymd}
\end{equation}
The \emph{continuous Prym variety} of the double cover $\pi:\tGa\to\Ga$ is the integral torus
\begin{equation}
\Prym_c(\tilde \Gamma/ \Gamma)=\Big( \Omega^1_{\tGa}(\tGa)/\Ker (\Id-\iota^*), \quad \Im (\Id-\iota_*), \quad [\cdot, \cdot]_P \Big).
\label{eq:Prymc}
\end{equation}
The pairing $[\cdot,\cdot]_P$ on both tori is induced by the integration pairing on $\Jac(\tGa)$. 
\end{definition}

Under the Abel--Jacobi identifications, the divisorial Prym variety is the connected component of the identity of the norm homomorphism $\pi_*:\Pic_0(\tGa)\to\Pic_0(\Ga)$. It is the obvious tropical analogue of the algebraic Prym variety and is the object studied in~\cite{JensenLen} and~\cite{LenUlirsch}. The continuous Prym variety, to the best of our knowledge, has not been considered before, and is not naturally a subset of the Jacobian (however, it is in some sense a tropical analogue of the variety constructed in~\cite[Proposition 12.1.8]{BirkenhakeLange}). We summarize the properties of the two Pryms and their relationship below. 

%\Dmitry{Do we say anything about chip-firing?}

%nor does it have an elementary description in terms of divisors and chip-firing \Felix{Are you sure? Shouldn't this just be the locus of anti-symmetric divisors $D + \iota D$ in $\Pic^0(\tilde \Gamma)$ with $\deg D$ even?}. \Dmitry{I'm not sure. $\Prym_c(\tGa/\Ga)$ does not even naturally sit inside $\Pic^0(\tGa)$}  

We first determine the relationship between the two Pryms and $\Jac(\tGa/\Ga)$. Consider the homomorphism $\Id-\iota:\Jac(\tGa)\to \Jac(\tGa)$. Looking at~\eqref{eq:im}, we see that $\Prym_c(\tGa/\Ga)$ is defined in the same way as the image torus $\Im(\Id-\iota)$, except that we do not saturate the second lattice. In particular, the surjective map $\Jac(\tGa)\to \Im(\Id-\iota)$ factors through a free isogeny $\Prym_c(\tGa/\Ga)\to \Im(\Id-\iota)$, whose geometric degree is the index of the lattice $\Im (\Id-\iota_*)$ in its saturation. On the other hand, $\pi\circ(\Id-\iota)=0$ and the universal property of the kernel (see Proposition~\ref{prop:ker}) implies that the map $\Im(\Id-\iota)\to \Jac(\tGa)$ factors through the kernel $\Prym_d(\tGa/\Ga)$. We now show that the latter is in fact the image of $\Id-\iota$, compute the geometric degree of the free isogeny, and determine the induced polarizations.

\begin{proposition} Let $\pi:\tGa\to \Ga$ be a double cover with dilation index $d(\tGa/\Ga)$.
\begin{enumerate}
    \item The group-theoretic kernel $\Ker \big(\pi:\Jac(\tGa)\to\Jac(\Ga)\big)$ has two connected components if $\pi$ is free and one if $\pi$ is dilated (see~\cite[Proposition 6.1]{JensenLen}).
    \item The divisorial Prym variety $\Prym_d(\tGa/\Ga)$ is the group-theoretic image of $\Id-\iota$:
\[
\Prym_d(\tilde \Gamma/ \Gamma)=\Im(\Id-\iota)=\Big( \Omega^1_{\tGa}(\tGa)/\Ker (\Id-\iota^*), \quad (\Im (\Id-\iota_*))^{\sat}, \quad [\cdot, \cdot]_K \Big).
\]    
Hence the map $\Id-\iota$ factors as
\begin{equation}
\Jac(\tGa)\xlongrightarrow{\epsilon} \Prym_c(\tGa / \Ga )\xlongrightarrow{\gamma} \Prym_d(\tGa / \Ga) \overset{i}{\hooklongrightarrow} \Jac(\tGa),
\label{eq:Prymmaps}
\end{equation}
and the middle map $\gamma:\Prym_c(\tGa /\Ga)\to \Prym_d(\tGa / \Ga)$ is a free isogeny of geometric degree $2^{d(\tGa/\Ga)-1}$ (and dilation degree 1). In particular, $\Prym_c(\tGa/\Ga)=\Prym_d(\tGa/\Ga)$ if $d(\tGa/\Ga)=1$, in other words if $\pi$ is free or if the dilation subgraph is connected. 

\item The polarization on $\Prym_d(\tGa/\Ga)$ induced from the principal polarization on $\Jac(\tGa)$ via $i$ has type $(1,\ldots,1,2,\ldots,2)$, where the number of $1$'s equals $d(\tGa/\Ga)-1$. 

\item There exists a principal polarization $\zeta_c$ on $\Prym_c(\tGa/\Ga)$ such that 
$\zeta=2\zeta_c$ is the polarization induced from the principal polarization on $\Jac(\tGa)$ via $i \circ \gamma$.
\end{enumerate}
\label{prop:Prymallpropetries}

\end{proposition}

The proof of this proposition is rather laborious, and involves finding formulas for the pushforward and pullback maps~\eqref{eq:pushpull} in terms of explicit bases. Before we give the proof, we propose three reasons why it is more natural to consider $\Prym_c(\tGa/\Ga)$ than $\Prym_d(\tGa / \Ga)$. 

\begin{enumerate}
    \item The induced polarization on $\Prym_d(\tGa/\Ga)$ is twice a principal polarization if and only if $d(\tGa/\Ga)=1$, in other words if $\pi$ is free or has a connected dilation subgraph. The induced polarization on $\Prym_c(\tGa/\Ga)$, however, is always twice a principal polarization, just as in the algebraic setting. 
    \item The divisorial Prym variety $\Prym_d(\tGa/\Ga)$ does not behave continuously under contractions of the edges of $\Ga$, specifically those that change the dilation index. This behavior was explored in detail in~\cite{2024GhoshZakharov}. Hence $\Prym_d(\tGa/\Ga)$ is unsuitable from a moduli-theoretic viewpoint. We do not explore this behavior here, but simply state without proof that $\Prym_c(\tGa/\Ga)$ is in fact continuous in families.
    \item The continuous Prym variety $\Prym_c(\tGa/\Ga)$ satisfies a natural universal property (see Proposition~\ref{prop:universal_property_Prymc}), while $\Prym_d(\tGa/\Ga)$ does not. 
\end{enumerate}

To simplify notation, we temporarily identify $H_1(\tGa,\ZZ)=\Om^1_{\tGa}(\tGa)$ and $H_1(\Ga,\ZZ)=\Om^1_{\Ga}(\Ga)$ using the principal polarizations. In terms of this identification, we have $\iota^*=\iota_*$ and there is a pullback map $\pi^*:H_1(\Ga,\ZZ)\to H_1(\tGa,\ZZ)$ on homology.  It is easy to verify that
\[
\pi^*\circ \pi_*=\Id+\iota_*.
\]
For a free double cover $\pi:\tGa\to \Ga$ we have $g(\tGa)=2g(\Ga)-1$ and therefore the dimension of the Prym is $g_0=g(\Ga)-1$. For a dilated double cover, we denote by $\tGa_{\dil}\subset \tGa$ the isomorphic preimage of the dilation subgraph $\Ga_{\dil}\subset \Ga$, and recall that $d(\tGa/\Ga)$ is the number of connected components of $\Ga_{\dil}$. Let $m_{\dil}=|E(\Ga_{\dil})|$ and $n_{\dil}=|V(\Ga_{\dil})|$ denote the number of dilated edges and dilated vertices, respectively. The numbers $m_{\dil}$ and $n_{\dil}$ depend on the choice of model but their difference does not, and we introduce the invariants
\begin{equation} \label{eq:ABC}
\begin{aligned}
    A &= g(\Ga)-m_{\dil}+n_{\dil}-d(\tGa/\Ga), \\
    B &= d(\tGa/\Ga)-1, \\
    C &= m_{\dil}-n_{\dil}+d(\tGa/\Ga).
\end{aligned} 
\end{equation}
It is easy to see that $|E(\tGa)|=2|E(\Ga)|-m_{\dil}$ and $|V(\tGa)|=2|V(\Ga)|-n_{\dil}$, therefore
\[
A+B=g(\Ga)-m_{\dil}+n_{\dil}-1=|E(\Ga)|-|V(\Ga)|-m_{\dil}+n_{\dil}=g(\tGa)-g(\Ga)=g_0.
\]
We first consider the case of free double covers.
%\Dmitry{what do we need from this paragraph?} In order to describe the polarization type of $(\Ker \Nm)_0$ associated to a double cover $\pi:\tGa\to \Ga$, and for other explicit calculations, we construct explicit bases for the homology groups $H_1(\tGa,\ZZ)$ and $H_1(\Ga,\ZZ)$. Propositions~\ref{prop:freebasis} and~\ref{prop:dilatedbasis} sharpen and improve on Lemma~1.5.4 in~\cite{LenUlirsch}\footnote{In our opinion, there are certain gaps in the proof of Lemma~1.5.4 in~\cite{LenUlirsch}.}.

\begin{proposition} 
    Let $\pi:\tGa\to \Ga$ be a connected free double cover and let $g(\Ga)=g$ and $g(\tGa)=2g-1$. Then there exists a basis $\al_1,\ldots,\al_{g-1},\ga_1$ of $H_1(\Ga,\ZZ)$ and a basis $\tal^{\pm}_1,\ldots,\tal^{\pm}_{g-1},\tga_1$ of $H_1(\tGa,\ZZ)$ such that
        \begin{align*}
        \iota_*(\tal^{\pm}_i) &= \tal^{\mp}_i, & \pi_*(\tal^{\pm}_i) &= \al_i, & 
        \pi^*(\al_i) &= \tal^+_i+\tal^-_i, & i=1,\ldots,g-1, \\
        \iota_*(\tga_1) &= \tga_1, & \pi_*(\tga_1) &= 2\ga_1, & \pi^*(\ga_1) &= \tga_1.
    \end{align*}
    
    \label{prop:freebasis}
\end{proposition}

\begin{proof} The proof of this statement is included as part of Constructions A and B in~\cite{LenZakharov}, and we briefly summarize it. Choose an orientation on $\Ga$ and a spanning tree $T\subset \Ga$, and denote the complementary edges by $E(\Ga)\backslash E(T)=\{e_0,\ldots,e_{g-1}\}$. Let $\te^{\pm}_i$ and $\tT^{\pm}\subset \tGa$ denote the preimages of $e_i$ and $T$, respectively. Let $S\subset\{e_0,\ldots,e_{g-1}\}$ denote the set of those complementary edges whose lifts connect the two trees $\tT^{\pm}$. The set $S$ is nonempty since $\tGa$ is connected, so we assume without loss of generality that $e_0\in S$. It follows that $\tT=\tT^+\cup \tT^-\cup\{\te^+_0\}$ is a spanning tree for $\tGa$.

For a cycle $\ga\in H_1(\Ga,\ZZ)$ and an oriented edge $e\in E(\Ga)$, denote by $\langle\ga,e\rangle$ the coefficient with which $e$ appears in $\ga$, and similarly for $H_1(\tGa,\ZZ)$. We denote $\ep_i\in H_1(\Ga,\ZZ)$ for $i=0,\ldots,g-1$ the unique cycle on the graph $T\cup\{e_i\}$ such that $\langle \ep_i,e_i\rangle=1$ Similarly, let $\tep_0\in H_1(\tGa,\ZZ)$ and $\tep^{\pm}_i\in H_1(\tGa,\ZZ)$ for $i=1,\ldots,g-1$ denote the unique cycles on $\tT\cup \{\te^-_0\}$ and $\tT\cup \{\te^{\pm}_i\}$ such that $\langle \tep_0,\te^-_0\rangle=1$ and $\langle \tep^{\pm}_i,\te^{\pm}_i\rangle=1$, respectively. The cycles $\ep_0,\ldots,\ep_{g-1}$ form a basis for $H_1(\Ga,\ZZ)$, and furthermore the coordinates of any $\ga\in H_1(\Ga,\ZZ)$ with respect to this basis are given by
\[
\ga=\langle \ga,e_0\rangle \ep_0+\cdots+\langle \ga,e_{g-1}\rangle\ep_{g-1},
\]
and a similar statement holds for $\tep_0,\tep^{\pm}_1,\ldots,\tep^{\pm}_{g-1}$.

The action of $\iota_*$, $\pi_*$, and $\pi^*$ on these bases is computed by looking at the coefficients of the edges $e_i$ and $\te^{\pm}_i$. Since $\langle \tep_0,\te^-_0\rangle=1$, we see that $\iota_*(\tep_0)=\tep_0$, $\pi_*(\tep_0)=2\ep_0$, and $\pi^*(\ep_0)=\tep_0$. Now denote $c_i=\langle \tep^+_i,\te^+_0\rangle$ for $i=1,\ldots,g-1$  (this number is equal to $0$ or $\pm 1$ since $\tep^+_i$ is a simple cycle). Comparing the coefficients of $\te^-_0$ and $\te^-_i$, we see that $\iota_*(\tep^+_i)=\tep^-_i+c_i\tep_0$ and $\iota_*(\tep^-_i)=\tep^+_i-c_i\tep_0$. Similarly, comparing the coefficients of $e_0$ and $e_i$, we see that $\pi_*(\tep^{\pm}_i)=\ep_i\pm c_i\ep_0$. Finally, using the relation $\pi^*\circ\pi_*=\Id+\iota_*$ we find that
\[
\pi^*(2\ep_0)=\pi^*(\pi_*(\tep_0))=(\Id+\iota_*)(\tep_0)=2\tep_0
\]
and
\[
\pi^*(2\ep_i)=\pi^*\pi_*(2\tep^+_i-c_i\tep_0)=(\Id+\iota_*)(2\tep^+_i-c_i\tep_0)=2(\tep^+_i+\tep^-_i),
\]
for $i=1,\ldots,g-1$. To complete the proof, we now set for $i=1,\ldots,g-1$
\begin{align*}
    \tal^+_i &= \tep^+_i, \\
    \tal^-_i &= \tep^-_i+c_i\tep_0,\\
    \al_i &= \ep_i+c_i\ep_0,
\end{align*}
and $\tga_1=\tep_0$, $\ga_1=\ep_0$. \qedhere
\end{proof}

We now consider the dilated case.

\begin{proposition} \label{prop:dilatedbasis} Let $\pi:\tGa\to\Ga$ be a dilated double cover. Then there exists a basis $\al_1,\ldots,\al_A$, $\ga_1,\ldots,\ga_C$
    of $H_1(\Ga,\ZZ)$ and a basis $\tal^{\pm}_1,\ldots,\tal^{\pm}_A$,  $\tbe_1,\ldots,\tbe_B$, $\tga_1,\ldots,\tga_C$
    of $H_1(\tGa,\ZZ)$ with $A$, $B$, and $C$ as defined in Equation~\eqref{eq:ABC}, such that
    \begin{align*}
        \iota_*(\tal^{\pm}_i) &= \tal^{\mp}_i,
        & \pi_*(\tal^{\pm}_i) &= \al_i,
        & \pi^*(\al_i) &= \tal^+_i+\tal^-_i,
        & i&=1,\ldots,A, \\
        \iota_*(\tbe_j) &= -\tbe_j,
        & \pi_*(\tbe_j) &= 0,
        &&& j&=1,\ldots,B,\\
        \iota_*(\tga_k) &= \tga_k,
        & \pi_*(\tga_k) &= \ga_k,
        & \pi^*(\ga_k) &= 2\tga_k,
        & k&=1,\ldots,C.
    \end{align*}
\end{proposition}

\begin{proof} Let $d = d(\tGa /\Ga)$ be the dilation index of $\pi$. We begin by contracting each dilated edge of $\Ga$ and the corresponding edge of $\tGa$. The result is a double cover $\pi':\tGa'\to \Ga'$ with associated involution $\iota':\tGa'\to \tGa'$, whose dilated vertices correspond to the connected components of $\Ga_{\dil}$. Denote these vertices by $v'_0,\ldots,v'_{d-1}\in V(\Ga')$ and their preimages by $\tv'_0,\ldots,\tv'_{d-1}\in V(\tGa')$. We now consider the free cover $p'':\tGa''\to \Ga''$ obtained from $p':\tGa'\to \Ga'$ in the following way. For each $i=0,\ldots,d-1$, we replace $v'_i$ with an undilated vertex $v''_i$ with an attached loop $e_i$, and replace $\tv'_i$ with a pair of vertices $\tv''^{\pm}_i$ connected by a pair of edges $\te''^{\pm}_i$. For each half-edge $h\in H(\Ga')$ rooted at $v'_i$, we attach its preimages $\thh^{\pm}$ to the vertices $\tv''^{\pm}_i$ in any manner. The result is a free cover $\pi'':\tGa''\to \Ga''$ whose contraction along the loops $e_0,\ldots,e_{d-1}$ is the edge-free cover $\pi':\tGa'\to \Ga'$. We denote $g(\Ga'')=g$, so that $g(\Ga)=g+m_{\dil}-n_{\dil}$.

We now pick a spanning tree $T\subset \Ga''$ and let $E(\Ga'')\backslash E(T)=\{e_0,\ldots,e_{g-1}\}$ be the complementary edges, where the first $d$ of the $e_i$ are the loops at the vertices $v_i$, as defined above. Let $\ep''_0,\ldots,\ep''_{g-1}$ and $\tep''_0,\tep''^{\pm}_1,\ldots,\tep''^{\pm}_{g-1}$ be the bases of $H_1(\Ga'',\ZZ)$ and $H_1(\tGa'',\ZZ)$ defined in Proposition~\ref{prop:freebasis}. The edges $e_0,\ldots,e_{d-1}$ are loops, so they form closed cycles and hence in fact $\ep_i=e_i$ for $i=0,\ldots,d-1$. Furthermore, the edges $\te''^+_i$ and $\te''^-_i$ have the same root vertices $\tv''^{\pm}_i$ for $i=0,\ldots,d-1$. This implies that $\tep''_0=\te''^+_0+\te''^-_0$, since the edge $\te''^+_0$ is contained in the spanning tree $\tT$. Also, for $i=1,\ldots,d-1$ the cycle $\tep''^-_i$ is obtained from $\tep''^+_i$ by replacing $\te''^+_i$ with $\te''^-_i$ and reversing the direction of the remaining path.

We now let $\ep'_0,\ldots,\ep'_{g-1}$ and $\tep'_0,\tep'^{\pm}_1,\ldots,\tep'^{\pm}_{g-1}$ denote the cycles in respectively $H_1(\Ga',\ZZ)$ and $H_1(\tGa',\ZZ)$ obtained by contracting the cycles $\ep''_0,\ldots,\ep''_{g-1}$ and $\tep''_0,\tep''^{\pm}_1,\ldots,\tep''^{\pm}_{g-1}$ defined above, in other words by setting $e''_i$ and $\te''^{\pm}_i$ to zero for $i=0,\ldots,d-1$. We see that $\tep'_0=0$ and $\ep'_i=0$ for $i=0,\ldots,d-1$. The remaining cycles $\ep'_d,\ldots,\ep'_{g-1}$ form a basis for $H_1(\Ga',\ZZ)$. Furthermore, we see that $\tep'^-_i=-\tep'^+_i$ for $i=1,\ldots,d-1$, and the cycles $\tep'^+_1,\ldots,\tep'^+_{d-1}$ and $\tep'^{\pm}_d,\ldots,\tep'^{\pm}_{g-1}$ form a basis for $H_1(\tGa',\ZZ)$. These bases satisfy the relations 
\begin{align*}
    \iota'_*(\tep'^+_i) &= -\tep'^+_i, 
    & \pi'_*(\tep'^+_i) &= 0,
    &&& i &= 1,\ldots,d-1, \\
    \iota'_*(\tep'^{\pm}_i) &= \tep'^{\mp}_i,
    & \pi'_*(\tep'^{\pm}_i) &= \ep_i,
    & \pi'^*(\ep'_i) &= \tep'^+_i+\tep'^-_i,
    & i &= d,\ldots,g-1.
\end{align*}

The edge set $E(\Ga')$ is identified with the set of non-dilated edges of $E(\Ga)$, hence we can view each cycle $\tep'^+_i$ as a simplicial chain $\tep^+_i$ in $\Ga$, with boundary $\partial(\tep^+_i)$ supported on the set of dilated vertices.  We claim that $\tep^+_i$ is in fact closed. Indeed, since $\iota_*(\tep^+_i)=-\tep^+_i$, it consists of a linear combination of expressions of the form $\te^+-\te^-$ for certain non-dilated pairs of edges, and if a root vertex of $\te^+$ is dilated, then it is also a root vertex of $\te^-$. Hence $\partial(\tep^+_i)=0$ and $\tep^+_i$ is a cycle, and we relabel $\tbe_i=\tep^+_i$ for $i=1,\ldots,B=d-1$.

Similarly, for each $i=d,\ldots,g-1$, let $\tep^{\pm}_i$ be the cycle $\tep'^{\pm}_i$, but viewed as a chain on $\Ga$. The boundaries $\partial(\tep^+_i)$ and $\partial(\tep^-_i)$ are equal and supported on the set of dilated vertices, so we can find a chain $\zeta_i$ supported on $\tGa_{\dil}$ such that $\tal^{\pm}_i=\tep^{\pm}_{i+d-1}+\zeta_{i+d-1}$ for $i=1,\ldots,A=g-d$ is a closed cycle on $\tGa$. Denoting $\al_i=\pi_*(\tal^{\pm}_i)$, we see that the $\tal^{\pm}_i$ and the $\al_i$ satisfy the required relations.

Finally, we let $\ga_1,\ldots,\ga_C$ be a basis for $H_1(\Ga_{\dil},\ZZ)$, and let $\tga_1,\ldots,\tga_C$ be the preimages of these cycles on $\tGa_{\dil}$. This completes the required basis. \qedhere
\end{proof}

The properties of the continuous and divisorial Prym varieties now follow directly.

\begin{proof} [Proof of Proposition~\ref{prop:Prymallpropetries}] We now again distinguish cycles and differential forms and denote the principal polarizations on $\Jac(\tGa)$ and $\Jac(\Ga)$ by $d:H_1(\tGa,\ZZ)\to \Omega^1_{\tGa}(\tGa)$ and $d:H_1(\Ga,\ZZ)\to \Omega^1_{\Ga}(\Ga)$, respectively, so that a $1$-form corresponding to a cycle $\ga$ is denoted $d\ga$. To describe the continuous and divisorial Pryms $\Prym_c(\tGa/\Ga)$ and $\Prym_d(\tGa/\Ga)$, we use Propositions~\ref{prop:freebasis} and~\ref{prop:dilatedbasis} to give explicit bases for their defining lattices (see~\eqref{eq:Prymd} and~\eqref{eq:Prymc}).

First, assume that $\pi$ is free, where we are reproving the results of~\cite{JensenLen}. Proposition~\ref{prop:freebasis} shows that the first lattices are equal and have basis
\[
(\Coker \pi^*)^{\torf}=\Omega^1_{\tGa}(\tGa)/\Ker (\Id-\iota^*)= \big\langle[d\tal^+_1],\ldots,[d\tal^+_{g-1}] \big\rangle,
\]
where $[d\tal]$ denotes the class of $d\tal\in \Omega^1_{\tGa}(\tGa)$. The second lattices are also equal:
\[
\quad \Ker \pi_*=\Im (\Id-\iota_*)= \big\langle\tal^+_1-\tal^-_1,\ldots,\tal^+_{g-1}-\tal^-_{g-1}\big\rangle.
\]
Hence $\Prym_c(\tGa/\Ga)=\Prym_d(\tGa/\Ga)$ and the dimension is $g_0=g(\Ga)-1$. Furthermore, $\Im(\Id-\iota_*)$ is equal to $\Ker \pi_*$ and is already saturated in $H_1(\tGa,\ZZ)$, hence $\Prym_d(\tGa/\Ga)$ is the group-theoretic image of $\Id-\iota$. By Proposition~\ref{prop:ker}, the number of connected components of the group-theoretic kernel of $\pi_*$ is the index of $\pi_* (H_1(\tGa,\ZZ) )$ in its saturation $\big(\pi_* (H_1(\tGa,\ZZ) ) \big)^{\sat}$ in $H_1(\Ga,\ZZ)$, and Proposition~\ref{prop:freebasis} shows that
\[
\big(\pi_* (H_1(\tGa,\ZZ) )\big)^{\sat} = H_1(\Ga,\ZZ),\qquad \big[H_1(\Ga,\ZZ):\pi_*(H_1(\tGa,\ZZ))\big]=2.
\]
Finally, the principal polarization $d:H_1(\tGa,\ZZ)\to \Omega^1_{\tGa}(\tGa)$ induces the polarization
\begin{align*}
    \zeta: \big\langle\tal^+_1-\tal^-_1,\ldots,\tal^+_{g-1}-\tal^-_{g-1} \big\rangle &\longrightarrow \big\langle[d\tal^+_1],\ldots,[d\tal^+_{g-1}] \big\rangle,\\ 
    \tal^+_i-\tal^-_i &\longmapsto 2[d\tal^+_i]
\end{align*}
on $\Prym_c(\tGa/\Ga)=\Prym_d(\tGa/\Ga)$, which is twice the principal polarization $\zeta_c:\tal^+_i-\tal^-_i\mapsto[d\tal^+_i]$.

We now consider a dilated double cover $\pi:\tGa\to \Ga$ with dilation index $d(\tGa/\Ga)$, so that $A=g_0-d(\tGa/\Ga)+1$ and $B=d(\tGa/\Ga)-1$. In terms of the bases given by Proposition~\ref{prop:dilatedbasis} we have
\begin{align*}
    \Ker (\Id-\iota^*) &= \big\langle d\tal^+_1+d\tal^-_1,\ldots,d\tal^+_A+d\tal^-_A,d\tga_1,\ldots,d\tga_C \big\rangle, \\
    \Im \pi^* &= \big\langle d\tal^+_1+d\tal^-_1,\ldots,d\tal^+_A+d\tal^-_A,2d\tga_1,\ldots,2d\tga_C \big\rangle,
\end{align*}
hence the first lattices are the same and have basis
\begin{equation}
    \label{eq:Prymcbasis1}
    (\Coker \pi^*)^{\mathrm{tf}}=\Omega^1_{\tGa}(\tGa)/\Ker (\Id-\iota^*)= \big\langle[d\tal^+_1],\ldots,[d\tal^+_A],[d\tbe_1],\ldots,[d\tbe_B] \big\rangle.
\end{equation}
The second lattices, however, are distinct. Indeed, the second lattice of $\Prym_c(\tGa/\Ga)$ 
\begin{equation}
\label{eq:Prymcbasis2}
\Im (\Id-\iota_*)= \big\langle \tal^+_1-\tal^-_1,\ldots,\tal^+_A-\tal^-_A,2\tbe_1,\ldots,2\tbe_B \big\rangle,
\end{equation}
while the second lattice of $\Prym_d(\tGa/\Ga)$ is
\[
\Ker \pi_*=(\Im (\Id-\iota_*))^{\mathrm{sat}}= \big\langle \tal^+_1-\tal^-_1,\ldots,\tal^+_A-\tal^-_A,\tbe_1,\ldots,\tbe_B \big\rangle.
\]
Hence $\Prym_d(\tGa/\Ga)$ is the image of $\Id-\iota$, and the natural map $\Prym_c(\tGa/\Ga)\to \Prym_d(\tGa/\Ga)$ is a free isogeny with geometric degree equal to the index of $\Im (\Id-\iota_*)$ in $\Ker \pi_*$, which is equal to $2^B=2^{d(\tGa/\Ga)-1}$.

In the dilated case the map $\pi_*$ is surjective, hence $\pi_*(H_1(\tGa,\ZZ))$ is saturated in its image and the group-theoretic kernel of $\pi$ has a single connected component. It remains to compute the induced polarizations. On $\Prym_d(\tGa/\Ga)$, the induced polarization $\Ker \pi_*\to (\Coker \pi^*)^{\mathrm{tf}}$ sends $\tal^+_i-\tal^-_i$ to $2[d\tal^+_i]$ and $\beta_j$ to $[d\beta_j]$, hence has type $(1^B,2^A)$. On the other hand, the induced polarization $\zeta:\Im(\Id-\iota_*)\to \Omega^1_{\tGa}(\tGa) / \Ker(\Id - \iota^\ast)$ on $\Prym_c(\tGa/\Ga)$ is equal to $\zeta=2\zeta_c$, where $\zeta_c$ is the principal polarization
\begin{equation}
    \begin{aligned}
        \zeta_c:\Im(\Id-\iota_*) &\longrightarrow \Omega^1_{\tGa}(\tGa) / \Ker(\Id - \iota^\ast), \\
        \tal^+_i-\tal^-_i &\longmapsto [d\tal^+_i], &i&=1,\ldots,A, \\
        \beta_j &\longmapsto [d\beta_j], &j&=1,\ldots,B.
    \end{aligned}
\label{eq:Prymcpp}
\end{equation}
This completes the proof. \qedhere
\end{proof}

We now define the tropical Abel--Prym map associated to a double cover $\pi:\tGa\to \Ga$. Choose a point $q\in \tGa$, and for any $p\in \tGa$ let $\ga_p$ be a path from $q$ to $p$. It is clear that the integral of a $1$-form $\omega\in \Omega^1_{\tGa}(\tGa)$ along the chain $\ga_p-\iota_*\ga_p$ depends only on the class of $\omega$ modulo $\Ker(\Id-\iota^*)$. Furthermore, this integral is well-defined modulo integration over elements of $\Im(\Id - \iota_*)$. Hence we make the following definition.
 
\begin{definition} Let $\pi:\tGa\to \Ga$ be a double cover and let $p\in \tGa$. The \emph{tropical Abel--Prym map} with base point $q$ is
    \begin{align*}
        \psi_q : \tilde \Gamma &\longrightarrow \Prym_c(\tGa/\Ga) \\
        p &\longmapsto \left( \omega \mapsto \int_{\gamma_p} \omega - \int_{\iota_\ast \gamma_p} \omega \right).
    \end{align*}
\end{definition}

In terms of the factorization~\eqref{eq:Prymmaps}, the Abel--Prym map is simply the composition $\psi_q=\epsilon\circ \phi_q$, where $\phi_q:\tGa\to \Jac(\tGa)$ is the Abel--Jacobi map, hence in particular it is a morphism of rational polyhedral spaces. The composition of $\psi_q$ with the isogeny $\gamma:\Prym_c(\tGa/\Ga)\to \Prym_d(\tGa/\Ga)$ is the divisorial Abel--Prym map $p\mapsto p-q-\iota(p-q)$. We note that if $\pi$ is a free double cover, then every divisor of the form $p-\iota(p)$ lies in the odd connected component of $\Ker \pi_*$, but the difference lies in the even connected component $\Prym_d(\tGa/\Ga)$.

The Abel--Prym map possesses a universal property analogous to Proposition~\ref{prop:universal_property_AbelJacobi_map}. For the algebraic version (including the case of a ramified double cover) see \cite{Masiewicki}.

\begin{proposition} \label{prop:universal_property_Prymc}
    Let $\pi:\tilde \Gamma\to \Ga$ be a double cover with associated involution $\iota:\tGa\to \tGa$. Let $\chi: \tilde\Gamma\to X$ be a morphism of rational polyhedral spaces to an integral torus $X=(\La,\La', [\cdot, \cdot])$. Assume that $\chi \circ \iota = - \chi$ and that the induced morphism on cotangent sheaves $\chi^{-1}\Omega^1_X \to \Omega^1_{\tilde\Gamma}$ takes values in $\Im \big(\Id - \iota^\ast : \Omega^1_{\tilde\Gamma} \to \Omega^1_{\tilde\Gamma} \big)$. Then there exists a unique homomorphism $\nu:\Prym_c(\tilde\Gamma/\Ga)\to X$ of integral tori such that the diagram
    \begin{center}
        \begin{tikzcd}
            \tilde\Gamma \arrow[r,"\chi"] \arrow[d,"\psi_q"] & X \arrow[d,"t_{-\chi(q)}"]\\
            \Prym_c(\tilde\Gamma/\Ga) \arrow[r,"\nu"] & X
        \end{tikzcd}
    \end{center}
    commutes for all $q\in \tilde \Gamma$.
\end{proposition}

The condition that $\chi^{-1}\Omega^1_X \to \Omega^1_{\tilde\Gamma}$ takes values in $\Im (\Id - \iota^\ast)$ is already implied by $\chi \circ \iota = - \chi$ in case the dilation index of the cover is equal to 1. Note further that the condition $\chi \circ \iota = -\chi$ implies that $\chi$ is constant on dilated edges of $\tilde\Gamma$. %\Dmitry{Do we need to prove this?}

\begin{proof} We first replace $\chi$ with $\chi_0 = t_{-\chi(q)} \circ \chi$ to avoid the translation.  As we noted above, the Abel--Prym map $\psi_q$ factors through the Abel--Jacobi map $\phi_q$ as $\psi_q=\epsilon \circ \phi_q$. By the universal property of the Jacobian, $\chi_0$ factors as 
    \begin{equation*}
        \begin{tikzcd}
            \tilde  \Gamma \arrow[rd, "\chi_0"] \arrow[d, "\phi_q"'] & \\
            \Jac(\tilde\Gamma) \arrow[r, "\mu"] & X
        \end{tikzcd}
    \end{equation*}
    where the map $\mu$ satisfies $\mu \circ \iota = -\mu$. Hence $\mu\circ (\Id-\iota)=2\mu$ and we may extend the commutative diagram as follows:
    \begin{equation*}
        \begin{tikzcd}
            \tilde  \Gamma \arrow[rd, "\chi_0"] \arrow[d, "\phi_q"'] & \\
            \Jac(\tilde\Gamma) \arrow[dd, bend right = 70, "\Id-\iota"'] \arrow[d,"\epsilon"] \arrow[r, "\mu"] & X \arrow[dd, "\times 2"] \\
            \Prym_c(\tilde\Gamma / \Gamma) \arrow[d] \arrow[ur,dashed,"\nu"] \arrow[dr, "\nu'"] &  \\
            \Jac(\tilde\Gamma) \arrow[r, "\mu"] & X.
        \end{tikzcd}
    \end{equation*}
    Here the left column is the factorization~\eqref{eq:Prymmaps} and the right vertical map on $X$ is multiplication by two.
    
    We claim that the composed map $\nu'$ is divisible by two, which implies the existence of the dashed map $\nu$ and proves the proposition. In terms of the basis~\eqref{eq:Prymcbasis2} (which subsumes the free and dilated case), the map $\nu'_\#:\Im(\Id-\iota_*)\to \Lambda'$ has the form
    \[
    \nu'_{\#}(\tal^+_i-\tal^-_i)=\mu_{\#}(\tal^+)-\mu_{\#}(\tal^-)=2\chi_*(\tal^+),\qquad \nu'_{\#}(2\tbe_j)=2\mu_{\#}(\tbe_j)=2\chi_*(\tbe_j),
    \]
    because $\mu\circ \iota=-\mu$. Hence there exists $\nu_\#:\Im(\Id-\iota_*)\to \Lambda'$ such that $\nu'_\#=2\nu_\#$. Similarly, the condition that $\chi^{-1}\Omega^1_X \to \Omega^1_{\tilde\Gamma}$ takes values in $\Im (\Id - \iota^\ast)$ implies that the image of $\mu^{\#}:\Lambda\to \Omega^1_{\tGa}(\tGa)$ lies in the submodule generated by the elements $d\tal^+_i-d\tal^-_i$ and $2d\tbe_j$. The classes of these elements are divisible by two in $\Omega^1_{\tGa}(\tGa)/\Ker (\Id-\iota^*)$, hence there exists a map $\nu^\#:\Lambda\to \Omega^1_{\tGa}(\tGa)/\Ker (\Id-\iota^*)$ such that $(\nu')^\#=2\nu^\#$. This completes the proof. \qedhere

\end{proof}

%\Felix{I believe this to be true if the dilation index is 1, but as I wrote in my email on April 29 it is not going to be true in general. In the paper by Masiewicki, Proposition 4.7 he presents a universal property for the Prym of a ramified double cover. He is asking for the additional condition on $\chi$ that the ramification points are identified under $\chi$. This does not work for us. }

%\begin{proof} \Felix{This is old and can go?} Let $\La_c=(\Coker \pi^*)^{tf}$ and $\La'_c=\Im(\iota_*-\Id)$. Let $[\cdot,\cdot]_c:\La_c\times \La'_c\to \RR$ be the pairing on $\Prym_c$. Let $X=(\La,\La',[\cdot,\cdot])$. We need maps:
%\[
%\nu^\#:\La\to \La_c,\nu_\#:\La'_c\to \La'.
%\]
%We have $\chi:\tGa\to X$ hence a map $\chi^*:\Omega_1(X)\to \Omega_1(\tGa)$, hence $\nu^\#$ is clearly the composition of $\chi^*$ with the quotient map $\Omega_1(\tGa)\to \La_c$. We also have $\chi_*:H_1(\tGa,\ZZ)\to H_1(X,\ZZ)=\La'$, hence $\nu_{\#}$ is the composition of $\chi_*$ with the inclusion $\La'_c\subset H_1(\tGa,\ZZ)$.

%What do we need to check:
%\[
%[\la,\nu_\#(\iota_*(\ga)-\ga)]=[\la,\chi_*(\iota_*(\ga))]-[\la,\chi_*(\ga)]=[\nu^\#(\la),\iota_*(\ga)-\ga]_c=\int_{\iota_*(\ga)}\nu^\#(\la)-\int_{\ga}\nu^\#(\la).
%\]
%for all $\la\in \La=\Omega^1_X(X)$ and all $\ga\in H^1(\tGa,\ZZ)$. This follows directly from (16), %nothing to check.

%\end{proof}

\begin{remark} \label{rem:AbelPrym_is_proper}
    The $d$-fold product of the tropical Abel-Prym map $\psi_q^d:\tGa^d\to \Prym_c(\tGa/\Ga)$, the free isogeny $\gamma : \Prym_c(\tGa / \Ga) \to \Prym_d(\tGa / \Ga)$, and the inclusion $i : \Prym_d(\tGa/\Ga) \to \Jac(\tilde \Gamma)$ are all proper morphisms of rational polyhedral spaces.
\end{remark}

\subsection{The tropical Poincar\'e--Prym formula}
Let $\Ga$ be a tropical curve of genus $g$, let $1\leq d\leq g$, let $\phi^d_q:\Ga^d\to \Jac(\Ga)$ be the $d$-fold product of the Abel--Jacobi map with an arbitrary base point $q$, and let $\tilde W_d$ denote the image of $\phi^d_q$. The tropical Poincar\'e formula \cite{GrossShokrieh_Poincare} states that 
\begin{equation}
    \cyc[\tilde W_d] = \frac{1}{(g-d)!} \xi^{g-d}\in H_{d,d}\big(\Jac(\Ga)\big).
\label{eq:tropicalPoincare}
\end{equation}
Here $\xi=\cyc[\Theta] \in H^{1,1}( \Jac(\Gamma))$ is the principal polarization on $\Jac(\Gamma)$ and we have identified $H_{d,d}(\Jac(\Ga))$ with $H^{g-d,g-d}(\Jac(\Ga))$ by Poincar\'e duality. The algebraic Poincar\'e formula has an analogue for Prym varieties, which is part of Welters' criterion and which we call the Poincar\'e--Prym formula. We conjecture that this formula holds in the tropical setting as well:

\begin{conjecture}(The tropical Poincar\'e--Prym formula) \label{conj:tropical_Welters_criterion}
    Let $\pi:\tilde \Gamma \to \Gamma$ be a double cover of tropical curves, let $q\in \tilde\Gamma$ be a base point, and let $g_0 =\dim \Prym_c(\tilde \Gamma / \Gamma)=g(\tGa)-g(\Ga)$. Then 
\[
(\psi^d_q)_*\cyc [\tGa^d]=\frac{2^d}{(g_0-d)!} \zeta^{g_0-d}\in H_{d,d}\big(\Prym_c(\tilde\Gamma / \Gamma)\big),
\]
for $1\leq d\leq g_0$, where $\zeta\in H^{1,1} \big(\Prym_c(\tGa/\Ga) \big)$ is the class of the principal polarization of $\Prym_c(\tilde\Gamma / \Gamma)$. 
\end{conjecture}

We only prove this result for $d=1$, which is all that we require to prove our main Theorems~\ref{thm:bigonal_construction} and~\ref{thm:tropical_Recillas_theorem}:

\begin{theorem} Conjecture~\ref{conj:tropical_Welters_criterion} holds for $d=1$:
\[
(\psi_q)_*\cyc [\tGa]=\frac{2}{(g_0-1)!} \zeta^{g_0-1}\in H_{1,1}\big(\Prym_c(\tilde\Gamma / \Gamma)\big).
\]
 \label{thm:tropical_Welters_criterion}   
\end{theorem}

\begin{proof} Since
\[(\psi_q)_*\cyc [\tGa]=\epsilon_* \big((\phi_q)_*\cyc[\tGa] \big) = \epsilon_* \big(\cyc [\widetilde{W}_1] \big),
\]
the tropical Poincar\'e formula for $d=1$ implies our result if we can show that 
\[
\epsilon_*\left(\frac{\xi^{\tg-1}}{(\tg-1)!}\right)=2\frac{\zeta^{g_0-1}}{(g_0-1)!}\in H_{1,1} \big(\Prym_c(\tGa/\Ga) \big),
\]
where $g_0=\tg-g$. We use the same arguments as in the proof of Proposition~\ref{prop:isomorphic_pptavs_via_homology}. First, assume that $\pi$ is free. In terms of the basis $H_1(\tGa,\ZZ)=\langle \tal^{\pm}_i,\tga_1\rangle$ given in Proposition~\ref{prop:freebasis} we see that the Poincar\'e dual of $\xi^{\tilde g-1}$ is given by
\[
\frac{\xi^{\tg-1}}{(\tg-1)!}=\sum_{i=1}^{g-1}\left[(d\tal^+_i)^*\otimes\tal^+_i+ (d\tal^-_i)^*\otimes\tal^-_i\right]+(d\tga_1)^*\otimes \tga_1.
\]
Similarly, the second lattice $\Im (\Id-\iota_*)$ of $\Prym_c(\tGa/\Ga)$ has basis $\langle \tal^+_i-\tal^-_i\rangle$, and in terms of the principal polarization~\eqref{eq:Prymcpp} we have
\[
\frac{\zeta^{g_0-1}}{(g_0-1)!}=\sum_{i=1}^{g-1}[d\tal^+_i]^*\otimes (\tal^+_i-\tal^-_i).
\]
The components $\epsilon_*=(\epsilon^\#)^*\otimes \epsilon_\#$ of the pushforward map 
\[
\epsilon_*: H_{1,1}\big(\Jac(\tGa)\big)=\left[\Omega^1_{\tGa}(\tGa)\right]^* \otimes H_1(\tGa,\ZZ)\longrightarrow H_{1,1} \big(\Prym_c(\tGa/\Ga) \big)=
\left[\Omega^1_{\tGa}(\tGa)/\Ker (\Id-\iota^*)\right]^*\otimes \Im (\Id-\iota_*)
\]
act on the basis elements as follows:
\[
(\epsilon^{\#})^*\left([d \tal^{\pm}_i]^*\right)=\pm [d\tal^+_i]^*,\qquad (\epsilon^{\#})^*\left(d\tga_1^*\right)=0,\qquad
\epsilon_\#(\tal^\pm_i)=\pm(\tal^+_i-\tal^-_i),\qquad \epsilon_\#(\tga_1)=0.
\]
Applying this to the above formula for $\xi^{\tg-1}/(\tg-1)!$, we obtain our result. The proof for the dilated case is similar, using instead the bases given in Proposition~\ref{prop:dilatedbasis}. \qedhere

\end{proof}

\section{Compatibility of the \texorpdfstring{$n$}{n}-gonal construction and tropical abelian varieties}
\label{sec:proof}
%%%%%%%%%%%%%%%%%%%%%%%%%%%%%%%%%%%%%%%%%%%%%%%%%%

We are now ready to prove the main theorems stated in the introduction. We restate the theorems for convenience, and begin with the trigonal construction. Throughout we assume that all tropical curves are smooth.

\begin{theorem} [Theorem~\ref{thm:tropical_Recillas_theorem}]

\label{thm:tropical_Recillas_theorem_restated}
    Let $K$ be a metric tree. The tropical trigonal and Recillas constructions establish a one-to-one correspondence
    \begin{equation*}
        \begin{tikzcd}[column sep = 4cm, every label/.append style = {font = \normalsize}]
        \left\{ 
        \begin{minipage}{0.3\textwidth}
            \begin{center}
                Tropical curves $\Pi$ with a harmonic map of degree 4 to $K$ with dilation profiles nowhere $(4)$ or $(2,2)$.
            \end{center}
        \end{minipage}
        \right\}
        \arrow[r, "\text{Recillas construction}", bend left = 20, start anchor = {[yshift = 1ex]east}, end anchor = {[yshift = 1ex]west}]
        &
        \left\{ 
        \begin{minipage}{0.3\textwidth}
            \begin{center}
                Free double covers $\tGa \to \Gamma$ with a harmonic map of degree 3 from $\Gamma$ to $K$.
            \end{center}    
        \end{minipage}
        \right\}
        \arrow[l, "\text{trigonal construction}", bend left = 20, start anchor = {[yshift = -1ex]west}, end anchor = {[yshift = -1ex]east}]
        \end{tikzcd}
    \end{equation*}
    and under this correspondence, the Prym variety of a double cover $\Prym(\tilde \Gamma / \Gamma)$ and the Jacobian $\Jac(\Pi)$ of the tetragonal curve are isomorphic as principally polarized tropical abelian varieties. 
\end{theorem}

We note that there is no difference between $\Prym_c(\tGa/\Ga)$ and $\Prym_d(\tGa/\Ga)$ because the double cover $\tGa\to \Ga$ is free, hence we simply refer to both objects as $\Prym(\tGa/\Ga)$. The techniques of tropical homology allow us to closely model the proof of the algebraic version of the theorem (see Theorem 12.7.2 in~\cite{BirkenhakeLange}). 

\begin{proof} Here we present the outline of the proof, and postpone  the necessary calculations and checks to a series of lemmas that are given later in this chapter. Recall that we have already established the bijection in Proposition~\ref{prop:trigonal_properties}, and it only remains to show that $\Prym(\tGa/\Ga)\cong \Jac(\Pi)$.

    We recall the setup and notation. Let $k : \Pi \to K$ be a generic tetragonal curve, so that $K$ does not have any points over which the degree profile of $k$ is $(2,2)$ or $(4)$. By the tropical Recillas construction (Definition~\ref{def:Recillas_construction}) we obtain a tower $\tilde \Gamma \overset{\pi}{\longrightarrow} \Gamma \overset{f}{\longrightarrow} K$, where $f:\Ga\to K$ is a trigonal curve and $\pi:\tGa\to \Ga$ is a free double cover. We denote $\iota:\tGa\to \tGa$ the associated involution. We choose graph models for our tropical curves, and by abuse of notation refer to edges and vertices of $\tGa$, $\Ga$, $\Pi$, and $K$.
    
\medskip

First, we \textbf{define a map} of rational polyhedral spaces $\chi : \tilde \Gamma \to \Jac (\Pi)$ such that $\chi\iota = -\chi$. To this end, we choose an Abel--Jacobi map $\Pi^2 \to \Jac(\Pi)$ suited to our purposes. Fix a point $x \in K$ and let
\begin{equation} \label{eq:bigonal_D}
    D = \sum_{y \in k^{-1}(x)} d_k(y)\cdot y \in \Div^+_4(\Pi)
\end{equation}
be the fiber of $k$ above $x$. The group $\Pic_0(\Pi)\cong \Jac(\Pi)$ is divisible because it is a real torus, hence we can find $M\in \Div_1(\Pi)$ such that $4M\sim D$ in $\Pic_4(\Pi)$. Let $\phi_M : \Pi \to \Jac(\Pi)$ be the Abel--Jacobi map associated to $M$ (note that $M$ may fail to be effective, in which case $\phi_M$ is actually a translation of an Abel--Jacobi map by a fixed divisor class). Moreover, define $L \in \Div_2(\Pi)$ as $L = 2M$ and denote
\[
\phi_L = \phi_M + \phi_M : \Pi^2 \longrightarrow \Jac(\Pi).
\]
The map $\phi_L$ is symmetric with respect to swapping the coordinates of $\Pi^2$, hence it descends to a map $\Div_2^+(\Pi) \to \Jac(\Pi)$. Define $\chi$ to be the composition of the inclusion $\tilde \Gamma \subseteq \Div^+_2(\Pi)$ from Equation \eqref{eq:tildeGamma_in_Sym2} with this descent of $\phi_L$. Since we are not using the polyhedral structure on $\Div_2^+(\Pi)$, we need to verify by hand that $\chi$ is a map of rational polyhedral spaces, and we do this in Lemma~\ref{lem:alpha} (where we also check that $\chi$ has the desired property $\chi\iota=-\chi$).

    We may now apply the universal property of the Prym variety from Proposition~\ref{prop:universal_property_Prymc} for any base point $q \in \tGa$ to obtain a commutative square
    \begin{equation*}
        \begin{tikzcd}
            \tilde \Gamma \arrow[r, "\chi"] \arrow[d, "\psi_q"] & \Jac(\Pi) \arrow[d, "t_{-\chi(q)}"] \\
            \Prym(\tilde \Gamma / \Gamma) \arrow[r, "\mu"] & \Jac(\Pi) .
        \end{tikzcd}
    \end{equation*}
    We note that it is not necessary to check the pullback condition on the cotangent sheaves, since the double cover $\pi:\tGa\to \Ga$ is free.
    Our goal is to show that $\mu$ is an isomorphism of principally polarized tropical abelian varieties. We saw already in Proposition~\ref{prop:trigonal_properties} that $\dim \Jac (\Pi) = g(\Pi) = g(\tilde \Gamma) -1 = \dim \Prym(\tilde \Gamma / \Gamma)$.
    By Theorem \ref{thm:tropical_Welters_criterion} we know that the pushforward of the fundamental class of $\tilde \Gamma$ along $\psi_q$ is equal to
    $$(\psi_q)_* \cyc[\tilde \Gamma] = \frac{2}{(g_0-1)!} \zeta^{g_0-1}\in H_{1,1}\big(\Prym(\tGa/\Ga)\big) ,$$
    where $g_0 = \dim \Prym(\tilde \Gamma / \Gamma)$ and $\zeta\in H^{1,1}\big(\Prym(\tGa/\Ga)\big)$ is the principal polarization.
    If we can show that 
\begin{equation}
\label{eq:mainformulatoprove}
\chi_* \cyc[\tilde \Gamma] = \frac{2}{(g_0-1)!} \xi^{g_0-1}\in H_{1,1}\big(\Jac(\Pi)\big),
\end{equation}
    where $\xi\in H^{1,1}\big(\Jac(\Pi)\big)$ is the principal polarization on $\Jac(\Pi)$, then $\mu_*(\zeta^{g_0-1})=\xi^{g_0-1}$ and therefore $\Prym(\tGa/\Ga)$ is isomorphic to $\Jac(\Pi)$ by the homological criterion of  Proposition~\ref{prop:isomorphic_pptavs_via_homology}. 
    \medskip
    
    We now \textbf{define a tropical 1-cycle} $A \in Z_1(\Pi^2)$ as follows. Recall that we can view $\tGa$ as a subset of $\Div_2^+(\Pi)$ (see Equation~\eqref{eq:tildeGamma_in_Sym2}), however, this does not induce the correct edge lengths on $\tGa$. Instead, we manually construct a cycle $A$ that represents the lift of $\tGa$ to $\Pi^2$ via the natural projection map $\Pi^2\to \Div_2^+(\Pi)$.     
    A tropical $1$-cycle on $\Pi^2$ is a map $\Pi^2\to \ZZ$, and for $x, y \in \Pi$ we set 
    \begin{equation}\label{eq:definition_cycle_A}
        A'(x, y) = \begin{cases}
            0 & \text{if }k(x)\neq k(y), \\
            1 & \text{if }k(x)=k(y)\text{ and }x \neq y, \\
            d - 1 & \text{if $x = y$, }
        \end{cases}
    \end{equation}
    where $d=d_k(x)=d_k(y)$ is the dilation factor of $k$ at $x=y$. 
    The support $|A'|$ is a purely 1-dimensional rational polyhedral space, contained in the preimage of the diagonal $\Delta_K$. 
    %The proof that $A$ is indeed a tropical 1-cycle is given in Lemma \ref{lem:A_cycle} below.
    One can show that $A'$ itself is a tropical 1-cycle. We will not carry this out and instead show in Lemma~\ref{lem:key_formula} below the slightly weaker statement that there exists a tropical 1-cycle $A$ which agrees with the map $A'$ away from a zero-dimensional locus. 
  \medskip
    
    \textbf{The key calculation} is the following formula
\begin{equation}
        \label{eq:key_equation}
[\De_{\Pi}]+A=4[\Pi\times p']+4[p'\times \Pi]\in H_{1,1}(\Pi^2),
\end{equation}
where $\De_{\Pi}$ is the diagonal, $p'\in \Pi$ is an arbitrary point, $[\cdot]$ denotes the fundamental class, and we identify cycles in $Z_1(\Pi^2)$ with their classes in $H_{1,1}(\Pi^2)$ via the cycle class map.

We prove this formula in several steps. First, we note that the spaces $K^2$ and $\Pi^2$ are smooth, so by Poincar\'e duality we can identify $H^{1,1}(K^2)\cong H_{1,1}(K^2)$ and $H^{1,1}(\Pi^2)\cong H_{1,1}(\Pi^2)$. Under this identification, pullback on cohomology induces a map $(k\times k)^*:H_{1,1}(K^2)\to H_{1,1}(\Pi^2)$, and applying the cycle class map to the main statement of Lemma~\ref{lem:key_formula} gives
    \begin{equation}
        \label{eq:key_equation2}
         [\Delta_\Pi] + A = (k \times k)^*[\Delta_K] \in H_{1,1}(\Pi^2) .
    \end{equation}
We then use Lemma~\ref{prop:diagonal_formula} to rewrite the right hand side of Equation~\eqref{eq:key_equation2} as
\begin{equation}
    \label{eq:key_equation3}
(k\times k)^*[\Delta_K] = (k\times k)^* \big( [K \times p] + [p \times K] \big)\in H_{1,1}(\Pi^2),
\end{equation}
where $p \in K$ is an arbitrary point. Finally, it is easy to show that $k^*[K]=[\Pi]$, which implies that the right hand sides of Equations~\eqref{eq:key_equation} and~\eqref{eq:key_equation3} are equal. %\Dmitry{We do need to say something more here}

To complete the proof, we compute the pushforward of Equation~\eqref{eq:key_equation} to $H_{1,1}(\Jac(\Pi))$ along the map $\phi_L:\Pi^2\to \Jac(\Pi)$.
    
    \begin{enumerate}
        \item \textbf{Claim:} $(\phi_L)_*[\Delta_\Pi] = 4(\phi_M)_* [\Pi]$. Indeed, by definition $\phi_L(x,y) = \phi_M(x) + \phi_M(y)$, so the restriction of $\phi_L$ to the diagonal $\Delta_\Pi$ is $\phi_M$ applied to the first coordinate followed by multiplication by $2$. But multiplication by $2$ on an integral torus induces multiplication by $4$ on $H_{1,1}$, so the claim follows.
        
        \item \textbf{Claim:} $(\phi_L)_*A = 2 \chi_* [\tilde\Gamma]$.
%        \Felix{This part is now full of technical subtleties but hopefully addresses referee (54) satisfactorily. } \Dmitry{I need to take a look at this}
        We define a continuous map $\tau : |A| \to \tilde\Gamma$ which factors the restriction of $\phi_L$ to $|A|$ as $\chi \circ \tau$ as follows. A point $(x,y)$ in $|A|$ is really a point on the diagonal of $\big[0, \ell_\Pi(e_1) \big] \times \big[0, \ell_\Pi(e_2) \big]$. Such a point is mapped under $\tau$ to the point with coordinates $d_\tau(e_1 \times e_2)\min\{x, y\}$ in the geometric realization $\big[0, \ell_{\tilde\Gamma}(e_1 + e_2)\big]$ of the edge $e_1 + e_2 \subseteq \tilde\Gamma$. Here we use the dilation factor
        \begin{equation} \label{eq:definition_tau}
            d_\tau(e_1 \times e_2) = \begin{cases}
                2 & \text{if $e_1 = e_2$ and } d_k(e_1) = d_k(e_2) = 2 \\
                1 & \text{else.} 
            \end{cases}
        \end{equation}
        We stress that $\tau$ is not a harmonic map of tropical curves and therefore not a morphism of (1-dimensional) rational polyhedral spaces! However, the problem only lies at vertices and restricting $\tau$ to the interior of any of the 1-dimensional faces of $|A|$ does give a morphism of rational polyhedral spaces onto an open edge of $\tilde \Gamma$. In particular, computing the lattice index in the definition of pushforward in Equation~\eqref{eq:definition_push_forward} can be done in two steps: Denote the image of $|A|$ in $\Jac(\Pi)$ under $\phi_L$ by $B$ and let
        \begin{itemize}
            \item $x \in |A|^\reg$ lie in the interior of the diagonal of $e_1 \times e_2$,
            \item $z = \tau(x) \in \tilde\Gamma^\reg$ in the interior of the edge $e_1 + e_2 \subseteq \tilde\Gamma$ and 
            \item $y = \chi(z) \in B^\reg$.
        \end{itemize}
        Then we may compute
        \begin{align*}
            \Big[ T_y^\ZZ B : d_x\phi_L \big(T_x^\ZZ |A| \big) \Big] 
            ={}&
            \Big[ T_z^\ZZ\tilde\Gamma : d_x\tau \big(T_x^\ZZ |A| \big) \Big] 
            \  \Big[ T_y^\ZZ B : d_z\chi \big(T_z^\ZZ \tilde\Gamma \big) \Big] \\
            ={}&
            d_\tau(e_1 \times e_2) \  \chi_\ast [\tilde\Gamma](y).
        \end{align*}
        Consequently, the value of $(\phi_L)_\ast A$ at $y$ is given as
        \begin{align*}
            (\phi_L)_\ast A(y) ={}& \sum_{x \in (\phi_L)^{-1}(y) \cap |A|^\reg} \Big[ T_y^\ZZ B : d_x\phi_L \big(T_x^\ZZ |A| \big) \Big] A(x)  \\
            ={}& \chi_\ast [\tilde\Gamma](y) 
            \underbrace{\sum_{x \in (\phi_L)^{-1}(y) \cap |A|^\reg} d_\tau(e_1 \times e_2)   A'(x)}_{(\ast)}
        \end{align*}
        and comparing the definitions of $d_\tau$ in Equation~\eqref{eq:definition_tau} and $A'$ in Equation~\eqref{eq:definition_cycle_A} shows that $(\ast) = 2$. Hence the tropical 1-cycles $(\phi_L)_\ast A$ and $2\chi_\ast [\tilde\Gamma]$ agree away from some 0-dimensional locus and thus define the same element of $Z_1(\Jac(\Pi))$.
        
        \item \textbf{Claim:} $(\phi_L)_*[\Pi \times p] = (\phi_L)_*[p \times \Pi] = (\phi_M)_* [\Pi]$. As above,  $\phi_L(x,p) = \phi_M(x) + \phi_M(p)$ for all $x \in \Pi$. But this means that the restriction of $\phi_L$ to $\Pi \times p$ is the composition of $\phi_M$ with a translation. Since translations induce the identity on homology, the claim follows.
    \end{enumerate}
    
    Summing up we have just shown $4(\phi_M)_*[\Pi] + 2\chi_*[\tilde\Gamma] = 8(\phi_M)_*[\Pi]$. 
    Solving for $\chi_*[\tilde \Gamma]$ and plugging in the Poincar\'e formula~\eqref{eq:tropicalPoincare}
    \[ (\phi_M)_*\cyc[\Pi] = \frac{1}{(g_0-1)!} \xi^{g_0-1} \]
    we obtain Equation~\eqref{eq:mainformulatoprove}, and the proof of Theorem~\ref{thm:tropical_Recillas_theorem} is complete. \qedhere
\end{proof}

\begin{example}
    \label{ex:trigonal_construction_computation} We consider the $(3,2)$-graph tower and the tetragonal graph shown on Figure~\ref{fig:example_trigonal}. We promote these to a tower $\tGa\to \Ga\to K$ of tropical curves and a generic tetragonal curve $\Pi\to K$ by assigning edge lengths to $K$, which we denote, going left to right, by $a$, $b$, $c$, $d$, and $e$.
    
    Let us compute the Prym and Jacobian varieties in this example explicitly to see that they are the same. We begin by computing $\Prym(\tilde\Gamma / \Gamma)$. It is possible to construct a basis for $\Ker \pi_*$ using Proposition~\ref{prop:freebasis}, but in this case it is easier to choose a basis by hand. It is clear that $\Ker \pi_*$ is spanned by the elements $\teta_1^+-\teta_1^-$ and $\teta_2^+-\teta_2^-$, where $\teta^{\pm}_i$ are the elements of $H_1(\tGa,\ZZ)$ shown on Figure~\ref{fig:example_trigonal_basis}. Furthermore, the module $(\Coker \pi^*)^{tf}$ (viewed as a quotient of $H_1(\tGa,\ZZ)$, which has been canonically identified with $\Omega^1_{\tGa}(\tGa)$) is spanned by the classes $[\teta_1^+]$ and $[\teta_2^+]$. Keeping in mind that the dilation factors of $\Gamma \to K$ also affect the lengths of edges in $\tilde\Gamma$, we obtain that the pairing of $\Prym(\tGa/\Ga)$ is
    \medskip
    
\begin{center}
    \begin{tabular}{c|c c}
        $[\cdot, \cdot]$ & $[\teta_1^+]$ & $[\teta_2^+]$ \\
        \hline
        $\teta_1^+ - \teta_1^-$ & $2(b+c+d)$ & $b+c+d$ \\
        $\teta_2^+ - \teta_2^-$ & $b+c+d$ & $\frac{3}{2}a + 2b + \frac{3}{2}c + 2d + \frac{3}{2} e$ .
    \end{tabular}
\end{center}    
\medskip
To compute the intersection matrix for $\Jac(\Pi)$, we choose the basis $\epsilon, \delta$ for $H_1(\Pi, \ZZ)$ depicted in Figure \ref{fig:example_tetragonal_basis}. The edge length pairing $[\cdot, \cdot] : H_1(\Pi, \ZZ) \times H_1(\Pi, \ZZ) \to \RR$ (where we have identified $\Omega^1_{\Pi}(\Pi)=H_1(\Pi,\ZZ)$ via the principal polarization) evaluated on the basis yields \medskip

\begin{center}
    \begin{tabular}{c|c c}
        $[\cdot, \cdot]$ & $\epsilon$ & $\delta$ \\
        \hline
        $\epsilon$ & $2(b+c+d)$ & $b+c+d$ \\
        $\delta$ & $b+c+d$ & $\frac{3}{2}a + 2b + \frac{3}{2}c + 2d + \frac{3}{2} e$ ,
    \end{tabular}
\end{center}
\medskip
which is evidently the same table as before. Hence 
\begin{align*}
    \mu^\# : \epsilon &\longmapsto [\tilde\eta_1^+] & \mu_\# :  \tilde\eta_1^+ - \tilde\eta_1^- &\longmapsto \epsilon \\
    \delta &\longmapsto [\tilde\eta_2^+] & \tilde\eta_2^+ - \tilde\eta_2^- &\longmapsto \delta
\end{align*}
defines an isomorphism $\Prym(\tilde\Gamma / \Gamma) \to \Jac(\Pi)$ of integral tori. To see that this is in fact an isomorphism of principally polarized tropical abelian varieties we quickly verify the compatibility with the principal polarizations which are given by $\epsilon \mapsto \epsilon$ and $\delta \mapsto\delta$ on $\Jac(\Pi)$ and $\tilde\eta_i^+ - \tilde\eta_i^- \mapsto [\tilde\eta_i^+]$ on $\Prym(\tGa / \Ga)$.

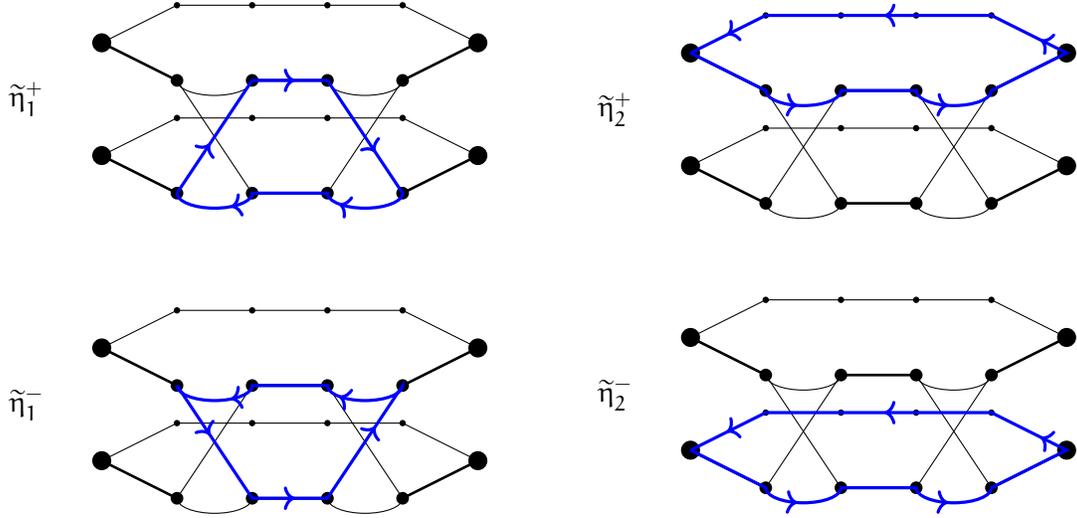
\begin{figure}[htb]
    \centering
    \begin{subfigure}{0.4\textwidth}
        \centering
        \begin{tikzpicture}
        	%Generator 1^+
        	\draw (-1, 1.25) node {$\teta_1^+$};
        	\drawGammaTilde
        	\drawCycle{
        		\draw[postaction = decorate] (1,0) -- ++ (1, 1.5) -- ++ (1,0) -- ++ (1, -1.5) arc (360:180:0.5 and 0.2) -- ++ (-1, 0) arc (360:180:0.5 and 0.2);
        	}
        \end{tikzpicture}
    \end{subfigure}
    \hspace{1cm}
    \begin{subfigure}{0.4\textwidth}
        \centering
        \begin{tikzpicture}
        	%Generator 2^+
        	\draw (-1, 1.25) node {$\teta_2^+$};
        	\drawGammaTilde
        	\drawCycle[0.15]{
        		\draw[postaction = decorate] (0, 2) -- ++ (1, -0.5) arc (180:360: 0.5 and 0.2) -- ++ (1,0) arc (180:360:0.5 and 0.2) -- ++ (1, 0.5) -- ++ (-1, 0.5) -- ++ (-3,0) -- ++ (-1, -0.5);
        	}
        \end{tikzpicture}
    \end{subfigure}
    \\[1cm]
    \begin{subfigure}{0.4\textwidth}
        \centering
        \begin{tikzpicture}
            %Generator 1^-
            \draw (-1, 1.25) node {$\teta_1^-$};
        	\drawGammaTilde
        	\drawCycle{
        		\draw[postaction = decorate] (1,1.5) -- ++ (1, -1.5) -- ++ (1,0) -- ++ (1, 1.5) arc (360:180:0.5 and 0.2) -- ++ (-1, 0) arc (360:180:0.5 and 0.2);	
        	}
        \end{tikzpicture}
    \end{subfigure}
    \hspace{1cm}
    \begin{subfigure}{0.4\textwidth}
        \centering
        \begin{tikzpicture}
        	%Generator 2^-
        	\draw (-1, 1.25) node {$\teta_2^-$};
        	\drawGammaTilde
        	\drawCycle[0.15]{
        		\draw[postaction = decorate] (0, 0.5) -- ++ (1, -0.5) arc (180:360: 0.5 and 0.2) -- ++ (1,0) arc (180:360:0.5 and 0.2) -- ++ (1, 0.5) -- ++ (-1, 0.5) -- ++ (-3,0) -- ++ (-1, -0.5);
        	}
        \end{tikzpicture}
    \end{subfigure}
    \caption{The four cycles in $H_1(\tilde\Gamma, \ZZ)$ that we use for our computation in Example~\ref{ex:trigonal_construction_computation}.}
    \label{fig:example_trigonal_basis}
\end{figure}

\begin{figure}[htb]
    \centering
    \begin{subfigure}{0.4\textwidth}
        \centering
        \begin{tikzpicture}
            \draw (-1, 0) node {$\epsilon$};
            \drawPi
            \drawCycle[0.15]{
                \draw[postaction = decorate] (1, -0.5) -- ++ (1, -0.5) -- ++ (1,0) -- ++ (1, 0.5) -- ++ (-1, 0.5) -- ++ (-1, 0) -- ++ (-1, -0.5);
            }
        \end{tikzpicture}
    \end{subfigure}
    \hspace{1cm}
    \begin{subfigure}{0.4\textwidth}
        \centering
        \begin{tikzpicture}
            \draw (-1, 0) node {$\delta$};
            \drawPi
            \drawCycle[0.05]{
                \draw[postaction = decorate] (0,0) -- ++ (2, -1) -- ++ (1, 0) -- ++ (2, 1) -- ++ (-2, 1) -- ++ (-1, 0) -- ++ (-2, -1);
            }
        \end{tikzpicture}
    \end{subfigure}
    \caption{Basis for $H_1(\Pi, \ZZ)$.}
    \label{fig:example_tetragonal_basis}
\end{figure}
\end{example}

We now present a series of lemmas which fill in the missing details in the proof of Theorem~\ref{thm:tropical_Recillas_theorem} above.

\begin{lemma}
    \label{lem:alpha}
    The map $\chi : \tilde\Gamma \to \Jac(\Pi)$ defined in the proof of Theorem \ref{thm:tropical_Recillas_theorem_restated} is a map of rational polyhedral spaces. Furthermore, it satisfies $\chi\iota = - \chi$.
\end{lemma}

\begin{proof}
    We start by verifying that $\chi$ is a map of rational polyhedral spaces. We have to check that the pullback of a 1-form $\omega \in \Omega^1_{\Jac(\Pi)} \cong \Omega^1_{\Pi}$ along $\chi$ is a 1-form on $\tilde \Gamma$. For this we need check that $\chi^*\omega$ has integer slopes on every edge of $\tilde \Gamma$, and that the sum of slopes around every vertex of $\tilde \Gamma$ is $0$.
    
    Recall that an edge $\te \in E(\tilde\Gamma)$ corresponds to a pair of edges of $f_1, f_2 \in E(\Pi)$ (which may be the same), with all three mapping to the same edge $e\in E(K)$ (see table in Definition~\ref{def:Recillas_construction}). The slope of $\chi^*\omega$ on $\te$ is equal to
    \begin{equation*} \label{eq:slope}
        \slope_{\te}(\chi^*\omega) = \frac{1}{\ell(\te)} \bigg( s_1 \ell(f_1) + s_2 \ell(f_2) \bigg) ,
    \end{equation*}
    where $s_1=\slope_{f_1}(\omega)$ and $s_2=\slope_{f_2}(\omega)$.
    If we denote the composed map $\tilde f = f \circ \pi$, then by harmonicity, we have
    \[
    \frac{\ell(f_i)}{\ell(\te)}=\frac{\ell(e)/d_k(f_i)}{\ell(e)/d_{\tf}(\te)}=\frac{d_{\tf}(\te)}{d_k(f_i)}
    \]
    for $i=1,2$. We can now verify case-by-case that $\slope_{\te}(\chi^*\omega)$ is an integer
    \medskip
    
    \begin{center}
    \begin{tabular}{c|c|c|c}
        $d_{\tf}(\te)$ & $d_k(f_1)$ & $d_k(f_2)$ & $\slope_{\te}(\chi^*\omega)$ \\
        \hline
        1 & 1 & 1 & $s_1 + s_2$ \\
        1 & 2 & 2 & $s_1 = s_2$ \\
        2 & 1 & 2 & $2s_1 + s_2$ \\
        3 & 1 & 3 & $3s_1 + s_2$ \\
        3 & 3 & 3 & $2s_1 = 2s_2$
    \end{tabular}
    \end{center}
    \medskip
    where in the second and fifth row the edges $f_1=f_2$ are the same and hence $s_1=s_2$.
    
    Now let $v\in V(\tGa)$ be a vertex of $\tilde\Gamma$, corresponding to two (not necessarily distinct) vertices $w_1,w_2\in V(\Pi)$. We want to show that the sum of outgoing slopes of $\chi^*\omega$ over the half-edges rooted at $v$ is 0, i.e. that
    \[ \sum_{\th\in T_v\tGa} \slope_{\th}(\chi^*\omega) 
        = \sum_{h\in T_{\tf(v)}K} \Bigg(
        \sum_{\th \in \tf^{-1}(h)\cap T_v\tGa} 
        \slope_{\th}(\chi^*\omega) \Bigg) = 0 ,
    \]
    where we have split the sum according to the image half-edge $h=\tf(\th)$. We do this by showing that for every vertex $v$ there are integers $a, b \in \ZZ$ such that for every $h \in T_{\tf(v)}K$ we have
    \begin{equation}
        \label{eq:harmonicity}
        \sum_{\th \in \tf^{-1}(h)\cap T_v\tGa} \slope_{\th}(\chi^*\omega) 
        = a \Bigg( \sum_{\th_1 \in k^{-1}(h)\cap T_{w_1}\Pi}
         \slope_{\th_1}(\omega) \Bigg) 
        + b \Bigg(\sum_{\th_2 \in k^{-1}(h)\cap T_{w_2}\Pi} \slope_{\th_2}(\omega) \Bigg) .
    \end{equation}
Summing Equation \eqref{eq:harmonicity} over all $h \in T_{\tf(v)}K$, the right hand side is 0 by harmonicity of $\omega$. The numbers $a$ and $b$ are determined by case distinction and direct computation in Figure \ref{tab:harmonicity}.
    
    \begin{figure}[thb]
        \centering
        {\setlength{\extrarowheight}{0.6cm}
\begin{tabular}{l|m{4cm}|m{3.5cm}|m{3cm}|m{2cm}}
    $d_{\tf}(v)$ 
    & \centering $w_1$ and $w_2$ 
    & Local picture of $\tilde\Gamma$ at $v$ over $h$ together with outgoing slopes of $\chi^*\omega$.
    & Local picture of $\Pi$ at $w_1$ and $w_2$ over $h$, and outgoing slopes of $\omega$
    & $(a, b)$
    \\
    %------------------
    %--- d_v = 1 ------
    %------------------
    \hline\hline
    \multirow{3}{1cm}{1}
    & $d_{k}(w_1) = d_k(w_2) = 1$ and $w_1 \neq w_2$
    & \begin{tikzpicture}
        \halfedge{0}
        \draw (1.2, 0) node[anchor = west] {$s_1 + s_2$};
    \end{tikzpicture}
    & \begin{tikzpicture}
        \halfedge{0}
        \draw (1.2, 0.0) node[anchor = west] {$s_2$};%
        \halfedge{0.4}
        \draw (1.2, 0.4) node[anchor = west] {$s_1$};%
    \end{tikzpicture}
    & $a = b = 1$ \\
    \cline{2-5}
    & \multirow{2}{4cm}{$d_k(w_1) = 2$ and $w_1 = w_2$}
    & \begin{tikzpicture}
        \halfedge{0}
        \draw (1.2, 0) node[anchor = west] {$s_1 + s_2$};%
    \end{tikzpicture}
    & \begin{tikzpicture}
        \twohalfedges{0}
        \draw (1.2, 0.15) node[anchor = west] {$s_1$};%
        \draw (1.2, -0.15) node[anchor = west] {$s_2$};%
    \end{tikzpicture}
    & \multirow{2}{3cm}{$a = 1$} \\
    & & \begin{tikzpicture}
        \halfedge{0}
        \draw (1.2, 0) node[anchor = west] {$s$};%
    \end{tikzpicture}
    & \begin{tikzpicture}
        \halfedge[2]{0}
        \draw (1.2, 0) node[anchor = west] {$s$};%
    \end{tikzpicture}
    & \\
    %------------------
    %--- d_v = 2 ------
    %------------------
    \hline
    \multirow{2}{1cm}{2} 
    & \multirow{2}{4cm}{$d_k(w_1) = 1$, $d_k(w_2) = 2$}
    & \begin{tikzpicture}
        \twohalfedges{0}
        \draw (1.2, 0.15) node[anchor = west] {$s_1 + s_2$};%
        \draw (1.2, -0.15) node[anchor = west] {$s_1 + s_3$};%
    \end{tikzpicture}
    & \begin{tikzpicture}
        \twohalfedges{0}
        \draw (1.2, 0.15) node[anchor = west] {$s_2$};%
        \draw (1.2, -0.15) node[anchor = west] {$s_3$};%
        \halfedge{0.4}
        \draw (1.2, 0.4) node[anchor = west] {$s_1$};%
    \end{tikzpicture}
    & \multirow{2}{3cm}{$a = 2, b = 1$} \\
    & & \begin{tikzpicture}
        \halfedge[2]{0}
        \draw (1.2, 0) node[anchor = west] {$2s_1 + s_2$};%
    \end{tikzpicture}
    & \begin{tikzpicture}
        \halfedge[2]{0}
        \draw (1.2, 0) node[anchor = west] {$s_2$};%
        \halfedge{0.4}
        \draw (1.2, 0.4) node[anchor = west] {$s_1$};%
    \end{tikzpicture}
    & \\
    %------------------
    %--- d_v = 3 ------
    %------------------
    \hline
    \multirow{6}{1cm}{3} 
    & \multirow{3}{4cm}{$d_k(w_1) = 1$, $d_k(w_2) = 3$}
    & \begin{tikzpicture}
        \halfedge[3]{0}
        \draw (1.2, 0) node[anchor = west] {$3s_1 + s_2$};%
    \end{tikzpicture}
    & \begin{tikzpicture}
        \halfedge[3]{0}
        \draw (1.2, 0) node[anchor = west] {$s_2$};%
        \halfedge{0.4}
        \draw (1.2, 0.4) node[anchor = west] {$s_1$};%
    \end{tikzpicture}
    & \multirow{3}{3cm}{$a = 3, b = 1$} \\
    & & \begin{tikzpicture}
        \twohalfedges[2]{0}
        \draw (1.2, 0.15) node[anchor = west] {$2s_1 + s_2$};%
        \draw (1.2, -0.15) node[anchor = west] {$s_1 + s_3$};%
    \end{tikzpicture}
    & \begin{tikzpicture}
        \twohalfedges[2]{0}
        \draw (1.2, 0.15) node[anchor = west] {$s_2$};%
        \draw (1.2, -0.15) node[anchor = west] {$s_3$};%
        \halfedge{0.5}
        \draw (1.2, 0.5) node[anchor = west] {$s_1$};%
    \end{tikzpicture}
    & \\
    & & \begin{tikzpicture}
        \threehalfedges{0}
        \draw (1.2, 0.25) node[anchor = west] {$s_1 + s_2$};%
        \draw (1.2, 0) node[anchor = west] {$s_1 + s_3$};%
        \draw (1.2, -0.25) node[anchor = west] {$s_1 + s_4$};%
    \end{tikzpicture}
    & \begin{tikzpicture}
        \threehalfedges{0}
        \draw (1.2, 0.25) node[anchor = west] {$s_2$};%
        \draw (1.2, 0) node[anchor = west] {$s_3$};%
        \draw (1.2, -0.25) node[anchor = west] {$s_4$};%
        \halfedge{0.5}
        \draw (1.2, 0.5) node[anchor = west] {$s_1$};%
    \end{tikzpicture}
    & \\
    \cline{2-5}
    & \multirow{3}{4cm}{$d_k(w_1) = 3$ and $w_1 = w_2$}
    & \begin{tikzpicture}
        \halfedge[3]{0}
        \draw (1.2, 0) node[anchor = west] {$2s$};%
    \end{tikzpicture}
    & \begin{tikzpicture}
        \halfedge[3]{0}
        \draw (1.2, 0) node[anchor = west] {$s$};%
    \end{tikzpicture}
    & \multirow{3}{3cm}{$a = 2$} \\
    & & \begin{tikzpicture}
        \twohalfedges[2]{0}
        \draw (1.2, 0.15) node[anchor = west] {$s_1 + 2s_2$};%
        \draw (1.2, -0.15) node[anchor = west] {$s_1$};%
    \end{tikzpicture}
    & \begin{tikzpicture}
        \twohalfedges[2]{0}
        \draw (1.2, 0.15) node[anchor = west] {$s_1$};%
        \draw (1.2, -0.15) node[anchor = west] {$s_2$};%
    \end{tikzpicture}
    & \\
    & & \begin{tikzpicture}
        \threehalfedges{0}
        \draw (1.2, 0.25) node[anchor = west] {$s_1 + s_2$};%
        \draw (1.2, 0) node[anchor = west] {$s_1 + s_3$};%
        \draw (1.2, -0.25) node[anchor = west] {$s_2 + s_3$};%
    \end{tikzpicture}
    & \begin{tikzpicture}
        \threehalfedges{0}
        \draw (1.2, 0.25) node[anchor = west] {$s_1$};%
        \draw (1.2, 0) node[anchor = west] {$s_2$};%
        \draw (1.2, -0.25) node[anchor = west] {$s_3$};%
    \end{tikzpicture}
    & 
\end{tabular}}
        \caption{Case-by-case proof that at every vertex $v$ the sum of outgoing slopes of $\chi^*\omega$ on edges in $\tilde\Gamma$ over $h \in E(K)$ is an integer linear combination of the corresponding sums at the vertices $w_1$ and $w_2$.}
        \label{tab:harmonicity}
    \end{figure}
    
    Now let us check that $\chi\iota = - \chi$. Let $p' \in \tilde\Gamma$ correspond to $p_1 + p_2$ with $k(p_1) = k(p_2)=p$. Then $\iota(p')$ corresponds to $p_3 + p_4$, where $p_1 + p_2 + p_3 + p_4$ is the fiber of $k$ over $p \in K$. Therefore,
    \[ \chi(\iota(p')) + \chi(p') = p_1 + p_2 + p_3 + p_4 - 2L \sim k^{-1}(p) - D \]
    by construction of $L$. The right hand side is now linearly equivalent to 0 because by definition in Equation~\eqref{eq:bigonal_D}, $D$ is a fiber of $k$ and all fibers of $k$ are linearly equivalent (here we use that $K$ is a tree). \qedhere
\end{proof}

\begin{lemma}
    \label{lem:key_formula}
    Keep the notation as in the proof of Theorem \ref{thm:tropical_Recillas_theorem} above. There is a tropical 1-cycle $A : \Pi^2 \to \RR$, which agrees with $A'$ away from a 0-dimensional locus. Moreover, we have
    \[
    [\Delta_\Pi] + A = (k \times k)^\ast [\Delta_K]\in Z_1(\Pi^2).
    \]
\end{lemma}

%\Dmitry{we actually show that these are equal as cycles, not just homology classes, hence $A$ is balanced and therefore Lemma 5.4 is superfluous}
%\Felix{I have modified the wording of the lemma and the corresponding part of the proof of the theorem to stress this. I think we are now good to remove the old Lemma 5.4, which I have already commented out.}

\begin{proof}
    In the proof of Theorem~\ref{thm:tropical_Recillas_theorem} we will think of $(k \times k)^\ast[\Delta_k]$ as a pullback of cohomology cycles. However, in order to compute $(k \times k)^\ast[\Delta_k]$ as a tropical 1-cycle, we give a different interpretation. Because $[\Delta_K]$ is a cycle in codimension 1 and $K \times K$ is smooth, we can identify $[\Delta_K]$ with a Cartier divisor $D \in \Div(K^2)$, which admits a pullback to $\Div(\Pi^2)$. This is a tropical 1-cycle in $Z_1(\Pi^2)$ after intersection with $[\Pi^2]$. Let us quickly verify, that the induced notion of pullback of cycles is the same as pullback on cohomology after applying the cycle class map. %\Felix{Am I being childish for checking this or is this really worth stating?} \Dmitry{let's leave it}
    
    Starting with $D \in \Div(K^2)$, we have 
    \[ \cyc[\Delta_K] = \cyc\big(D \cdot [K \times K] \big) = \cyc[K \times K] \frown c_1(\calL(D)) , \]
    in other words, the Poincar\'e dual of $\cyc[\Delta_K]$ is $c_1(\calL(D))$. Similarly, we treat the pullback $(k\times k)^*D$
    \[ \cyc \Big( \big( (k\times k)^*D \big) \cdot [\Pi \times \Pi] \Big) = \cyc[\Pi \times \Pi] \frown c_1 \Big(\calL \big((k\times k)^* D \big) \Big) . \]
    But now pullback of Cartier divisors commutes with taking the line bundle \cite[Proposition~3.15]{GrossShokrieh_homology} and that in turn commutes with the Chern class map \cite[Proposition~5.11]{GrossShokrieh_homology}. Hence the two interpretations of the pullback coincide. \medskip
    
    Let us now prove the formula that we claimed.
    Let $e$ be an edge of $K$, which we identify with the interval $[0, \ell(e)]$, and consider the cell $e \times e \subseteq K \times K$. The diagonal of this cell is defined as a Cartier divisor by the piecewise affine (i.e. rational) function 
    \begin{align*}
        f : e \times e &\longrightarrow \RR, \\
        (x,y) &\longmapsto \max\{y-x, 0\} . 
    \end{align*}
    Let $k^{-1}(e) = \{ e_1, \ldots, e_m \}$ be the fiber over $e$ with $m \in \{2,3,4\}$. Without loss of generality we assume that the dilation factors along $e_2, \ldots, e_m$ are 1. Consider the cell $e_i \times e_j \subseteq \Pi \times \Pi$ and the pullback of $f$ restricted to this cell is
    \begin{align*}
        g = \big(f \circ (k \times k) \big) \big|_{e_i \times e_j} : e_i \times e_j &\longrightarrow \RR \\ 
        (x,y) &\longmapsto \max \big\{ d_k(e_j) y - d_k(e_i) x, \ 0 \big\} .
    \end{align*}
    In order to determine the multiplicity of $\div(g)$ along the support of $g$, we have to evaluate $g$ on a lattice normal vector of the support. 
    The support of $\div(g)$ has primitive integer tangent vector $\frac{1}{\gcd \big(d_k(e_i), d_k(e_j) \big)} \begin{pmatrix} d_k(e_j) \\ d_k(e_i) \end{pmatrix}$ and a lattice normal vector is determined by completing this vector into a $\ZZ$-basis for $\ZZ^2$. This is symmetric in $i$ and $j$, so we may assume without loss of generality $i \leq j$ and find the lattice normal vector $\begin{pmatrix} 0 \\ 1 \end{pmatrix}$ and the slope of $g$ in this direction is $d_k(e_j)$. This means that the value of $(k \times k)^*[\Delta_K]$ on $\big\{ (x,y) \in e_i \times e_j \bigmid k(x) = k(y) \big\}$ is $d_k(e_1)$ for $i = j = 1$ and $1$ otherwise. This is precisely the same as the value of $[\Delta_\Pi] + A'$ on this locus. By the definition of the group structure of $Z_1(\Pi^2)$ this shows that there is a tropical 1-cycle $A = (k\times k)^\ast [\Delta_K] - [\Delta_\Pi] \in Z_1(\Pi^2)$ which differs from $A'$ only in a 0-dimensional locus. \qedhere
\end{proof}

\begin{lemma}
    \label{prop:diagonal_formula}
    Let $K$ be a smooth metric tree (so that each extremal edge is infinite), and let $p \in K$ be any point. 
    Then
    \[ [\Delta_K] = [K \times p] + [p \times K] \in H_{1,1}(K^2) , \]
    where $\Delta_K$ is the diagonal and $[\cdot]$ denotes the fundamental class.
\end{lemma}

\begin{proof}
    We note that the terms in the above equation are tropical homology classes associated to 1-cycles (see Example \ref{ex:diagonal_cycle} for the diagonal cycle), and that the equation certainly does not hold in  in $Z_1(K^2)$.
    
    We proceed in two steps.
    As a base case, we first prove the claim for the compactified real line $L = \RR \cup \{-\infty, \infty\}$.
    We orient $L$ from $-\infty$ to $\infty$. 
    Let $\eta$ be the generator of $\Omega^1_{L}(\RR) \cong \ZZ$ compatible with this orientation, more precisely
    \begin{equation*}
        \Omega^1_{L, x} = \begin{cases}
            \langle \eta\rangle & \text{if } x \in \RR, \\
            0 & \text{otherwise.}
        \end{cases}
    \end{equation*}
    Then $\eta_1 = \begin{pmatrix} \eta \\ 0 \end{pmatrix}$ and $\eta_2 = \begin{pmatrix} 0 \\ \eta \end{pmatrix}$ provide generators for the stalks of $\Omega^1_{L^2}$, or again more precisely
    \begin{equation*}
        \Omega^1_{L^2, (x,y)} = \begin{cases}
            \langle \eta_1, \eta_2 \rangle & \text{if } (x,y) \in \RR^2 \\
            \langle \eta_1 \rangle & \text{if } x \in \RR \text{ and } y = \pm \infty \\
            \langle \eta_2 \rangle & \text{if } y \in \RR \text{ and } x = \pm \infty \\
            0 & \text{otherwise.}
        \end{cases}
    \end{equation*}
    Let $\sigma : \Delta^2 \to \big\{(x,y) \in L^2 \bigmid x \geq y \big\}$ be a singular 2-simplex parametrizing the area below the diagonal of $L^2$ with orientation compatible with that on $L$. Write
    \begin{align*}
        \tau_1 &: \Delta^1 \longrightarrow \big\{(x, -\infty) \in L^2 \big\} \\
        \tau_2 &: \Delta^1 \longrightarrow \big\{(\infty, y) \in L^2 \big\} \\
        \delta &: \Delta^1 \longrightarrow \big\{(x,x) \in L^2 \big\}
    \end{align*}
    for the restriction of $\sigma$ to the faces of $\Delta^2$, so that the boundary of $\sigma$ as a singular chain is $\partial \sigma = \tau_1 - \delta + \tau_2$.
    Then $B = \sigma \otimes (\eta_1^* + \eta_2^*)$ is a $(1,2)$-chain and
    \begin{align*}
        \partial B &= \tau_1 \otimes (\eta_1^* + \underbrace{\eta_2^*}_{\mathrlap{= 0 \text{ on } \Im \tau_1}}) - \delta \otimes (\eta_1^* + \eta_2^*) + \tau_2 \otimes (\underbrace{\eta_1^*}_{\mathrlap{= 0 \text{ on } \Im \tau_2}} + \eta_2^*) \\
        &= \cyc [L \times -\infty] - \cyc [\Delta_L] + \cyc [\infty \times L] .
    \end{align*}
    Note that the second equality only holds because the images of $\tau_1$ and $\tau_2$ are contained in the boundary of the rational polyhedral space $L^2$. If we did the same computation with a finite interval $L = [a,b]$, we would not see $\cyc [L \times a]$ and $\cyc [b \times L]$ in $\partial \big(\sigma \otimes (\eta_1^* + \eta_2^*) \big)$ because $\eta_1^*$ and $\eta_2^*$ do not vanish on the topological boundary.
    
    To finish the proof for $L$, we want to argue that $\cyc [L \times \infty ] = \cyc [L \times p]$ for any $p \in \RR$. Let $\sigma'$ be a singular 2-complex (consisting of two simplices) parametrizing $\big\{(x,y) \in L^2 \bigmid -\infty \leq y \leq p \big\}$. %\Dmitry{We probably need a pair of simplices? Two triangles making a square?} \Felix{No, why would we?}
    Then 
    \[ \partial (\sigma' \otimes \eta_1^*) = \cyc [L \times -\infty] - \cyc [L \times p] . \]
    The key here is that $\eta_1^*$ vanishes on $[-\infty, p] \times \{\pm \infty\}$ so that there is no contribution to $\partial(\sigma' \otimes \eta_1^*)$ from the remaining two edges of $\partial \sigma'$. 
    Similarly we see $[\infty \times L] = [p \times L]$ which completes the proof for the base case. \medskip
    
    Now let $K$ be any smooth tree. Fix two distinct points $-\infty, \infty \in K$ on the boundary and let $L \subseteq K$ be the path from $-\infty$ to $\infty$. As a rational polyhedral space, $L$ is isomorphic to the compactified real line from the base case. Furthermore, $L$ is a deformation retract of $K$. Denote the retraction map $\rho : K \to L$. This is a proper map of rational polyhedral spaces, so there is a pushforward map $\rho_*$ in homology.
    
    We claim that the induced map $\rho^2_* : H_{1,1}(K^2) \to H_{1,1}(L^2)$ is an isomorphism. We first show that $\rho_*$ is an isomorphism. Indeed, $H_{1,0}$ and $H_{0,1}$ of $K$ and $L$ are trivial because $K$ and $L$ are trees. On the other hand, the pushforward map on  $H_{0,0}$ is an isomorphism sending the class of a point to the class of a point, and similarly on $H_{1,1}$ we have an isomorphism given by sending the fundamental class of $K$ to the fundamental class of $L$. By the K\"unneth formula, $H_{1,1}(K^2)$ decomposes as
\[        H_{0,0}(K) \otimes H_{1,1}(K) 
        \quad \oplus \quad
        \underbrace{H_{0,1}(K) \otimes H_{1,0}(K)}_{=0}
        \quad \oplus \quad
        \underbrace{H_{1,0}(K) \otimes H_{0,1}(K)}_{=0} 
        \quad \oplus \quad
        H_{1,1}(K) \otimes H_{0,0}(K) ,
\]    and similarly for $L$. Hence $\rho_*$ being an isomorphism implies that $\rho^2_* : H_{1,1}(K^2) \to H_{1,1}(L^2)$ is an isomorphism as well.

We have already proved our claim in $H_{1,1}(L^2)$, so it suffices to show that
\begin{align*}
    \rho^2_* \cyc [\Delta_K] &= \cyc [\Delta_L], \\
    \rho^2_* \cyc [K \times p] &= \cyc [L \times \rho(p)], \quad \text{ and } \\
    \rho^2_* \cyc [p \times K] &= \cyc [\rho(p) \times L] .
\end{align*}
    But the cycle class map and pushforward commute, so we may verify this on the level of tropical cycles, where it is easy to see. On every edge $e$ of $K$, the map $\rho$ is either the identity (if $e$ is part of $L$) or constant. Hence the index of tangent spaces in Equation \eqref{eq:definition_push_forward} is either 1 or 0, i.e. the part of the cycles $[\Delta_K]$, $[K \times p]$, and $[p \times K]$ that is already in $L$ survives the pushforward while the rest is weighted with 0. This finishes the proof. \qedhere
\end{proof}

We now consider the bigonal construction. Unlike the trigonal construction, the two double covers involved in the bigonal construction are always dilated, hence it is necessary to distinguish the divisorial and continuous Pryms. The latter carry natural principal polarizations, while the former do not. It turns out that the non-principally polarized divisorial Pryms are related by taking the dual.

\begin{theorem}[Theorem~\ref{thm:bigonal_construction}] \label{thm:bigonal_construction_restated}
    Let $\tilde \Gamma \overset{\pi}{\longrightarrow} \Gamma \overset{f}{\longrightarrow} K$ be a tower of harmonic morphisms of metric graphs of degrees $\deg \pi = \deg f = 2$, where $K$ is a metric tree. Assume that there is no point $x\in K$ with the property that $|f^{-1}(x)|=2$ and $|(f\circ\pi)^{-1}(x)|=2$. Then the output $\tilde \Pi \overset{\pi'}{\longrightarrow} \Pi \overset{f'}{\longrightarrow} K$ of the bigonal construction has the same property, and applying the bigonal construction to it reproduces the original tower. If moreover $\tilde \Gamma$ and $\tilde \Pi$ are both connected, then there is an isomorphism of polarized tropical abelian varieties
    \[ \Prym_d(\tPi/\Pi)^\vee\cong \Prym_d(\tGa/\Ga),\]
    where the polarization on $\Prym_d(\tPi/\Pi)^\vee$ is the pullback of the principal polarization on $\Prym_c(\tPi/\Pi)^\vee$. 
\end{theorem}

The algebraic version of Theorem~\ref{thm:bigonal_construction_restated} \cite{Pantazis} was proved by reduction to the tetragonal construction. Since this is not an option for us, we give the following proof which is an adaptation of the proof of Theorem~\ref{thm:tropical_Recillas_theorem_restated} that we just discussed.

\begin{proof} Recall that we have already established the properties of the bigonal construction in Propositions~\ref{prop:bigonal_properties} and~\ref{prop:bigonalconnected}, and it only remains to compare the Pryms. As for Theorem~\ref{thm:tropical_Recillas_theorem_restated}, we give an outline of the proof and spin off the necessary technical calculations into a series of lemmas. 

By construction, the curve $\tilde \Pi$ comes with an embedding in $\Div_2^+(\tilde \Gamma)$, compare Equation~\eqref{eq:tildePi_in_Symn}. As before, let $D \in \Div_4^+(\tilde\Gamma)$ be a fiber of $f \circ \pi$ and choose $M \in \Div_1(\tilde\Gamma)$ such that $4M = D$. Using $M$ as the base point we define $\chi : \tilde \Pi \to \Prym_c(\tilde \Gamma / \Gamma)$ to be the restriction of the second power of the Abel--Prym map $\psi_M : \tilde \Gamma \to \Prym_c(\tilde \Gamma / \Gamma)$ to $\tilde\Pi$. 
    Similarly to Lemma~\ref{lem:alpha}, one can check that $\chi$ is a morphism of rational polyhedral spaces and the choice of $M$ ensures $\chi \circ \iota = -\chi$. Moreover, in Lemma~\ref{lem:pullback_condition} we verify that the pullback of 1-forms along $\chi$ takes values in the image of $\Id - \iota^* : \Omega^1_{\tilde\Pi} \to \Omega^1_{\tilde\Pi}$. Fixing a base point $q' \in \tilde\Pi$, we obtain a commutative diagram
    \begin{equation*}
        \begin{tikzcd}
            \tilde \Pi \arrow[r, "\chi"] \arrow[d, "\psi_{q'}"] & \Prym_c(\tilde \Gamma/\Gamma) \arrow[d, "t_{-\chi(q')}"] \\
            \Prym_c(\tilde \Pi / \Pi) \arrow[r, "\nu"] & \Prym_c(\tilde \Gamma / \Gamma) 
        \end{tikzcd}
    \end{equation*}
    by the universal property of the continuous Prym variety (Proposition~\ref{prop:universal_property_Prymc}). 
    
    We recall (see Proposition~\ref{prop:Prymallpropetries}) that the continuous and divisorial Prym varieties associated to the double covers $\pi:\tGa\to \Ga$ and $\pi':\tPi\to \Pi$ are related by natural free isogenies 
    \[
    \gamma: \Prym_c(\tGa/\Ga) \longrightarrow \Prym_d(\tGa/\Ga),\qquad \gamma' : \Prym_c(\tPi/\Pi) \longrightarrow \Prym_d(\tPi/\Pi).
    \]
    The divisorial Pryms carry polarizations $\zeta\in H^{1,1}\big(\Prym_d(\tGa/\Ga)\big)$ and $\zeta'\in H^{1,1}\big(\Prym_d(\tPi/\Pi) \big)$ induced from their respective Jacobians. The continuous Pryms carry principal polarizations $\zeta_c\in H^{1,1}\big(\Prym_c(\tGa/\Ga)\big)$ and $\zeta'_c\in H^{1,1}\big(\Prym_c(\tPi/\Pi) \big)$ such that $\gamma^\ast \zeta = 2\zeta_c$ and $(\gamma')^*(\zeta')=2\zeta'_c$.

    Let $\Prym_c(\tGa/\Ga)^{\vee}$ be the dual pptav with principal polarization $\zeta^{-1}_c$. To avoid overloading notation, we denote by $\xi : \Prym_c(\tGa/\Ga) \to \Prym_c(\tGa/\Ga)^\vee$ the isomorphism induced by the principal polarization $\zeta_c$. The maps $\nu$, $\xi$, $\ga'$, and the dual $\ga^{\vee}$ of $\ga$ fit into the diagram
            \begin{equation} \label{eq:diagram_beta}
            \begin{tikzcd}
                \Prym_c(\tPi/\Pi) \arrow[r, "\nu"] \arrow[d, "\gamma'"] & \Prym_c(\tGa/\Ga) \arrow[r, "\xi"] & \Prym_c(\tGa/\Ga)^{\vee} \\
                \Prym_d(\tPi / \Pi) \arrow[rr, dotted,"\delta"] && \Prym_d(\tilde \Gamma / \Gamma)^{\vee}  \arrow[u, "\gamma^{\vee}"].
            \end{tikzcd}
        \end{equation}
    We claim that there exists an isomorphism $\delta:\Prym_d(\tPi / \Pi)\to \Prym_d(\tGa / \Ga)^{\vee}$ that makes this diagram commute. Furthermore, we claim that the pullback along $\delta$ of the induced polarization $(\gamma^{\vee})^*(\zeta_c^{-1})$ on $\Prym_d(\tGa/\Ga)^{\vee}$ is equal to $\zeta'$. 

    First, we show in Lemma~\ref{lem:beta} that there exists a morphism $\delta$ that makes~\eqref{eq:diagram_beta} commute. To prove that it is an isomorphism, we compute its total degree (the product of the dilation and geometric degrees, see Definition~\ref{def:types_of_homomorphisms}). In Lemma~\ref{lem:pushbigonal}, we use similar ideas as in the proof of Theorem~\ref{thm:tropical_Recillas_theorem_restated} to show that 
    \[
    \nu_\ast ((\zeta'_c)^{h-1}) = 2\zeta_c^{h-1}\in H_{1,1}\big(\Prym_c(\tGa/\Ga)),
    \]
    where $h$ is the common dimension of the Pryms. Proposition~\ref{prop:degree_pptavs_via_homology} then implies that the total degree of $\nu$ is equal to $\deg \nu = 2^h$. 

    On the other hand, Proposition~\ref{prop:Prymallpropetries} states that the morphisms $\ga$ and $\ga'$ are free isogenies of degrees $2^{d-1}$ and $2^{d'-1}$, respectively, where $d$ and $d'$ are the dilation indices of $\pi:\tGa\to \Ga$ and $\pi':\tPi\to \Pi$. Taking the dual exchanges dilation and geometric degrees, therefore 
    \begin{align*}
        \deg \ga'&=\deg_g\ga' \cdot \deg_d \ga'=2^{d'-1}\cdot 1=2^{d'-1},\\ 
        \deg \ga^{\vee}&=\deg_g\ga^{\vee} \cdot \deg_d \ga^{\vee}=1\cdot 2^{d-1}=2^{d-1}.
    \end{align*}
    By Proposition~\ref{prop:bigonalconnected}, we have $d+d'-2=h$, therefore
    \[
    \deg (\xi\circ \nu)=\deg \ga'\deg \delta\deg \ga^{\vee}=2^h\deg \delta.
    \]
    Since $\xi$ is an isomorphism, we see that $\deg \delta=1$ and hence $\delta$ is an isomorphism.
    
    Finally, we compute the pullbacks of the polarizations along all maps in~\eqref{eq:diagram_beta} . The principal polarization $\zeta_c^{-1}$ on $\Prym_c(\tGa/\Ga)^{\vee}$ pulls back to $\xi^*(\zeta_c^{-1})=\zeta_c$ on $\Prym_c(\tGa/\Ga)$ and hence to
    \[(\xi\circ \nu)^*(\zeta_c^{-1})=\nu^*(\zeta_c)=2\zeta'_c
    \]
    on $\Prym_c(\tPi/\Pi)$ by Proposition~\ref{prop:degree_pptavs_via_homology}. Since $(\ga')^*(\zeta')=2\zeta'_c$, it must be that $\zeta'=\delta^*\big((\gamma^\vee)^*(\zeta^{-1}_c)\big)$, therefore $\delta$ is in fact an isomorphism of polarized tropical abelian varieties. \qedhere

\end{proof}

\begin{example}
    Let us return to the towers of graphs depicted in Figure~\ref{fig:bigonal_example_2}. We assign real edge lengths $a$, $b$, and $c$ to the three edges of $K$ and thus obtain two $(2,2)$-towers $\tilde \Gamma_1 \to \Gamma_1 \to K$ and $\tilde \Gamma_2 \to \Gamma_2 \to K$ of tropical curves. To compute the (non-principally polarized) divisorial Prym varieties $\Prym_d(\tGa_1/\Ga_1)$ and $\Prym_d(\tGa_2 / \Ga_2)$, we work with the bases
    \begin{align*}
        \Ker \left(\pi_*:H_1(\tilde \Gamma_1,\ZZ) \to H_1(\Gamma_1,\ZZ)\right) &= \langle \teta_1^+ - \teta_1^-, \teta_2 \rangle \\
        \Ker\left(\pi'_*:H_1(\tilde \Gamma_2,\ZZ) \to H_1(\Gamma_2,\ZZ)\right) &= \langle \tilde \epsilon_1, \tilde\epsilon_2^+ - \tilde \epsilon_2^- \rangle
    \end{align*}
    depicted in Figure~\ref{fig:basis_bigonal_example}. For $i  = 1,2$ let $\zeta_i : H_1(\tilde\Gamma_i, \ZZ) \to \Omega_{\tilde\Gamma_i}^1(\tilde\Gamma_i)$ be the principal polarization of $\Jac(\tilde\Gamma_i)$. With this, the pairings induced from the integration pairings on the $\tilde \Gamma_i$ are \medskip
    
    \begin{center}
        \begin{tabular}{c|c c}
            $[\cdot, \cdot]_1$ & $\zeta_1(\teta_1^+ - \teta_1^-)$ & $\zeta_1(\teta_2)$ \\
            \hline
            $\teta_1^+ - \teta_1^-$ & $2b$ & $b$ \\
            $\teta_2$ & $2b$ & $a+2b+c$
        \end{tabular}
        \qquad and \qquad
        \begin{tabular}{c|c c}
            $[\cdot, \cdot]_2$ & $\zeta_2(\tilde \epsilon_1)$ & $\zeta_2(\tilde \epsilon_2^+ - \tilde\epsilon_2^-)$ \\
            \hline
            $\tilde \epsilon_1$ & $2b$ & $2b$ \\
            $\tilde \epsilon_2^+ - \tilde\epsilon_2^-$ & $b$ & $a+2b+c$
        \end{tabular}
    \end{center}  
    \medskip
    and we clearly see that transposing one of these matrices gives the other. This shows that $\Prym_d(\tGa_1/\Ga_1) \cong \Prym_d(\tGa_2/\Ga_2)^\vee$ on the level of integral tori. But even more is true: the polarization of $\Prym_d(\tGa_2/\Ga_2)^\vee$ induced by the principal polarization on $\Prym_c(\tGa_2/\Ga_2)^{\vee}$ is given by 
    \begin{align*}
        \xi_2^\vee : \Omega^1_{\tGa_2}(\tilde\Gamma_2) &\longrightarrow H_1(\tilde\Gamma_2, \ZZ) \\
        \zeta_2(\tilde \epsilon_1) &\longmapsto 2\tilde\epsilon_1 \\
        \zeta_2(\tilde\epsilon_2^+ - \tilde\epsilon_2^-) &\longmapsto \tilde\epsilon_2^+ - \tilde\epsilon_2^-
    \end{align*}
    and hence the pull-back of $\xi_2^\vee$ along the isomorphism $\Prym_d(\tilde\Gamma_1 / \Gamma_1) \to \Prym_d (\tilde\Gamma_2 / \Gamma_2)^\vee$ yields the induced polarization $\xi_1$ on $\Prym_d(\tGa_1/\Ga_1)$. %Hence, we have an isomorphism of principally polarized tropical abelian varieties.
    
    \begin{figure}[htb]
        \centering
        \begin{subfigure}{0.4\textwidth}
            \centering
            \begin{tikzpicture}[x=15mm, y=15mm]
                %Generator 1^+
                \draw (-0.5, 0.5) node {$\teta_1^+$};
            	\drawExampleBigonal
            	\drawCycle{
            		\draw[postaction = decorate] (1,0) arc (180:540:0.5 and 0.2);
            	}
            \end{tikzpicture}
        \end{subfigure}
        \qquad
        \begin{subfigure}{0.4\textwidth}
            \centering
            \begin{tikzpicture}[x=15mm, y=15mm]
                %Generator 1
                \draw(-0.5, 0.5) node {$\tilde \epsilon_1$};
                \drawExampleBigonalOutput
                \drawCycle{
                    \draw[postaction = decorate] (1,0) arc (180:540:0.5 and 0.2);
                }
            \end{tikzpicture}
        \end{subfigure}
        \hspace{1cm}
        \begin{subfigure}{0.4\textwidth}
            \centering
            \begin{tikzpicture}[x=15mm, y=15mm]
                %Generator 1^-
                \draw (-0.5, 0.5) node {$\teta_1^-$};
            	\drawExampleBigonal
            	\drawCycle{
            		\draw[postaction = decorate] (1,-0.5) arc (180:-180:0.5 and 0.2);
            	}
            \end{tikzpicture}
        \end{subfigure}
        \qquad
        \begin{subfigure}{0.4\textwidth}
            \centering
            \begin{tikzpicture}[x=15mm, y=15mm]
                %Generator 2^+
                \draw(-0.5, 0.5) node {$\tilde \epsilon_2^+$};
                \drawExampleBigonalOutput
                \drawCycle{
                    \path[draw, postaction = decorate] (1, 0) -- ++ (-1, 0) -- ++ (1, 0.5) -- ++ (1,0) -- ++ (1, -0.5) -- ++ (-1, 0) arc (0:180:0.5 and 0.2);
                }
            \end{tikzpicture}
        \end{subfigure}
        \hspace{1cm}
        \begin{subfigure}{0.4\textwidth}
            \centering
            \begin{tikzpicture}[x=15mm, y=15mm]
                %Generator 2
                \draw (-0.5, 0.5) node {$\teta_2$};
            	\drawExampleBigonal
            	\drawCycle[0.41]{
            		\path[draw, postaction = decorate] (1,0) -- ++ (-1, -0.25) -- ++ (1, -0.25) arc (180:360:0.5 and 0.2) -- ++ (1, 0.25) -- ++ (-1, 0.25) arc (0:180:0.5 and 0.2);
            	}
            \end{tikzpicture}
        \end{subfigure}
        \qquad
        \begin{subfigure}{0.4\textwidth}
            \centering
            \begin{tikzpicture}[x=15mm, y=15mm]
                %Generator 2^-
                \draw(-0.5, 0.5) node {$\tilde \epsilon_2^-$};
                \drawExampleBigonalOutput
                \drawCycle{
                    \path[draw, postaction = decorate] (1,0) -- ++ (-1, 0) -- ++ (1, -0.5) -- ++ (1, 0) -- ++ (1, 0.5) -- ++ (-1, 0) arc (180:360:-0.5 and 0.2);
                }
            \end{tikzpicture}
        \end{subfigure}
        \caption{Bases of the first homology $H_1(\tilde\Gamma_1, \ZZ)$ and $H_1(\tilde \Gamma_2, \ZZ)$ used in the computation of the Prym varieties. Recall that the involutions on $\tilde\Gamma_1$ and $\tilde\Gamma_2$ are given by flipping the pictures along their horizontal symmetry.}
        \label{fig:basis_bigonal_example}
    \end{figure}
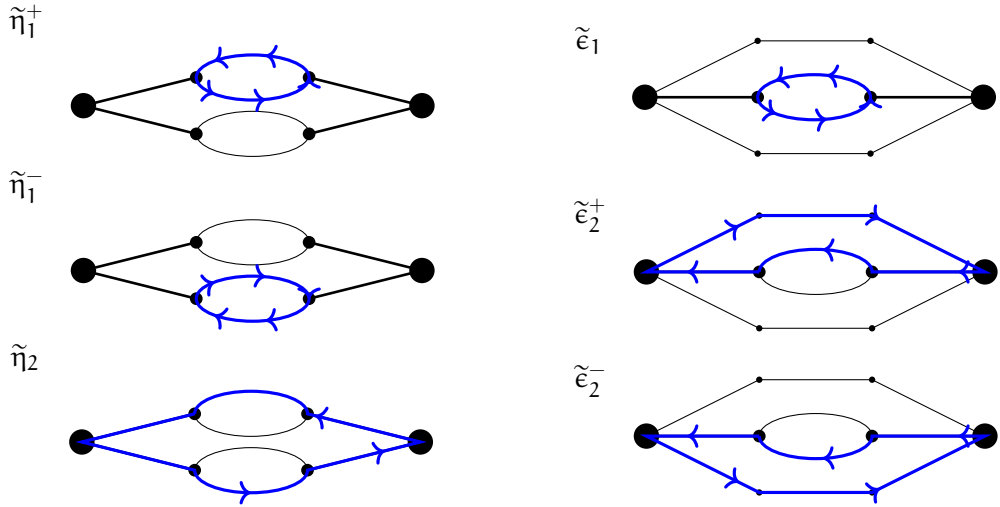
\end{example}

We now verify the details in the proof of Theorem~\ref{thm:bigonal_construction_restated}.

\begin{lemma} \label{lem:pullback_condition}
    In the notation of Theorem~\ref{thm:bigonal_construction_restated}, the map $\chi : \tilde\Pi \to \Prym_c(\tilde\Gamma / \Gamma)$ satisfies the pullback condition on cotangent sheaves from Proposition~\ref{prop:universal_property_Prymc}, in other words $\chi^{-1}\Omega^1_{\Prym_c(\tilde\Gamma/\Gamma)} \to \Omega^1_{\tilde\Pi}$ takes values in the image of $\Id - \iota'^\ast : \Omega^1_{\tilde\Pi} \to \Omega^1_{\tilde\Pi}$.
\end{lemma}

\begin{proof}
    We denote by $\iota'$ the involution on $\tilde \Pi$ and we denote by $\iota$ the involution on $\Jac(\tilde\Gamma)$ induced by the involution of $\tilde\Gamma$. Note that the map $\chi$ factors as $\chi_0 : \tilde\Pi \to \Jac(\tilde\Gamma)$ followed by $\Id - \iota$. Using Figure~\ref{fig:bigonaltypes} it is easy to verify that 
    \begin{equation} \label{eq:chi_iota_commute}
        \iota\circ \chi_0 = \chi_0 \circ \iota'
    \end{equation}
    and we claim that this implies the statement of the lemma.     
    Indeed, let $\omega \in \Omega^1_{\Prym_c(\tilde\Gamma / \Gamma)}$. Then 
    \begin{align*}
        \chi^\ast \omega &= \chi_0^\ast ((\Id - \iota)^\ast \omega) \\
        &= \chi_0^\ast \omega - \chi_0^\ast (\iota^\ast \omega) \\
        &= \chi_0^\ast \omega - \iota'^\ast (\chi_0^\ast \omega) \\
        &= (\Id - \iota')^\ast (\chi_0^\ast \omega). \qedhere
    \end{align*}
\end{proof}

\begin{lemma} \label{lem:beta}
    There exists a homomorphism $\delta:\Prym_d(\tPi / \Pi)\to \Prym_d(\tilde \Gamma / \Gamma)^{\vee}$ that makes Diagram~\eqref{eq:diagram_beta} commute.
\end{lemma}

\begin{proof} We recall that the lattices $\Prym_c(\tPi/\Pi)=(\La_1,\La_1', [\cdot, \cdot]_\Pi)$, $\Prym_c(\tGa/\Ga)=(\La_2,\La_2', [\cdot, \cdot]_\Gamma)$ of the continuous Pryms are
\[
\La_1=(\Coker (\pi'^*))^{\mathrm{tf}},\qquad \La_1'=\Im(\Id - \iota'_*),\qquad 
\La_2=(\Coker (\pi^*))^{\mathrm{tf}},\qquad \La_2'=\Im(\Id - \iota_*).
\]
Here $\iota:\tGa\to \tGa$ and $\iota':\tPi\to \tPi$ are the involutions associated to the double covers $\pi$ and $\pi'$, respectively. The divisorial Pryms are
\[
\Prym_d(\tPi/\Pi)= \big(\La_1,(\La_1')^{\sat}, [\cdot, \cdot]_\Pi \big),\qquad  \Prym_d(\tGa/\Ga)= \big(\La_2,(\La_2')^{\sat}, [\cdot, \cdot]_\Gamma \big),
\]
where the saturations are
\[
(\La'_1)^{\sat}=\Ker \pi'_*,\qquad (\La'_2)^{\sat}=\Ker \pi_*.
\]
The homomorphism $\nu$ consists of the homomorphisms
\[
\nu^\#:\La_2\longrightarrow \La_1 \qquad \text{and} \qquad \nu_\#:\La'_1\longrightarrow \La'_2.
\]
Hence the homomorphisms $\delta^\#$ and $\delta_\#$ comprising $\delta$, if they exist, are the extensions of respectively $\nu^\#\circ \zeta_c$ and $\nu_\#$ to the saturations of their domains: this determines them uniquely, and the consistency condition~\eqref{eq:homomorphism_condition} follows automatically. In other words, we need to show the following:
\begin{enumerate}
    \item The homomorphism $\nu_\#:\La'_1\to \La'_2$ extends to $(\La'_1)^{\sat}=\Ker \pi'_*$.
    \item The homomorphism $\nu^\#\circ \zeta_c:\La'_2\to\La_1$ extends to $(\La'_2)^{\sat}=\Ker \pi_*$.
\end{enumerate}

We start with the first claim. According to Proposition~\ref{prop:dilatedbasis} (see also~\eqref{eq:Prymcbasis2}), $\Lambda_1'$ has a basis consisting of elements of the form $\tal^+_i-\tal^-_i$ and $2\tbe_j$, so that the $\tal^+_i-\tal^-_i$ and $\tbe_j$ form a basis for $(\Lambda_1')^{\sat}$. Recall from the proof of Proposition~\ref{prop:universal_property_Prymc} that $\nu_\#(2\tbe_j)=\mu_\#(\tbe_j)$, where $\mu_\#:\Jac(\tPi)\to \Prym_c(\tGa/\Ga)$ is the homomorphism corresponding the map $\chi:\tPi\to \Prym_c(\tGa/\Ga)$. Therefore, to show that $\nu_\#$ extends to $(\Lambda_1')^{\sat}$, it is sufficient to show that each $\mu_\#(\tbe_j)$ is divisible by $2$.

We consider an arbitrary simple cycle $\tbe\in H_1(\tPi,\ZZ)$ such that $\iota_*\tbe=-\tbe$. Such a cycle does not contain dilated edges, and we can assume without loss of generality that it contains exactly two dilated vertices $\tu$ and $\tv$ of $\tPi$, so that $\iota$ folds it in half and that the image of its support in $\Pi$ is a simple path from $\pi'(\tu)$ to $\pi'(\tv)$ (see Figure~\ref{fig:edgelabels_in_beta}). Recall that points of $\tPi$ are pairs of points of $\tGa$ whose images in $\Ga$ form a fiber of $f:\Ga\to K$. Hence we can write
\[
\tbe = \sum_{i = 1}^n \left[\te_i + \te'_i - \iota(\te_i)-\iota(\te'_i)\right],
\]
where $\te_i$ and $\te'_i$ are edges of $\tGa$ such that $\te_i+\te'_i$ defines an edge of $\tPi$ (in other words, $\pi(\te_i)+\pi(\te'_i)=f^*(e_i)$ for some edge $e_i\in E(K)$), and $\iota$ is the involution on $\tGa$. We note that potentially $\te_i=\te_j$ or $\te_i=\te'_j$ for $i\neq j$. We further assume that the edges are labeled and oriented in such a way that 
\begin{equation}
    t(\te_i) = s(\te_{i+1})\qquad \text{and} \qquad  t(\te'_i) = s(\te'_{i+1}). 
    \label{eq:matching_edges}
\end{equation}

The map $\chi:\tPi\to \Prym_c(\tGa/\Ga)$ is the map $\tPi\to \Jac(\tGa)$ (the restriction of the second power of the Abel--Jacobi map on $\tGa$ to $\tPi\subset \Div_2^+(\tGa)$) followed by $\Id-\iota_\ast$. Hence $\nu_\#(\tbe)$ simply reinterprets $\tbe\in H_1(\tPi,\ZZ)$ as a $1$-chain in $H_1(\tGa,\ZZ)$ and applies $\Id-\iota_\ast$ to it. To finish the proof of claim (1) we write $\tbe = \tep - \iota_\ast(\tep)$ for $\tep = \sum_{i = 1}^{n} \left[\te_i - \iota (\te'_i)\right]$ and show that $\tep \in H_1(\tilde  \Gamma, \ZZ)$; this will imply that $\nu_\#(\tbe)$ is divisible by two:
\[
\nu_\#(\tbe)=(\Id-\iota_\ast)(\tbe)=(\Id-\iota_\ast)^2(\tep)=2(\Id-\iota_\ast)(\tep).
\]
Indeed, from Equation~\eqref{eq:matching_edges} and its $\iota$-pushforward it is already clear that the boundary is
\[ \partial \tep = s(\iota(\te'_1))-s(\te_1)+t(\te_n)-t(\iota(\te'_n)),
 \]
and we verify that the right hand side is zero. At the leftmost point of $\tbe$ (see Figure~\ref{fig:edgelabels_in_beta}) we obtain the condition
\[ \tu=s(\te_1)+s(\te'_1)=s(\iota(\te_1))+s(\iota(\te'_1)), \]
which implies either
\begin{align*} 
    s(\te_1) &= s(\iota (\te_1)) \qquad \text{and} \qquad  s(\te'_1) = s(\iota(\te'_1)) \\
    \text{or}\qquad  s(\te_1) &= s(\iota(\te'_1)) \qquad \text{and} \qquad s(\te'_1) = s(\iota (\te_1)).
\end{align*}
If we are in the first case, then in fact $s(\te_1)= s(\iota (\te_1))= s(\te'_1) = s(\iota(\te'_1))$ or else the tower $\tilde\Gamma \to \Gamma \to K$ would have a type V point. Hence we see that always $s(\iota(\te'_1)) - s(\te_1) = 0$. An analogous computation at the rightmost vertex $\tv$ shows that $t(\te_n)-t(\iota(\te'_n))=0$, and therefore $\partial \tep = 0$ and part (1) is established.
\medskip

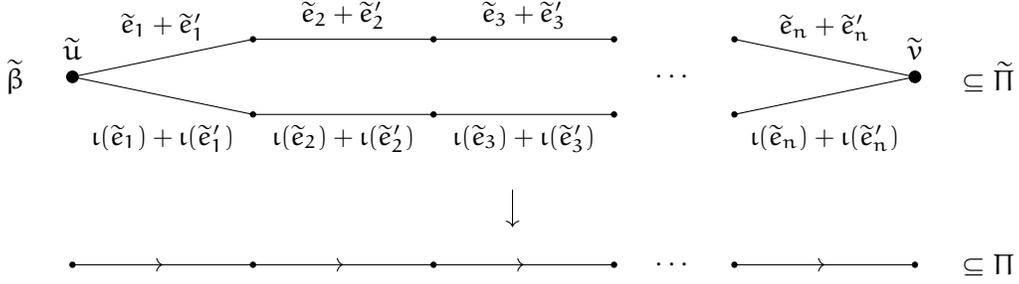
\begin{figure}
    \centering
    \begin{tikzpicture}
        \draw (-0.5, 0) node[anchor = east] {$\tbe$};
        \draw (11.7, 0) node[anchor = west] {$\subseteq \tilde\Pi$};
        \draw (11.7, -2.5) node[anchor = west] {$\subseteq \Pi$};
        \draw (0, 0.1) node[anchor = south] {$\tilde u$};
        \vertex[2]{0,0}

        \draw (0, 0) -- (2.4, 0.5) node [midway, above = 0.1] {\small $\te_1 + \te'_1$};
        \draw (0, 0) -- (2.4, -0.5) node [midway, below = 0.25] {\small $\iota (\te_1) + \iota (\te'_1)$};
        \vertex{2.4, 0.5}
        \vertex{2.4, -0.5}

        \draw (2.4, 0.5) -- (4.8, 0.5) node [midway, above] {\small $\te_2 + \te'_2$};
        \draw (2.4, -0.5) -- (4.8, -0.5) node [midway, below] {\small $\iota (\te_2) + \iota (\te'_2)$};
        
        \vertex{4.8, 0.5}
        \vertex{4.8, -0.5}

        \draw (4.8, 0.5) -- (7.2, 0.5) node [midway, above] {\small $\te_3 + \te'_3$};
        \draw (4.8, -0.5) -- (7.2, -0.5) node [midway, below] {\small $\iota (\te_3) + \iota (\te'_3)$};
        
        \vertex{7.2, 0.5}
        \vertex{7.2, -0.5}

        \draw (8, 0) node[anchor = center] {$\cdots$};
        
        \vertex{8.8, 0.5}
        \vertex{8.8, -0.5}

        \draw (8.8, 0.5) -- (11.2, 0) node [midway, above = 0.1] {\small $\te_n + \te'_n$};
        \draw (8.8, -0.5) -- (11.2, 0) node [midway, below = 0.25] {\small $\iota (\te_n) + \iota (\te'_n)$};

        \draw (11.2, 0.1) node[anchor = south] {$\tilde v$};
        \vertex[2]{11.2, 0}

        \path[draw, ->] (5.85, -1.5) -- (5.85, -2);
        \begin{scope}[decoration={
    			markings,
    			mark= at position 0.5 with {\arrow{>}}}
    		] 
    		\path[draw, postaction = decorate] (0, -2.5) -- (2.4, -2.5);
            \path[draw, postaction = decorate] (2.4, -2.5) -- (4.8, -2.5);
            \path[draw, postaction = decorate] (4.8, -2.5) -- (7.2, -2.5);
            \path[draw, postaction = decorate] (8.8, -2.5) -- (11.2, -2.5);
    	\end{scope}	        
        \draw (8, -2.5) node[anchor = center] {$\cdots$};
        \foreach \i in {0, 2.4, 4.8, 7.2, 8.8, 11.2} {
            \vertex{\i, -2.5}
        }
        
    \end{tikzpicture}
    \caption{Labeling convention for the edges in a cycle of type $\tbe$. The orientation for the edges in $\tilde\Pi$ is induced from the indicated orientation of the edges in $\Pi$.}
    \label{fig:edgelabels_in_beta}
\end{figure}

\iffalse

Now unfolding the construction of $\nu$ we see that $\nu_\#$ really just reinterprets the simplicial 1-chain $\beta \in C_1(\tilde\Pi, \ZZ)$ as a 1-chain in $C_1(\tilde\Gamma, \ZZ)$ in the obvious way and then applies $\Id - \iota$ to it. To finish the proof of claim (1) we write $\beta = \gamma - \iota_\ast(\gamma)$ for $\gamma = \sum_{i = 1}^{n} e_1^{(i)} - \iota e_2^{(i)}$ and show that $\gamma \in H_1(\tilde  \Gamma, \ZZ)$. Indeed, from Equation~\eqref{eq:matching_edges} and its $\iota$ pushforward it is already clear that
\[ \partial \gamma = s\iota e_2^{(1)} - se_1^{(1)} + te_1^{(n)} - t\iota e_2^{(n)} \]
and we verify that the right hand side is zero. From the shape of $\beta$ (see Figure~\ref{fig:edgelabels_in_beta}) we obtain the condition
\[ se_1^{(1)} + se_2^{(1)} = s\iota e_1^{(1)} + s\iota e_2^{(1)} \]
which translates into 
\begin{align*} 
    se_1^{(1)} &= s\iota e_1^{(1)} \qquad \text{and} \qquad  se_2^{(1)} = s \iota e_2^{(1)} \\
    \text{or}\qquad  s e_1^{(1)} &= s \iota e_2^{(1)} \qquad \text{and} \qquad s e_2^{(1)} = s \iota e_1^{(1)}.
\end{align*}
If we are in the first case, then in fact $se_1^{(1)} = s\iota e_1^{(1)} =  se_2^{(1)} = s \iota e_2^{(1)}$ or else the tower $\tilde\Gamma \to \Gamma \to K$ would have a type V point. Hence we see that always $s \iota e_2^{(1)} - se_1^{(1)} = 0$. The analogous consideration for the end vertex of $e_1^{(n)} + e_2^{(n)}$ finishes the proof that $\partial \gamma = 0$ hence part (1) is established.

\fi

To prove claim (2), we need to show that for any generator of $\Lambda_2$ of the form $d\tbe$ we have that $\nu^\#(d\tbe)$ is divisible by 2. Since $\nu^\#$ is essentially nothing but the descent of $\mu^\# : \Omega_{\tilde\Gamma}^1(\tGa) \to \Omega_{\tilde\Pi}^1(\tPi)$ to cokernels, it suffices to show that $\mu^\#(d\tbe) = d\tep - \iota^\ast d\tep$ for some $d\tep \in \Omega^1_{\tilde\Pi}(\tPi)$. For a description of $\mu^\#$ refer to the proof of Lemma~\ref{lem:alpha}. Now using the self-duality of the bigonal construction and the natural principal polarization of Jacobians, this computation is reduced to the one we already carried out in part (1) of this proof. Indeed, let $\tilde\Delta \to \Delta \to K$ be the bigonal construction of $\tilde\Pi \to \Pi \to K$, which we already know to be isomorphic to $\tilde\Gamma \to \Gamma \to K$. Then it is easy to see that
\begin{equation*}
    \begin{tikzcd}
        \Omega_{\tilde\Delta}^1(\tilde \Delta) \cong \Omega_{\tilde\Gamma}^1(\tGa) \arrow[r, "\mu^\#"] & \Omega_{\tilde \Pi}^1(\tPi)   \\
        H_1(\tilde\Delta, \ZZ) \arrow[u, "\cong"] \arrow[r, "\mu_\#"] & H_1(\tilde\Pi, \ZZ) \arrow[u, "\cong"]
    \end{tikzcd}
\end{equation*}
commutes, where the vertical maps are the principal polarizations of $\Jac(\tilde\Pi)$ and $\Jac(\tilde\Delta)$ and the map $\mu_\#$ is precisely the map studied in the proof of part (1) of this lemma with $\tilde\Delta$ in the role of $\tilde\Gamma$. \qedhere

\end{proof}

    \begin{lemma}  In the notation of Theorem~\ref{thm:bigonal_construction_restated}, we have
    \label{lem:pushbigonal}
    \[
    \nu_\ast ((\zeta'_c)^{h-1}) = 2\zeta_c^{h-1}\in H_{1,1}\big(\Prym_c(\tGa/\Ga)).
    \]
    \end{lemma}
    \begin{proof} Theorem~\ref{thm:tropical_Welters_criterion} states that
    \[ (\psi_{q'})_*\cyc [\tPi]=\frac{2}{(h-1)!} (\zeta_c')^{h-1}\in H_{1,1}\big(\Prym_c(\tilde\Pi / \Pi)\big). \]
    To prove the lemma, it suffices to show that 
    \begin{equation} \label{eq:step1_bigonal}
        \chi_*\cyc [\tPi]=\frac{4}{(h-1)!} \zeta_c^{h-1}\in H_{1,1}\big(\Prym_c(\tilde\Gamma / \Gamma)\big).
    \end{equation}
    Let $k = f \circ \pi : \tilde \Gamma \to K$. 
    Just as we did in the proof of Theorem~\ref{thm:tropical_Recillas_theorem_restated}, we define a tropical 1-cycle $B \in Z_1(\tilde\Gamma^2)$ which represents the lift of $\tilde \Pi \subseteq \Div_2^+(\tilde\Gamma)$ to $\tilde\Gamma^2$. We claim that 
    \[ B'(x, y) = \begin{cases}
        2 & \text{if } x = y \text{ and } d_k(x) = 4 \\
        1 & \text{if } k(x) = k(y), \ x \neq y, \text{ and } x \neq \iota(y) \\
        1 & \text{if } k(x) = k(y), \ \iota(x) \neq x, \text{ and } d_k(x) = 2 \\
        0 & \text{else}        
    \end{cases} \]
    satisfies balancing, i.e. $B'$ defines a tropical 1-cycle. As before, we refrain from proving this and instead we show that there is a tropical 1-cycle $B$ which agrees with $B'$ away from a 0-dimensional locus such that the following key formula 
    \begin{equation} \label{eq:proof_bigonal}
        [\Delta_{\tilde\Gamma}] + (\Id \times \iota)_*[\Delta_{\tilde\Gamma}] + B = (k\times k)^*[\Delta_K] = 4 [\tilde \Gamma \times p'] + 4 [p' \times \tilde \Gamma]
    \end{equation}
    holds in $H_{1,1}(\tilde\Gamma^2)$ for an arbitrary point $p' \in \tilde\Gamma$. 
    The second equality is the $(k \times k)$-pullback of Lemma~\ref{prop:diagonal_formula}. For the first equality we compute $(k \times k)^\ast [\Delta_K]$ as in the proof of Lemma~\ref{lem:key_formula} and we see that the result for the value of $(k\times k)^\ast [\Delta_K]$ at a point $(x, y)$ with $k(x) = k(y)$ is
    \begin{center}
        \begin{tabular}{c|c}
            type of $k(x)$ & $(k\times k)^\ast [\Delta_K](x,y)$ \\ \hline
            I & 4 \\
            II & 2 \\
            III & 1, unless $x = y$ with $d_k(x) = 2$, in which case it is 2 \\
            IV & 1 \\
        \end{tabular}
    \end{center}
    and the first equality of Equation~\eqref{eq:proof_bigonal} follows. Next we apply $(\psi_M^2)_*$ to Equation~\eqref{eq:proof_bigonal} to obtain an expression in $H_{1,1}\big(\Prym(\tilde\Gamma / \Gamma)\big)$ consisting of the following terms. 
    \begin{enumerate}
        \item $(\psi_M^2)_*[\Delta_{\tilde \Gamma}] = 4(\psi_M)_\ast[\tilde \Gamma]$ because again the pushforward along $\psi_M^2$ acts as $\psi_M$ followed by multiplication by 2 on the Prym variety which then induces multiplication by 4 on homology. 
        
        \item $(\psi_M^2)_*(\Id \times \iota)_*[\Delta_{\tilde\Gamma}] = 0$. This is because
        \[ \psi_M^2\big(p, \iota(p)\big) = \big(p-M-\iota(p)+\iota(M)\big)+\big(\iota(p)-M-p+\iota(M)\big) = -2M + 2\iota(M) \]
        is constant.
        
        \item $(\psi_M^2)_* B = 2\chi_\ast[\tilde \Pi]$, because at interior points of top-dimensional faces of $|B|$, the pushforward may be computed by first pushing forward along the (topological) quotient map identifying $(x,y)$ and $(y,x)$ and then pushing forward along $\chi$. A case-by-case analysis shows that the first step always introduces a factor of 2, either through the quotient map being 2:1 or because it is the value of $B$.
        
        \item $(\psi_M^2)_*[\tilde\Gamma \times p'] = (\psi_M^2)_*[p' \times \tilde\Gamma] = (\psi_M)_*[\tilde\Gamma]$.
    \end{enumerate}
    Solving the $(\psi_M^2)$-pushforward of Equation~\eqref{eq:proof_bigonal} for $\chi_*[\tilde\Pi]$ gives
    \[ \chi_\ast [\tPi] = 2 (\psi_M)_\ast [\tGa] = \frac{4}{(h-1)!} \zeta_c^{h-1}. \]
    The second equality is once more Theorem~\ref{thm:tropical_Welters_criterion} and this concludes the proof of Equation~\eqref{eq:step1_bigonal}. \qedhere

    \end{proof}

\FloatBarrier
\bibliographystyle{alpha}
\bibliography{references}

\end{document}